\documentclass[11pt]{article}
\usepackage{amscd}
\usepackage{amssymb,amsfonts,amsmath,psfrag,amscd,cite}
\usepackage{amsmath}
  \usepackage{paralist}
  \usepackage{graphics}
  \usepackage{epsfig}
  \usepackage[colorlinks=true]{hyperref}
  \usepackage[all]{xy}
\usepackage{graphicx}
\newtheorem{theo}{Theorem}
\newtheorem{prop}[theo]{Proposition}

\newtheorem{defi}[theo]{Definition}
\newtheorem{lemm}[theo]{Lemma}
\newtheorem{rema}[theo]{Remark}

\makeatletter \@addtoreset{equation}{section}
\@addtoreset{theo}{section} \makeatother

\begin{document}
\date{}
\title{Hamilton-Jacobi Equations for Nonholonomic Reducible\\ Hamiltonian Systems on a Cotangent Bundle}
\author{Manuel de Le\'{o}n \\
Instituto de Ciencias Matem\'{a}ticas ,\\
Consejo Superior de Investigaciones Cient\'{\i}ficas ,\\
c/Nicol\'{a}s Cabrera 13-15, 28049 Madrid, Spain.\\
E-mail: mdeleon@icmat.es \\\\
Hong Wang \thanks { Corresponding author Hong Wang (E-mail: hongwang@math.nankai.edu.cn).} \\
School of Mathematical Sciences and LPMC,\\
Nankai University,  Tianjin 300071, P.R.China \\
E-mail: hongwang@nankai.edu.cn\\\\
June 16, 2021} \maketitle

{\bf Abstract:} In this paper, for a variety of nonholonomic (reducible)
Hamiltonian systems, we first give to various
distributional Hamiltonian systems, by analyzing carefully
the dynamics and structures of the nonholonomic Hamiltonian systems.
Secondly, we derive precisely
the geometric constraint conditions of the induced distributional two-form
for the nonholonomic dynamical vector field,
which are called the Type I and Type II of Hamilton-Jacobi equations.
Thirdly, we generalize the above results for the nonholonomic
reducible Hamiltonian systems with symmetries, as well as with momentum maps,
and prove two types of Hamilton-Jacobi theorems
for various nonholonomic reduced
distributional Hamiltonian systems.
Finally, as an application, we give two
examples to illustrate the theoretical results.
These researches reveal the deeply internal
relationships of the nonholonomic constraints,
the induced (resp. reduced) distributional two-forms
and the dynamical vector fields of the nonholonomic Hamiltonian system
and its various distributional Hamiltonian systems.\\



\tableofcontents

\section{Introduction}

It is well-known that Hamilton-Jacobi theory is an important research subject
in mathematics and analytical mechanics.
see Abraham and Marsden \cite{abma78}, Arnold
\cite{ar89} and Marsden and Ratiu \cite{mara99},
and the Hamilton-Jacobi equation is
also fundamental in the study of the quantum-classical relationship
in quantization, and it also plays an important role
in the study of stochastic dynamical systems, see
Woodhouse \cite{wo92}, Ge and Marsden \cite{gema88},
and L\'{a}zaro-Cam\'{i} and Ortega \cite{laor09}.
Hamilton-Jacobi theory from the variational
point of view is originally developed by Jacobi in 1866, which states
that the integral of Lagrangian of a system along the solution of
its Euler-Lagrange equation satisfies the Hamilton-Jacobi equation.
The classical description of this problem from the generating function and the geometrical point
of view is given by Abraham and Marsden in \cite{abma78} as follows:
Let $Q$ be a smooth manifold and $TQ$ the tangent bundle, $T^* Q$
the cotangent bundle with the canonical symplectic form $\omega$,
and the projection $\pi_Q: T^* Q \rightarrow Q $ induces the map $
T\pi_{Q}: TT^* Q \rightarrow TQ. $
\begin{theo}
Assume that the triple $(T^*Q,\omega,H)$ is a Hamiltonian system
with Hamiltonian vector field $X_H$, and $W: Q\rightarrow
\mathbb{R}$ is a given generating function. Then the following two assertions
are equivalent:\\
\noindent $(\mathrm{i})$ For every curve $\sigma: \mathbb{R}
\rightarrow Q $ satisfying $\dot{\sigma}(t)= T\pi_Q
(X_H(\mathbf{d}W(\sigma(t))))$, $\forall t\in \mathbb{R}$, then
$\mathbf{d}W \cdot \sigma $ is an integral curve of the Hamiltonian
vector field $X_H$.\\
\noindent $(\mathrm{ii})$ $W$ satisfies the Hamilton-Jacobi equation
$H(q^i,\frac{\partial W}{\partial q^i})=E, $ where $E$ is a
constant.
\end{theo}

From the proof of the above theorem given in
Abraham and Marsden \cite{abma78}, we know that
the assertion $(\mathrm{i})$ with equivalent to
Hamilton-Jacobi equation $(\mathrm{ii})$ by the generating function,
gives a geometric constraint condition of the canonical symplectic form
on the cotangent bundle $T^*Q$
for Hamiltonian vector field of the system.
Thus, the Hamilton-Jacobi equation reveals the deeply internal relationships of
the generating function, the canonical symplectic form
and the dynamical vector field of a Hamiltonian system.\\

On the other hand, we have known that, in mechanics,
it is very often that many systems have constraints,
and usually, the constraints in dynamics are restrictions on
positions and velocities of the system.
There are two types of constraints, and the first one is holonomic,
which is that imposed on the configuration space of a system; and the second one
is nonholonomic, which involves the conditions on the velocities of a system,
such as rolling constraints. Thus, the nonholonomic mechanics describes
the motion of systems constrained by nonintegrable constraints,
i.e., constraints on the system velocities that do not arise from
constraints on the configurations alone. For a
nonholonomic Hamiltonian system, from Bates and $\acute{S}$niatycki \cite{basn93},
we know that, under the restriction given nonholonomic constraints, in general,
we can derive a distributional Hamiltonian system, which is called a semi-Hamiltonian
system in Patrick \cite{pa07}. But,
the leading distributional Hamiltonian system may not be a
Hamiltonian system, and it has no generating function,
then we cannot give the Hamilton-Jacobi theorem for the distributional
Hamiltonian system just like same as the above Theorem 1.1. Moreover,
for the nonholonomic Hamiltonian system with symmetry,  the
leading reduced distributional Hamiltonian system by nonholonomic reduction
may not be yet a Hamiltonian system, and we cannot give the Hamilton-Jacobi theorem
for the reduced distributional Hamiltonian system as the above Theorem 1.1.
We have to look for a new way. \\

It is worthy of noting that the regular point
symplectic reduction for the Hamiltonian system with symmetry and
coadjoint equivariant momentum map was set up by famous professors
Jerrold E. Marsden and Alan Weinstein, which is called
Marsden-Weinstein reduction, and great developments have been
obtained around the work in the theoretical study and applications
of mathematics, mechanics and physics; see
Abraham and Marsden \cite{abma78}, Arnold \cite{ar89}, 
Koiller \cite{ko92}, Libermann and Marle
\cite{lima87}, Marsden \cite{ma92},
Marsden et al. \cite{mamiorpera07, mamora90}, 
Marsden and Perlmutter \cite{mape00},
Marsden and Ratiu \cite{mara99}, Marsden and
Weinstein \cite{mawe74}, Meyer \cite{me73},
Nijmeijer and Van der Schaft \cite {nivds90} and
Ortega and Ratiu \cite{orra04} for more details and development.
But, in Marsden et al.\cite{mawazh10} and Wang \cite{wa18},
the authors found that the symplectic reduced space of a Hamiltonian system
defined on the cotangent bundle of a configuration manifold may not be a
cotangent bundle, and hence the set of Hamiltonian systems with
symmetries on the cotangent bundle is not complete under the
Marsden-Weinstein reduction. Thus, the symplectic reduced system of a
Hamiltonian system with symmetry defined on the cotangent bundle
may not be a Hamiltonian system on a cotangent bundle,
then we cannot give the Hamilton-Jacobi theorem for the Marsden-Weinstein
reduced system as the above Theorem 1.1. \\

Now, it is a natural problem how to generalize Theorem 1.1 to fit
the above nonholonomic systems and their reduced systems.
Note that if take that $\gamma=\mathbf{d}W$ in
the above Theorem 1.1, then $\gamma$ is a closed one-form on $Q$, and
the equation $\mathbf{d}(H \cdot \mathbf{d}W)=0$ is equivalent to
the Hamilton-Jacobi equation $H(q^i,\frac{\partial W}{\partial
q^i})=E$, where $E$ is a constant, which is called the classical
Hamilton-Jacobi equation. This result is used the
formulation of a geometric version of Hamilton-Jacobi theorem for
Hamiltonian system, see Cari\~{n}ena et al \cite{cagrmamamuro06,
cagrmamamuro10}. Moreover, note that Theorem 1.1 is also
generalized in the context of time-dependent Hamiltonian system by
Marsden and Ratiu in \cite{mara99}, and the Hamilton-Jacobi equation may
be regarded as a nonlinear partial differential equation for some
generating function $S$. Thus, the problem is become how to choose a
time-dependent canonical transformation $\Psi: T^*Q\times \mathbb{R}
\rightarrow T^*Q\times \mathbb{R}, $ which transforms the dynamical
vector field of a time-dependent Hamiltonian system to equilibrium,
such that the generating function $S$ of $\Psi$ satisfies the
time-dependent Hamilton-Jacobi equation. In particular, for the
time-independent Hamiltonian system, ones may look for a symplectic
map as the canonical transformation. This work offers an important
idea that one can use the dynamical vector field of a Hamiltonian
system to describe Hamilton-Jacobi equation. In consequence, if assume that
$\gamma: Q \rightarrow T^*Q$ is a closed one-form on $Q$, and define
that $X_H^\gamma = T\pi_{Q}\cdot X_H \cdot \gamma$, where $X_{H}$ is
the dynamical vector field of Hamiltonian system $(T^*Q,\omega,H)$,
then the fact that $X_H^\gamma$ and $X_H$ are $\gamma$-related, that
is, $T\gamma\cdot X_H^\gamma= X_H\cdot \gamma$ is equivalent that
$\mathbf{d}(H \cdot \gamma)=0, $ which is given in Cari\~{n}ena et
al \cite{cagrmamamuro06, cagrmamamuro10}.
Motivated by the above research work, Wang in \cite{wa17} prove
an important lemma, which is
a modification for the corresponding result of Abraham and Marsden
in \cite{abma78}, such that we can derive precisely
the geometric constraint conditions of
the regular reduced symplectic forms for the
dynamical vector fields of a regular reducible Hamiltonian
system on the cotangent bundle of a configuration manifold,
which are called the Type I and Type II of Hamilton-Jacobi equation,
because they are the development of the above classical Hamilton-Jacobi equation
given by Theorem 1.1, see Abraham and Marsden \cite{abma78}
and Wang \cite{wa17}. \\

Since the Hamilton-Jacobi theory is developed based on the
Hamiltonian picture of dynamics, it is natural idea to extend the
Hamilton-Jacobi theory to the nonholonomic Hamiltonian systems,
and they with symmetry and momentum map. Our idea is that
to use a variety of dynamical vector fields of the (reduced) distributional
Hamiltonian systems to describe a variety of Hamilton-Jacobi equations.
In this paper, we first derive precisely
the geometric constraint conditions of the induced distributional two-form
and the reduced distributional two-form
for the nonholonomic dynamical vector fields.
Moreover, we consider
the nonholonomic reductions compatible with Marsden-Weinstein reduction
and regular orbit reduction, and lead to the $\mathbf{J}$-nonholonomic
regular point and orbit reduced distributional Hamiltonian systems,
by analyzing carefully the dynamics and structures of the nonholonomic
Hamiltonian systems. These systems are not yet Hamiltonian, but, we can
give their two types of Hamilton-Jacobi equations, as an
extension of two types of Hamilton-Jacobi equations
for the Marsden-Weinstein reduced Hamiltonian system
and the regular orbit reduced Hamiltonian system
given in \cite{wa17} to the nonholonomic context.\\

The paper is organized as follows. In section 2 we first recall the
main facts about the dynamics of a nonholonomic Hamiltonian system,
including the influence of symmetries, which are helpful for us to
understand the constructions of a distributional Hamiltonian system
and its a variety of the nonholonomic reduced distributional Hamiltonian systems.
In section 3, we first prove an important lemma, which is a tool for our research.
Then derive precisely
the geometric constraint conditions of the distributional two-form
for the nonholonomic dynamical vector field,
that is, the two types of Hamilton-Jacobi equation
for the distributional Hamiltonian system.
The nonholonomic reducible Hamiltonian systems with
symmetries, as well as momentum maps, are considered respectively in
section 4 and section 5, and derive precisely
the geometric constraint conditions of a variety of nonholonomic reduced
distributional two-forms
for the nonholonomic reducible Hamiltonian vector fields,
that is, the two types of Hamilton-Jacobi
equations for a variety of nonholonomic reduced distributional Hamiltonian systems
(in particular, when the Lie group is not Abelian). As the
applications of the theoretical results, we consider the motions of
the constrained particle in space $\mathbb{R}^3$ and the vertical
rolling disk in section 6, and derive two types of Hamilton-Jacobi
equations for the distributional Hamiltonian systems and their
reduced distributional Hamiltonian systems corresponding to the two
nonholonomic systems. These research work develop the nonholonomic
reduction and Hamilton-Jacobi theory of the nonholonomic Hamiltonian
systems with symmetries, as well as momentum maps,
and make us have much deeper understanding
and recognition for the structures of the nonholonomic Hamiltonian
systems.

\section{Dynamics of Nonholonomic Mechanical System}

In this section, we first review briefly some basic facts about
nonholonomic mechanical systems and give the descriptions of
dynamics of a nonholonomic Hamiltonian system and the nonholonomic
Hamiltonian system with symmetry, as well as momentum maps, which are helpful for us in
subsequent sections to understand the constructions of
distributional Hamiltonian system and the nonholonomic reduced
distributional Hamiltonian system. We shall follow the notations and
conventions introduced in Cantrijn et al. \cite{calemama99}, Bates
and $\acute{S}$niatycki in \cite{basn93}, Cushman et al.
\cite{cudusn10} and \cite{cukesnba95}, Montgomery \cite{mo02},
de Le\'{o}n et al. \cite{lema96}, Marsden et al. \cite{mawazh10} and Wang \cite{wa17}.\\

In order to describe the dynamics of a nonholonomic mechanical system,
we need some restriction conditions for nonholonomic constraints of
the system. At first, we note that the set of Hamiltonian vector fields
forms a Lie algebra with respect to the Lie bracket, since
$X_{\{f,g\}}=-[X_f, X_g]. $ But, the Lie bracket operator, in
general case, may not be closed on the restriction of a nonholonomic
constraint. Thus, we have to give the following completeness
condition for nonholonomic constraints of a system.\\

{\bf $\mathcal{D}$-completeness } Let $Q$ be a smooth manifold and
$TQ$ its tangent bundle. A distribution $\mathcal{D} \subset TQ$ is
said to be {\bf completely nonholonomic} (or bracket-generating) if
$\mathcal{D}$ along with all of its iterated Lie brackets
$[\mathcal{D},\mathcal{D}], [\mathcal{D}, [\mathcal{D},\mathcal{D}]],
\cdots ,$ spans the tangent bundle $TQ$. Moreover, we consider a
mechanical system on $Q$. Then nonholonomic constraints of the
system are said to be {\bf completely nonholonomic} if the
distribution $\mathcal{D} \subset TQ$ defined by the nonholonomic
constraints is completely nonholonomic.\\

In this paper we consider that a nonholonomic mechanical system is
given by a Lagrangian function $L: TQ \rightarrow \mathbb{R}$
subject to constraints determined by a completely nonholonomic
distribution $\mathcal{D}\subset TQ$ on the configuration manifold
$Q$. We denote by $D$ the total space of $\mathcal{D}$ in $TQ$,
which is a constraint submanifold. For simplicity we always assume
that $\tau_Q(D)=Q, $ where $\tau_Q: TQ \rightarrow Q$ is the
canonical projection, that is, the constraints are purely
kinematical in the sense that they do not impose restrictions on the
allowable positions. The motions of the system are forced to take
place on $D$ and this requires the introduction of some "reaction
force". In order to describe the constraint submanifold in the phase
space and the dynamics of system, we have to give the following
regularity condition.\\

{\bf $\mathcal{D}$-regularity } In the following we always assume
that $Q$ is a smooth manifold with coordinates $(q^i)$, and $TQ$ its
tangent bundle with coordinates $(q^i,\dot{q}^i)$, and $T^\ast Q$
its cotangent bundle with coordinates $(q^i,p_j)$, which are the
canonical cotangent coordinates of $T^\ast Q$ and $\omega=
dq^{i}\wedge dp_{i}$ is canonical symplectic form on $T^{\ast}Q$. If
the Lagrangian $L: TQ \rightarrow \mathbb{R}$ is hyperregular, that
is, the Hessian matrix
$(\partial^2L/\partial\dot{q}^i\partial\dot{q}^j)$ is nondegenerate
everywhere, then the Legendre transformation $FL: TQ \rightarrow T^*
Q$ is a diffeomorphism. In this case the Hamiltonian $H: T^* Q
\rightarrow \mathbb{R}$ is given by $H(q,p)=\dot{q}\cdot
p-L(q,\dot{q}) $ with Hamiltonian vector field $X_H$,
which is defined by the Hamilton's equation
$\mathbf{i}_{X_H}\omega=\mathbf{d}H$, and
$\mathcal{M}=\mathcal{F}L(\mathcal{D})$ is a constraint submanifold
in $T^* Q$. In particular, for the nonholonomic constraint
$\mathcal{D}$, the Lagrangian $L$ is said to be {\bf
$\mathcal{D}$-regular}, if the restriction of Hessian matrix
$(\partial^2L/\partial\dot{q}^i\partial\dot{q}^j)$ on $\mathcal{D}$
is nondegenerate everywhere. Moreover, a nonholonomic system is said
to be {\bf $\mathcal{D}$-regular}, if its Lagrangian $L$ is {\bf
$\mathcal{D}$-regular}. Note that the restriction of a positive
definite symmetric bilinear form to a subspace is also positive
definite, and hence nondegenerate. Thus, for a simple nonholonomic
mechanical system, that is, whose Lagrangian is the total kinetic
energy minus potential energy, it is {\bf $\mathcal{D}$-regular }
automatically, which is coincident with the sense of regularity of
nonholonomic system given by de L\'{e}on and
Mart\'{\i}n de Diego \cite{lema96}.\\

A nonholonomic Hamiltonian system is a 4-tuple
$(T^\ast Q,\omega,\mathcal{D},H)$, which is a Hamiltonian system with a
$\mathcal{D}$-completely and $\mathcal{D}$-regularly nonholonomic
constraint $\mathcal{D} \subset TQ$.
In the following we shall describe the dynamics of the nonholonomic
Hamiltonian system $(T^*Q,\omega,\mathcal{D},H)$. We define the
distribution $\mathcal{F}$ as the pre-image of the nonholonomic
constraints $\mathcal{D}$ for the map $T\pi_Q: TT^* Q \rightarrow
TQ$, that is, $\mathcal{F}=(T\pi_Q)^{-1}(\mathcal{D})\subset TT^*Q,
$ which is a distribution along $\mathcal{M}$, and
$\mathcal{F}^\circ:=\{\alpha \in T^*T^*Q | <\alpha,v>=0, \; \forall
v\in TT^*Q \}$ is the annihilator of $\mathcal{F}$ in
$T^*T^*Q_{|\mathcal{M}}$. We consider the following nonholonomic
constraints condition
\begin{align} (\mathbf{i}_X \omega -\mathbf{d}H) \in \mathcal{F}^\circ,
\;\;\;\;\;\; X \in T\mathcal{M},
\label{2.1} \end{align} from Cantrijn et al.
\cite{calemama99}, we know that there exists an unique nonholonomic
vector field $X_n$ satisfying the above condition $(2.1)$, if the
admissibility condition $\mathrm{dim}\mathcal{M}=
\mathrm{rank}\mathcal{F}$ and the compatibility condition
$T\mathcal{M}\cap \mathcal{F}^\bot= \{0\}$ hold, where
$\mathcal{F}^\bot$ denotes the symplectic orthogonal of
$\mathcal{F}$ with respect to the canonical symplectic form
$\omega$ on $T^*Q$. In particular, when we consider the Whitney sum
decomposition $T(T^*Q)_{|\mathcal{M}}=T\mathcal{M}\oplus
\mathcal{F}^\bot$ and the canonical projection $P:
T(T^*Q)_{|\mathcal{M}} \rightarrow T\mathcal{M}$,
we have that $X_n= P(X_H)$.\\

If the Lagrangian $L: TQ \rightarrow \mathbb{R}$ is singular, in
this case the Hessian matrix
$(\partial^2L/\partial\dot{q}^i\partial\dot{q}^j)$ is degenerate. By
using the Gotay-Nester presymplectic constraint algorithm, see
\cite{gone79}, we can find a final constraint submanifold
$\mathcal{M}_f \subset T^*Q$, such that on which there exists a
nonholonomic vector field $X_n$ satisfying the following
nonholonomic constraints condition
\begin{align} (\mathbf{i}_X \omega -\mathbf{d}H)_{|\mathcal{M}_f} \in \mathcal{F}^\circ,
\;\;\;\;\;\; X_{|\mathcal{M}_f} \in T\mathcal{M}_f.
\label{2.2} \end{align}
Therefore, without loss of generality, we shall henceforth always
assume that there exists a nonholonomic vector field $X_n$
satisfying the nonholonomic constraints condition.\\

From the condition (2.1) we know that the nonholonomic vector field,
in general case, may not be Hamiltonian, because of the restriction
of nonholonomic constraints. But, we hope to study the dynamical
vector field of nonholonomic Hamiltonian system by using the similar
method of studying Hamiltonian vector field. On the other hand, we
also note that Bates and $\acute{S}$niatycki in \cite{basn93} give a
method to study the nonholonomic Hamiltonian system and nonholonomic
reduction. In fact, for a nonholonomic Hamiltonian system
$(T^*Q,\omega,\mathcal{D},H)$, by using their method, we know that
there exist a distribution $\mathcal{K}=\mathcal {F}\cap
T\mathcal{M}$, a non-degenerate distributional two-form
$\omega_\mathcal{K}$, which is the
restriction of the induced symplectic form $\omega_{\mathcal{M}}$ on
$T^*\mathcal{M}$ fibrewise to the distribution $\mathcal{K}$,
and a vector field $X_\mathcal {K}$ on the
constraint submanifold $\mathcal{M}=\mathcal{F}L(\mathcal{D})\subset
T^*Q$, such that the distributional Hamiltonian equation
$\mathbf{i}_{X_\mathcal{K}}\omega_\mathcal{K}=\mathbf{d}H_\mathcal
{K}$ holds, where the function $H_{\mathcal{K}}$ satisfies
$\mathbf{d}H_{\mathcal{K}}= \tau_{\mathcal{K}}\cdot \mathbf{d}H_{\mathcal {M}}$,
and $H_{\mathcal {M}}$ is the
restriction of the Hamiltonian function $H$ to the constraint submanifold $\mathcal{M}$,
and $\tau_{\mathcal{K}}$ is the restriction map to distribution $\mathcal{K}$.
Then the triple $(\mathcal{K},\omega_{\mathcal{K}},H_{\mathcal{K}})$
is called a distributional Hamiltonian system, and $X_\mathcal {K}$ is its
nonholonomic vector field.\\

Moreover, we consider the nonholonomic Hamiltonian system with
symmetry and nonholonomic reduction. Assume that Lie group $G$ acts
smoothly by the left on $Q$, its tangent lifted acts on $TQ$ and its
cotangent lifted acts on $T^\ast Q$, which is free, proper and
symplectic. The orbit space $T^* Q/ G$ is a smooth manifold and the
canonical projection $\pi_{/G}: T^* Q \rightarrow T^* Q /G $ is a
surjective submersion. In the following we shall describe the
dynamics of the nonholonomic Hamiltonian system with symmetry
$(T^*Q,G,\omega,\mathcal{D},H)$, where $H: T^*Q \rightarrow
\mathbb{R}$ is a $G$-invariant Hamiltonian, and the completely
nonholonomic constraints $\mathcal{D}\subset TQ$ is a $G$-invariant
distribution, that is, the tangent of the group action maps
$\mathcal{D}_q$ to $\mathcal{D}_{gq}$ for any $q\in Q $. Since the
Legendre transformation $\mathcal{F}L: TQ \rightarrow T^*Q$ is a
fiber-preserving map, then
$\mathcal{M}=\mathcal{F}L(\mathcal{D})\subset T^*Q$ is
$G$-invariant, and the quotient space
$\bar{\mathcal{M}}=\mathcal{M}/G$ of the $G$-orbit in $\mathcal{M}$
is a smooth manifold with projection $\pi_{/G}:
\mathcal{M}\rightarrow \bar{\mathcal{M}}( \subset T^* Q /G)$ which
is a surjective submersion. From Bates and $\acute{S}$niatycki
\cite{basn93}, we know that there exists a distribution
$\bar{\mathcal{K}}$, a non-degenerate distributional two-form
$\omega_{\bar{\mathcal{K}}}$,
and a vector field $X_{\bar{\mathcal {K}}}$
on $\bar{\mathcal{M}}$ which takes values in the constraint
distribution $\bar{\mathcal{K}}$, such that the following equation
holds, that is,
$\mathbf{i}_{X_{\bar{\mathcal{K}}}}\omega_{\bar{\mathcal{K}}}
=\mathbf{d}h_{\bar{\mathcal{K}}}$,
where the function $h_{\bar{\mathcal{K}}}$ satisfies
$\mathbf{d}h_{\bar{\mathcal{K}}}= \tau_{\bar{\mathcal{K}}}\cdot \mathbf{d}h_{\bar{\mathcal{M}}}$,
and $h_{\bar{\mathcal{M}}}\cdot \pi_{/G}= H_{\mathcal{M}}$. In this case, the triple
$(\bar{\mathcal{K}},\omega_{\bar{\mathcal {K}}},h_{\bar{\mathcal{K}}})$ is called a
nonholonomic reduced distributional Hamiltonian system, and
$X_{\bar{\mathcal {K}}}$ is its nonholonomic reduced dynamical vector field.\\

In particular, we assume that the Lie group $G$ is not Abelian,
and the cotangent lifted $G$-action on
$T^*Q$ is free, proper and symplectic, and admits a
$\operatorname{Ad}^\ast$-equivariant momentum map $\mathbf{J}:
T^\ast Q\rightarrow \mathfrak{g}^\ast$, where $\mathfrak{g}$ is a
Lie algebra of $G$ and $\mathfrak{g}^\ast$ is the dual of
$\mathfrak{g}$. Let $\mu \in\mathfrak{g}^\ast$ be a regular value of
$\mathbf{J}$ and denote by $G_\mu$ the isotropy subgroup of the
coadjoint $G$-action at the point $\mu\in\mathfrak{g}^\ast$, which
is defined by $G_\mu=\{g\in G|\operatorname{Ad}_g^\ast \mu=\mu \}$.
Since $G_\mu (\subset G)$ acts freely and properly on $Q$ and on
$T^\ast Q$, then $G_\mu$ acts also freely and properly on
$\mathbf{J}^{-1}(\mu)$, so that the space $(T^\ast
Q)_\mu=\mathbf{J}^{-1}(\mu)/G_\mu$ is a symplectic manifold with
symplectic form $\omega_\mu$ uniquely characterized by the relation
\begin{equation}\pi_\mu^\ast \omega_\mu=i_\mu^\ast
\omega. \label{2.3}\end{equation} The map
$i_\mu:\mathbf{J}^{-1}(\mu)\rightarrow T^\ast Q$ is the inclusion
and $\pi_\mu:\mathbf{J}^{-1}(\mu)\rightarrow (T^\ast Q)_\mu$ is the
projection. The pair $((T^\ast Q)_\mu,\omega_\mu)$ is
the Marsden-Weinstein reduced space of $(T^\ast Q,\omega)$ at $\mu$,
(see Marsden and Weinstein \cite{mawe74}, Marsden \cite{ma92},
and Marsden et al.\cite{mamiorpera07}). In
the following we assume that for the regular value
$\mu\in\mathfrak{g}^\ast$, the constraint submanifold $\mathcal{M}$
is clean intersection with $\mathbf{J}^{-1}(\mu)$, that is,
$\mathcal{M} \cap \mathbf{J}^{-1}(\mu)\neq \emptyset$. Note that
$\mathcal{M}$ is also $G_\mu (\subset G)$ action invariant, and so
is $\mathbf{J}^{-1}(\mu)$, because $\mathbf{J}$ is
$\operatorname{Ad}^\ast$-equivariant. It follows that the quotient
space $\mathcal{M}_\mu =(\mathcal{M}\cap \mathbf{J}^{-1}(\mu))
/G_\mu \subset (T^\ast Q)_\mu$ of the $G_\mu$-orbit in
$\mathcal{M}\cap \mathbf{J}^{-1}(\mu)$, is a smooth manifold with
projection $\pi_\mu: \mathcal{M}\cap \mathbf{J}^{-1}(\mu)
\rightarrow \mathcal{M}_\mu$ which is a surjective submersion.\\

In the following we shall describe the dynamics of the nonholonomic
Hamiltonian system with symmetry and momentum map
$(T^*Q,G,\omega,\mathbf{J},\mathcal{D},H)$ by using the method given
by Bates and $\acute{S}$niatycki in \cite{basn93}. Assume that the
distribution $T(\mathbf{J}^{-1}(\mu))\cap \mathcal{F}$ pushes down
to a distribution $\mathcal{F}_\mu= T\pi_\mu
(T(\mathbf{J}^{-1}(\mu))\cap \mathcal{F})$ on $(T^\ast Q)_\mu$ along
$\mathcal{M}_\mu$, and $h_\mu$ is the Marsden-Weinstein reduced Hamiltonian function
$h_\mu: (T^* Q)_\mu \rightarrow \mathbb{R}$ defined by $h_\mu\cdot
\pi_\mu= H\cdot i_\mu$. We consider the following nonholonomic
constraints condition
\begin{align} (\mathbf{i}_{X_\mu} \omega_\mu -\mathbf{d}h_\mu)_{|\mathcal{M}_\mu}
\in \mathcal{F}_\mu^\circ, \;\;\;\;\;\; X_\mu \in T\mathcal{M}_\mu.
\label{2.4} \end{align}
Thus, there exists an unique
nonholonomic vector field $X_\mu$ satisfying the above condition
$(2.4)$, if the admissibility condition
$\mathrm{dim}\mathcal{M}_\mu= \mathrm{rank}\mathcal{F}_\mu$ and the
compatibility condition $T\mathcal{M}_\mu\cap \mathcal{F}_\mu^\bot=
\{0\}$ hold, where $\mathcal{F}_\mu^\bot$ is denoted the symplectic
orthogonal of $\mathcal{F}_\mu$ with respect to the Marsden-Weinstein reduced
symplectic form $\omega_\mu$. In consequence, we know that there exists a
distribution $\mathcal{K}_\mu=\mathcal {F}_\mu\cap
T\mathcal{M}_\mu$, a non-degenerate reduced
distributional two-form $\omega_{\mathcal{K}_\mu},$
which is the restriction of the induced symplectic form $\omega_{\mathcal{M}_\mu}$ on
$T^*\mathcal{M}_\mu$ fibrewise to the distribution $\mathcal{K}_\mu$,
and a vector field $X_{\mathcal {K}_\mu}$ on the reduced constraint submanifold
$\mathcal{M}_\mu=(\mathcal{M}\cap \mathbf{J}^{-1}(\mu)) /G_\mu, $
such that the equation
$\mathbf{i}_{X_{\mathcal{K}_\mu}}\omega_{\mathcal{K}_\mu} =
\mathbf{d}h_{\mathcal{K}_\mu}$ holds, where the function $h_{\mathcal {K}_\mu}$ satisfies
$\mathbf{d}h_{\mathcal{K}_\mu}= \tau_{\mathcal{K}_\mu}\cdot \mathbf{d}h_\mu $,
that is, the restriction condition of the $\mathbf{d}h_\mu$
to the reduced distribution $\mathcal{K}_\mu$.
Then the triple
$(\mathcal{K}_\mu,\omega_{\mathcal {K}_\mu},h_{\mathcal {K}_\mu})$ is called a
$\mathbf{J}$-nonholonomic regular point reduced distributional Hamiltonian system,
and $X_{\mathcal {K}_\mu}$ is its $\mathbf{J}$-nonholonomic regular point reduced
dynamical vector field.\\

It is worthy of noting that the orbit reduction of a Hamiltonian system is an alternative
approach to symplectic reduction given by Marle \cite{ma76}
and Kazhdan, Kostant and Sternberg \cite{kakost78}, which is different from the
Marsden-Weinstein reduction. For the nonholonomic
Hamiltonian system with symmetry and momentum map
$(T^*Q,G,\omega,\mathbf{J},\mathcal{D},H)$,
if $\mu\in \mathfrak{g}^\ast$ is a
regular value of the momentum map $\mathbf{J}$ and
$\mathcal{O}_\mu=G\cdot \mu\subset \mathfrak{g}^\ast$ is the
$G$-orbit of the coadjoint $G$-action through the point $\mu$,
by using the above method, we know that there exists a
distribution $\mathcal{K}_{\mathcal{O}_\mu}$, a non-degenerate reduced
distributional two-form $\omega_{\mathcal{K}_{\mathcal{O}_\mu}}$ and a vector
field $X_{\mathcal {K}_{\mathcal{O}_\mu}}$ on the regular orbit reduced constraint submanifold
$\mathcal{M}_{\mathcal{O}_\mu}=(\mathcal{M}\cap \mathbf{J}^{-1}(\mathcal{O}_\mu)) /G, $
such that the equation
$\mathbf{i}_{X_{\mathcal{K}_{\mathcal{O}_\mu}}}\omega_{\mathcal{K}_{\mathcal{O}_\mu}} =
\mathbf{d}h_{\mathcal{K}_{\mathcal{O}_\mu}}$ holds.
Here the regular orbit reduced space is
$((T^\ast Q)_{\mathcal{O}_\mu}=\mathbf{J}^{-1}(\mathcal{O}_\mu)/G,
\omega_{\mathcal{O}_\mu}),$ in which the the symplectic form $\omega_{\mathcal{O}_\mu}$
uniquely characterized by the relation
\begin{equation}i_{\mathcal{O}_\mu}^\ast \omega=\pi_{\mathcal{O}_{\mu}}^\ast
\omega_{\mathcal{O}
_\mu}+\mathbf{J}_{\mathcal{O}_\mu}^\ast\omega_{\mathcal{O}_\mu}^+,
\label{2.5}
\end{equation}
where $\mathbf{J}_{\mathcal{O}_\mu}$ is
the restriction of the momentum map $\mathbf{J}$ to
$\mathbf{J}^{-1}(\mathcal{O}_\mu)$, that is,
$\mathbf{J}_{\mathcal{O}_\mu}=\mathbf{J}\cdot i_{\mathcal{O}_\mu}$
and $\omega_{\mathcal{O}_\mu}^+$ is the $+$-symplectic structure on
the orbit $\mathcal{O}_\mu$ given by
\begin{equation}\omega_{\mathcal{O}_\mu}^
+(\nu)(\xi_{\mathfrak{g}^\ast}(\nu),\eta_{\mathfrak{g}^\ast}(\nu))
=<\nu,[\xi,\eta]>,\;\; \forall\;\nu\in\mathcal{O}_\mu, \;
\xi,\eta\in \mathfrak{g}. \label{2.6}
\end{equation}
The maps
$i_{\mathcal{O}_\mu}:\mathbf{J}^{-1}(\mathcal{O}_\mu)\rightarrow
T^\ast Q$ and
$\pi_{\mathcal{O}_\mu}:\mathbf{J}^{-1}(\mathcal{O}_\mu)\rightarrow
(T^\ast Q)_{\mathcal{O}_\mu}$ are natural injection and the
projection, respectively. The
distribution $T(\mathbf{J}^{-1}(\mathcal{O}_\mu))\cap \mathcal{F}$ pushes down
to a distribution $\mathcal{F}_{\mathcal{O}_\mu}= T\pi_{\mathcal{O}_\mu}
(T(\mathbf{J}^{-1}(\mathcal{O}_\mu))\cap \mathcal{F})$ on $(T^\ast Q)_{\mathcal{O}_\mu}$ along
$\mathcal{M}_{\mathcal{O}_\mu}$, and
$\mathcal{K}_{\mathcal{O}_\mu}=\mathcal {F}_{\mathcal{O}_\mu}\cap T\mathcal{M}_{\mathcal{O}_\mu}.$
Then the triple
$(\mathcal{K}_{\mathcal{O}_\mu},\omega_{\mathcal {K}_{\mathcal{O}_\mu}},
h_{\mathcal{K}_{\mathcal{O}_\mu}})$ is called a
$\mathbf{J}$-nonholonomic regular orbit reduced distributional Hamiltonian system,
in which $\omega_{\mathcal {K}_{\mathcal{O}_\mu}}$ is the restriction of the induced symplectic form $\omega_{\mathcal{M}_{\mathcal{O}_\mu}}$ on
$T^*\mathcal{M}_{\mathcal{O}_\mu}$ fibrewise to the distribution
$\mathcal{K}_{\mathcal{O}_\mu}$, and the function
$h_{\mathcal {K}_{\mathcal{O}_\mu}}$ satisfies
$\mathbf{d}h_{\mathcal{K}_{\mathcal{O}_\mu}}= \tau_{\mathcal{K}_{\mathcal{O}_\mu}}\cdot \mathbf{d}h_{\mathcal{O}_\mu} $,
that is, the restriction condition of the $\mathbf{d}h_{\mathcal{O}_\mu}$
to the reduced distribution $\mathcal{K}_{\mathcal{O}_\mu}$,
where
$h_{\mathcal{O}_\mu}$ is regular orbit reduced Hamiltonian function
$h_{\mathcal{O}_\mu}: (T^* Q)_{\mathcal{O}_\mu} \rightarrow \mathbb{R}$ defined by
$h_{\mathcal{O}_\mu}\cdot \pi_{\mathcal{O}_\mu}= H\cdot i_{\mathcal{O}_\mu}$,
and $X_{\mathcal {K}_{\mathcal{O}_\mu}}$ is the $\mathbf{J}$-nonholonomic
regular orbit reduced dynamical vector field.\\

In the following we shall derive precisely
the geometric constraint conditions of the induced distributional two-form
and the reduced distributional two-forms
for a variety of nonholonomic dynamical vector fields,
that is, the two types of Hamilton-Jacobi equations
for the various distributional Hamiltonian systems.

\section{Hamilton-Jacobi Equations for a Distributional Hamiltonian System }

In this section, for a nonholonomic Hamiltonian system
$(T^*Q,\omega,\mathcal{D},H)$, where $\omega$ is the canonical
symplectic form on $T^* Q$, and $\mathcal{D}\subset TQ$ is a
$\mathcal{D}$-completely and $\mathcal{D}$-regularly nonholonomic
constraint of the system, we first give its distribution
$\mathcal{K}$, an associated non-degenerate distributional two-form
$\omega_\mathcal{K}$ induced by the canonical symplectic form and a
distributional Hamiltonian system, then derive precisely
the geometric constraint conditions of the distributional two-form
for the nonholonomic dynamical vector field,
that is, the two types of Hamilton-Jacobi equation
for the distributional Hamiltonian system. In order to do this, we need
first to analyze carefully the dynamics and structure of the nonholonomic
Hamiltonian system following the results given by Bates
and $\acute{S}$niatycki in \cite{basn93}, (see also Cushman et al.
\cite{cudusn10} and
\cite{cukesnba95} for more details).\\

From now on, we assume that $L:TQ \rightarrow \mathbb{R}$ is a
hyperregular Lagrangian, and the Legendre transformation
$\mathcal{F}L: TQ \rightarrow T^*Q$ is a diffeomorphism. As above,
our nonholonomic constraint $\mathcal{D} \subset TQ$ is $\mathcal{D}$-completely
and $\mathcal{D}$-regularly, and let $\mathcal{D}^0 \subset T^*Q$ its annihilator. From
$\S 2$, we can define the constraint submanifold
$\mathcal{M}=\mathcal{F}L(\mathcal{D})\subset T^*Q$,
$i_{\mathcal{M}}: \mathcal{M}\rightarrow T^*Q, $ and
$\omega_{\mathcal{M}}= i_{\mathcal{M}}^* \omega $, that is, the
symplectic form $\omega_{\mathcal{M}}$ is induced from the canonical
symplectic form $\omega$ on $T^* Q$, where $i_{\mathcal{M}}^*:
T^*T^*Q \rightarrow T^*\mathcal{M}. $ For the distribution
$\mathcal{F}=(T\pi_Q)^{-1}(\mathcal{D})\subset TT^*Q,$ we define the
distribution $ \mathcal{K}=\mathcal {F}\cap T\mathcal{M}.$ Note that
$\mathcal{K}^\bot=\mathcal {F}^\bot\cap T\mathcal{M}, $ where
$\mathcal{K}^\bot$ denotes the symplectic orthogonal of
$\mathcal{K}$ with respect to the canonical symplectic form
$\omega$, and the admissibility condition $\mathrm{dim}\mathcal{M}=
\mathrm{rank}\mathcal{F}$ and the compatibility condition
$T\mathcal{M}\cap \mathcal{F}^\bot= \{0\}$ hold, then we know that the
restriction of the symplectic form $\omega_{\mathcal{M}}$ on
$T^*\mathcal{M}$ fibrewise to the distribution $\mathcal{K}$, that
is, $\omega_\mathcal{K}= \tau_{\mathcal{K}}\cdot
\omega_{\mathcal{M}}$ is non-degenerate, where $\tau_{\mathcal{K}}$
is the restriction map to distribution $\mathcal{K}$. It is worthy
of noting that $\omega_\mathcal{K}$ is not a true two-form on a
manifold, so it does not make sense to speak about it being closed.
We call $\omega_\mathcal{K}$ as a distributional two-form to avoid
any confusion. Because $\omega_\mathcal{K}$ is non-degenerate as a
bilinear form on each fibre of $\mathcal{K}$, there exists a vector
field $X_\mathcal{K}$ on $\mathcal{M}$ which takes values in the
constraint distribution $\mathcal{K}$, such that the following
nonholonomic constraints condition holds, that is,
\begin{align}\mathbf{i}_{X_\mathcal{K}}\omega_{\mathcal{K}}=\mathbf{d}H_{\mathcal{K}},
\label{3.1} \end{align}
where $\mathbf{d}H_\mathcal{K}$ is the restriction of
$\mathbf{d}H_\mathcal{M}$ to $\mathcal{K}$,
and the function $H_{\mathcal{K}}$ satisfies
$\mathbf{d}H_{\mathcal{K}}= \tau_{\mathcal{K}}\cdot \mathbf{d}H_{\mathcal {M}}$,
and $H_\mathcal{M}=
\tau_{\mathcal{M}}\cdot H$ is the restriction of $H$ to
$\mathcal{M}$. Then $(3.1)$ is called the distributional Hamiltonian
equation, see  Bates and $\acute{S}$niatycki \cite{basn93}. Thus,
the geometric formulation of a distributional Hamiltonian system may
be summarized as follows.

\begin{defi} (Distributional Hamiltonian System)
Assume that the 4-tuple $(T^*Q,\omega,\mathcal{D},H)$ is a nonholonomic
Hamiltonian system, where $\omega$ is the canonical
symplectic form on $T^* Q$, and $\mathcal{D}\subset TQ$ is a
$\mathcal{D}$-completely and $\mathcal{D}$-regularly nonholonomic
constraint of the system. If there exist a distribution
$\mathcal{K}$, an associated non-degenerate distributional two-form
$\omega_\mathcal{K}$ induced by the canonical symplectic form
and a vector field $X_\mathcal {K}$ on the
constraint submanifold $\mathcal{M}=\mathcal{F}L(\mathcal{D})\subset
T^*Q$, such that the distributional Hamiltonian equation
$\mathbf{i}_{X_\mathcal{K}}\omega_\mathcal{K}=\mathbf{d}H_\mathcal
{K}$ holds, where $\mathbf{d}H_\mathcal{K}$ is the restriction of
$\mathbf{d}H_\mathcal{M}$ to $\mathcal{K}$, and
the function $H_{\mathcal{K}}$ satisfies
$\mathbf{d}H_{\mathcal{K}}= \tau_{\mathcal{K}}\cdot \mathbf{d}H_{\mathcal {M}}$
as defined above,
then the triple $(\mathcal{K},\omega_{\mathcal{K}},H_{\mathcal{K}})$
is called a distributional Hamiltonian system of the nonholonomic
Hamiltonian system $(T^*Q,\omega,\mathcal{D},H)$, and $X_\mathcal
{K}$ is called a nonholonomic dynamical
vector field of the distributional Hamiltonian system
$(\mathcal{K},\omega_{\mathcal {K}},H_{\mathcal{K}})$. Under the above circumstances, we refer to
$(T^*Q,\omega,\mathcal{D},H)$ as a nonholonomic Hamiltonian system
with an associated distributional Hamiltonian system
$(\mathcal{K},\omega_{\mathcal {K}},H_{\mathcal{K}})$.
\end{defi}

Since the non-degenerate
distributional two-form $\omega_{\mathcal{K}}$ is not symplectic,
and the distributional Hamiltonian system
$(\mathcal{K},\omega_{\mathcal {K}},H_{\mathcal{K}})$ is not yet a Hamiltonian system,
and has no yet generating function,
and hence we can not describe the Hamilton-Jacobi equation for a
distributional Hamiltonian system just like as in Theorem 1.1.
But, for a given nonholonomic Hamiltonian system
$(T^*Q,\omega,\mathcal{D},H)$ with an associated distributional Hamiltonian
system $(\mathcal{K},\omega_{\mathcal {K}},H_{\mathcal{K}})$,
we can derive precisely the geometric constraint conditions of
the non-degenerate distributional two-form $\omega_\mathcal{K}$
for the nonholonomic dynamical vector field $X_\mathcal {K}$,
that is, the two types of Hamilton-Jacobi equation for the distributional Hamiltonian
system $(\mathcal{K},\omega_{\mathcal {K}},H_{\mathcal{K}})$. In order to do this, we need
first give two important notions and a key lemma, (see also Wang \cite{wa17}),
which is obtained by a careful modification for the
corresponding results of Abraham and Marsden in \cite{abma78}.
This lemma offers also an important tool for the proofs of the two types of Hamilton-Jacobi
theorems for the distributional Hamiltonian system and the nonholonomic
reduced distributional Hamiltonian system.\\

Let $Q$ be a smooth manifold and $TQ$ its tangent bundle, $T^* Q$
its cotangent bundle with the canonical symplectic form $\omega$,
and $\mathcal{D}\subset TQ$ is a $\mathcal{D}$-regularly nonholonomic
constraint, and
the projection $\pi_Q: T^* Q \rightarrow Q $ induces the map $
T\pi_{Q}: TT^* Q \rightarrow TQ. $ Assume that $\gamma: Q
\rightarrow T^*Q$ is an one-form on $Q$, if $\gamma$ is closed,
then $\mathbf{d}\gamma(x,y)=0, \; \forall\;
x, y \in TQ$. In the following we introduce two weaker notions.

\begin{defi}
\noindent $(\mathrm{i})$ The one-form $\gamma$ is called to be closed with respect to $T\pi_{Q}:
TT^* Q \rightarrow TQ, $ if for any $v, w \in TT^* Q, $ we have
$\mathbf{d}\gamma(T\pi_{Q}(v),T\pi_{Q}(w))=0; $\\

\noindent $(\mathrm{ii})$ The one-form $\gamma$ is called to be closed
on $\mathcal{D}$ with respect to $T\pi_{Q}:
TT^* Q \rightarrow TQ, $ if for any $v, w \in TT^* Q, $
and $T\pi_{Q}(v), \; T\pi_{Q}(w) \in \mathcal{D},$  we have
$\mathbf{d}\gamma(T\pi_{Q}(v),T\pi_{Q}(w))=0. $
\end{defi}

From the above definition we know that, the notion that
$\gamma$ is closed on $\mathcal{D}$ with respect to $T\pi_{Q}:
TT^* Q \rightarrow TQ, $ is weaker than the notion that $\gamma$ is closed
with respect to $T\pi_{Q}: TT^* Q \rightarrow TQ. $
From Wang \cite{wa17} we also know that the latter, that is, $\gamma$ is closed
with respect to $T\pi_{Q}: TT^* Q \rightarrow TQ, $
is weaker than the notion that $\gamma$ is closed. Thus,
the notion that $\gamma$ is closed on $\mathcal{D}$ with respect to $T\pi_{Q}:
TT^* Q \rightarrow TQ, $ is weaker than that $\gamma$ is closed on $\mathcal{D}$,
that is, $\mathbf{d}\gamma(x,y)=0, \; \forall\; x, y \in \mathcal{D}$.
In fact, if $\gamma$ is a closed one-form on $\mathcal{D}$,
then it must be closed on $\mathcal{D}$ with respect to
$T\pi_{Q}: TT^* Q \rightarrow TQ. $
Conversely, if $\gamma$ is closed on $\mathcal{D}$ with respect to
$T\pi_{Q}: TT^* Q \rightarrow TQ, $ then it may not be closed on $\mathcal{D}$.
We can prove a general result as follows.

\begin{prop}
Assume that $\gamma: Q \rightarrow T^*Q$ is an one-form on $Q$ and
it is not closed on $\mathcal{D}$. We define the set $N$, which is a subset of $TQ$,
such that the one-form $\gamma$ on $N$ satisfies the condition that
for any $x,y \in N, \; \mathbf{d}\gamma(x,y)\neq 0. $ Denote
$Ker(T\pi_Q)= \{u \in TT^*Q| \; T\pi_Q(u)=0 \}, $ and $T\gamma: TQ
\rightarrow TT^* Q .$ If $T\gamma(N)\subset Ker(T\pi_Q), $ then
$\gamma$ is closed with respect to $T\pi_{Q}: TT^* Q \rightarrow TQ.$
and hence $\gamma$ is closed on $\mathcal{D}$ with respect to
$T\pi_{Q}: TT^* Q \rightarrow TQ.$
\end{prop}

\noindent{\bf Proof: } In fact, for any $v, w \in TT^* Q, $ if
$T\pi_{Q}(v) \notin N, $ or $T\pi_{Q}(w))\notin N, $ then by the
definition of $N$, we know that
$\mathbf{d}\gamma(T\pi_{Q}(v),T\pi_{Q}(w))=0; $ If $T\pi_{Q}(v)\in
N, $ and $T\pi_{Q}(w))\in N, $ from the condition $T\gamma(N)\subset
Ker(T\pi_Q), $ we know that $T\pi_{Q}\cdot T\gamma \cdot
T\pi_{Q}(v)= T\pi_{Q}(v)=0, $ and $T\pi_{Q}\cdot T\gamma \cdot
T\pi_{Q}(w)= T\pi_{Q}(w)=0, $ where we have used the
relation $\pi_Q\cdot \gamma\cdot \pi_Q= \pi_Q, $ and hence
$\mathbf{d}\gamma(T\pi_{Q}(v),T\pi_{Q}(w))=0. $ Thus, for any $v, w
\in TT^* Q, $ we have always that
$\mathbf{d}\gamma(T\pi_{Q}(v),T\pi_{Q}(w))=0. $
In particular, for any $v, w \in TT^* Q, $
and $T\pi_{Q}(v), \; T\pi_{Q}(w) \in \mathcal{D},$  we have
$\mathbf{d}\gamma(T\pi_{Q}(v),T\pi_{Q}(w))=0. $
that is, $\gamma$ is closed on $\mathcal{D}$
with respect to $T\pi_{Q}: TT^* Q \rightarrow TQ. $
\hskip 0.3cm $\blacksquare$\\

Now, we prove the following Lemma 3.4. It is worthy of noting that
this lemma is an extension of Lemma 2.4 given in Wang \cite{wa17} to the
nonholonomic context.

\begin{lemm}
Assume that $\gamma: Q \rightarrow T^*Q$ is an one-form on $Q$, and
$\lambda=\gamma \cdot \pi_{Q}: T^* Q \rightarrow T^* Q .$ Then
we have that\\
\noindent $(\bf \mathrm{i})$ for any $x, y \in TQ, \;
\gamma^*\omega (x,y)= -\mathbf{d}\gamma (x,y),$ and for any $v, w \in
TT^* Q, $ \\ $ \lambda^*\omega(v,w)=
-\mathbf{d}\gamma(T\pi_{Q}(v), \; T\pi_{Q}(w)),$
since $\omega$ is the canonical symplectic form on $T^*Q$; \\
\noindent $(\bf \mathrm{ii})$ for any $v, w \in TT^* Q, \;
\omega(T\lambda \cdot v,w)= \omega(v, w-T\lambda \cdot
w)-\mathbf{d}\gamma(T\pi_{Q}(v), \; T\pi_{Q}(w))$ ;\\
\noindent $(\bf \mathrm{iii})$ If $L$ is $\mathcal{D}$-regular, and
$\textmd{Im}(\gamma)\subset \mathcal{M}=\mathcal{F}L(\mathcal{D}), $
then we have that $ X_{H}\cdot \gamma \in \mathcal{F}$ along
$\gamma$, and $ X_{H}\cdot \lambda \in \mathcal{F}$ along
$\lambda$, that is, $T\pi_{Q}(X_H\cdot\gamma(q))\in
\mathcal{D}_{q}, \; \forall q \in Q $, and $T\pi_{Q}(X_H\cdot\lambda(q,p))\in
\mathcal{D}_{q}, \; \forall q \in Q, \; (q,p) \in T^* Q. $
\end{lemm}

\noindent{\bf Proof:} The proofs of $(\mathrm{i})$ and $(\mathrm{ii})$
are given in Wang \cite{wa17}.
Now, we prove $(\mathrm{iii})$. For any $q \in Q
, \; (q,p)\in T^* Q, $ we have that
$$
X_H\cdot \gamma(q)= (\frac{\partial H}{\partial
p_i}\frac{\partial}{\partial q^i}-\frac{\partial H}{\partial
q^i}\frac{\partial}{\partial p_i})\gamma(q).
$$
and
$$
X_H\cdot \lambda(q,p)= (\frac{\partial H}{\partial
p_i}\frac{\partial}{\partial q^i}-\frac{\partial H}{\partial
q^i}\frac{\partial}{\partial p_i})\gamma\cdot \pi_Q(q,p).
$$
Then,
$$
T\pi_Q(X_H\cdot \gamma(q))=T\pi_Q(X_H\cdot \lambda(q,p))=(\frac{\partial H}{\partial
p_i}\frac{\partial}{\partial q^i})\gamma(q)= \gamma^*(\frac{\partial
H(q,p)}{\partial p_i})\frac{\partial}{\partial q^i},
$$
where $\gamma^*: T^*T^* Q\rightarrow T^* Q. $ Since
$\textmd{Im}(\gamma)\subset \mathcal{M}, $ and
$\gamma^*(\frac{\partial H(q,p)}{\partial p_i})\in
\mathcal{M}_{(q,p)}=\mathcal{F}L(\mathcal{D}_q), $ from $L$ is
$\mathcal{D}$-regular, $\mathcal{F}L$ is a diffeomorphism, then
there exists a $v_q\in \mathcal{D}_q, $ such that
$\mathcal{F}L(v_q)=\gamma^*(\frac{\partial H(q,p)}{\partial p_i}). $
Thus,
$$
T\pi_Q(X_H\cdot \gamma(q))=T\pi_Q(X_H\cdot
\lambda(q,p))=\mathcal{F}L(v_q)\frac{\partial}{\partial q^i} \in
\mathcal{D},
$$
it follows that $ X_{H}\cdot \gamma \in \mathcal{F}$ along
$\gamma$, and $ X_{H}\cdot \lambda \in \mathcal{F}$ along
$\lambda$. \hskip 0.3cm $\blacksquare$\\

By using the above Lemma 3.4,
we can derive precisely the geometric constraint conditions of
the non-degenerate distributional two-form $\omega_{\mathcal{K}}$
for the nonholonomic dynamical vector field $X_\mathcal {K}$,
that is, the following two types of
Hamilton-Jacobi equation for the distributional Hamiltonian
system $(\mathcal{K},\omega_{\mathcal {K}},H_{\mathcal{K}})$. At first, by using
the fact that the one-form $\gamma: Q
\rightarrow T^*Q $ is closed on $\mathcal{D}$ with respect to
$T\pi_Q: TT^* Q \rightarrow TQ, $
$\textmd{Im}(\gamma)\subset \mathcal{M}, $ and
$\textmd{Im}(T\gamma)\subset \mathcal{K}, $
we can prove the Type I of
Hamilton-Jacobi theorem for the distributional Hamiltonian system.
For convenience, the
maps involved in the following theorem and its proof are shown in
Diagram-1.

\begin{center}
\hskip 0cm \xymatrix{& \mathcal{M} \ar[d]_{X_{\mathcal{K}}}
\ar[r]^{i_{\mathcal{M}}} & T^* Q \ar[d]_{X_{H}}
 \ar[r]^{\pi_Q}
& Q \ar[d]_{X^\gamma} \ar[r]^{\gamma} & T^*Q \ar[d]^{X_H} \\
& \mathcal{K}  & T(T^*Q) \ar[l]_{\tau_{\mathcal{K}}} & TQ
\ar[l]_{T\gamma} & T(T^* Q)\ar[l]_{T\pi_Q}}
\end{center}
$$\mbox{Diagram-1}$$

\begin{theo} (Type I of Hamilton-Jacobi Theorem for a Distributional Hamiltonian System)
For the nonholonomic Hamiltonian system $(T^*Q,\omega,\mathcal{D},H)$
with an associated distributional Hamiltonian system
$(\mathcal{K},\omega_{\mathcal {K}},H_{\mathcal{K}})$, assume that $\gamma: Q
\rightarrow T^*Q$ is an one-form on $Q$, and $X^\gamma =
T\pi_{Q}\cdot X_H \cdot \gamma$, where $X_{H}$ is the dynamical
vector field of the corresponding unconstrained Hamiltonian system
$(T^*Q,\omega,H)$. Moreover, assume that $\textmd{Im}(\gamma)\subset
\mathcal{M}=\mathcal{F}L(\mathcal{D}), $ and $
\textmd{Im}(T\gamma)\subset \mathcal{K}. $ If the
one-form $\gamma: Q \rightarrow T^*Q $ is closed on $\mathcal{D}$ with respect to
$T\pi_Q: TT^* Q \rightarrow TQ, $ then $\gamma$ is a
solution of the equation $T\gamma \cdot
X^\gamma= X_{\mathcal{K}} \cdot \gamma. $ Here
$X_{\mathcal{K}}$
is the dynamical vector field of the distributional Hamiltonian system
$(\mathcal{K},\omega_{\mathcal {K}},H_{\mathcal{K}})$. The equation $T\gamma \cdot
X^\gamma= X_{\mathcal{K}} \cdot \gamma $ is called the Type I of
Hamilton-Jacobi equation for the distributional Hamiltonian system
$(\mathcal{K},\omega_{\mathcal {K}},H_{\mathcal{K}})$.
\end{theo}

\noindent{\bf Proof: } At first, we note that
$\textmd{Im}(\gamma)\subset \mathcal{M}, $ and
$\textmd{Im}(T\gamma)\subset \mathcal{K}, $ in this case,
$\omega_{\mathcal{K}}\cdot
\tau_{\mathcal{K}}=\tau_{\mathcal{K}}\cdot \omega_{\mathcal{M}}=
\tau_{\mathcal{K}}\cdot i_{\mathcal{M}}^* \cdot \omega, $ along
$\textmd{Im}(T\gamma)$. Thus, using the non-degenerate
distributional two-form $\omega_{\mathcal{K}}$, from Lemma 3.4 (ii) and (iii),
if we take that $v= X_{H}\cdot \gamma \in \mathcal{F},$ and for any $w
\in \mathcal{F}, \; T\lambda(w)\neq 0, $ and
$\tau_{\mathcal{K}}\cdot w \neq 0, $ then we have that
\begin{align*}
& \omega_{\mathcal{K}}(T\gamma \cdot X^\gamma, \;
\tau_{\mathcal{K}}\cdot w)=
\omega_{\mathcal{K}}(\tau_{\mathcal{K}}\cdot T\gamma \cdot
X^\gamma, \; \tau_{\mathcal{K}}\cdot w)\\ & =
\tau_{\mathcal{K}}\cdot i_{\mathcal{M}}^* \cdot\omega(T\gamma \cdot
X^\gamma, \; w ) = \tau_{\mathcal{K}}\cdot
i_{\mathcal{M}}^* \cdot\omega(T(\gamma \cdot \pi_Q)\cdot X_H \cdot \gamma, \; w)\\
& =\tau_{\mathcal{K}}\cdot i_{\mathcal{M}}^* \cdot(\omega(X_H \cdot
\gamma, \; w-T(\gamma
\cdot \pi_Q)\cdot w)-\mathbf{d}\gamma(T\pi_{Q}(X_H\cdot \gamma), \; T\pi_{Q}(w)))\\
& = \tau_{\mathcal{K}}\cdot i_{\mathcal{M}}^* \cdot\omega(X_H \cdot
\gamma, \; w) - \tau_{\mathcal{K}}\cdot i_{\mathcal{M}}^* \cdot
\omega(X_H \cdot \gamma, \; T(\gamma
\cdot \pi_Q) \cdot w)\\
& \;\;\;\;\;\; -\tau_{\mathcal{K}}\cdot i_{\mathcal{M}}^* \cdot\mathbf{d}\gamma(T\pi_{Q}(X_H\cdot \gamma), \; T\pi_{Q}(w))\\
& = \omega_{\mathcal{K}}( \tau_{\mathcal{K}}\cdot X_H \cdot \gamma,
\; \tau_{\mathcal{K}}\cdot w) -
\omega_{\mathcal{K}}(\tau_{\mathcal{K}}\cdot X_H \cdot \gamma, \;
\tau_{\mathcal{K}}\cdot T(\gamma
\cdot \pi_Q) \cdot w)\\
& \;\;\;\;\;\;
-\tau_{\mathcal{K}}\cdot i_{\mathcal{M}}^* \cdot\mathbf{d}\gamma(T\pi_{Q}(X_H\cdot \gamma), \; T\pi_{Q}(w))\\
& = \omega_{\mathcal{K}}(X_{\mathcal{K}}\cdot \gamma, \;
\tau_{\mathcal{K}} \cdot w) -
\omega_{\mathcal{K}}(X_{\mathcal{K}} \cdot \gamma, \;
\tau_{\mathcal{K}}\cdot T\gamma \cdot T\pi_{Q}(w))\\
& \;\;\;\;\;\; - \tau_{\mathcal{K}}\cdot
i_{\mathcal{M}}^* \cdot\mathbf{d}\gamma(T\pi_{Q}(X_H\cdot \gamma),
\; T\pi_{Q}(w)),
\end{align*}
where we have used that $ \tau_{\mathcal{K}}\cdot T\gamma= T\gamma, $ and
$\tau_{\mathcal{K}}\cdot X_H\cdot \gamma = X_{\mathcal{K}}\cdot
\gamma, $ since $\textmd{Im}(T\gamma)\subset \mathcal{K}. $
If the one-form $\gamma: Q \rightarrow T^*Q $ is closed on $\mathcal{D}$ with respect to
$T\pi_Q: TT^* Q \rightarrow TQ, $ then we have that
$\mathbf{d}\gamma(T\pi_{Q}(X_H\cdot \gamma), \; T\pi_{Q}(w))=0, $
since $X_{H}\cdot \gamma, \; w \in \mathcal{F},$
and $T\pi_{Q}(X_H\cdot \gamma), \; T\pi_{Q}(w) \in \mathcal{D}, $ and hence
$$
\tau_{\mathcal{K}}\cdot
i_{\mathcal{M}}^* \cdot\mathbf{d}\gamma(T\pi_{Q}(X_H\cdot \gamma),
\; T\pi_{Q}(w))=0,
$$
and
\begin{equation}
\omega_{\mathcal{K}}(T\gamma \cdot X^\gamma, \;
\tau_{\mathcal{K}}\cdot w)- \omega_{\mathcal{K}}(X_{\mathcal{K}}\cdot \gamma, \;
\tau_{\mathcal{K}} \cdot w)
= -\omega_{\mathcal{K}}(\tau_{\mathcal{K}}\cdot X_H \cdot \gamma, \;
\tau_{\mathcal{K}}\cdot T\gamma \cdot T\pi_{Q}(w)).
\label{3.2} \end{equation}
If $\gamma$ satisfies the equation $T\gamma\cdot X^\gamma= X_{\mathcal{K}}\cdot \gamma ,$
from Lemma 3.4(i) we deduce that
\begin{align*}
 -\omega_{\mathcal{K}}(X_{\mathcal{K}} \cdot \gamma, \;
\tau_{\mathcal{K}}\cdot T\gamma \cdot T\pi_{Q}(w))
& = -\omega_{\mathcal{K}}(T\gamma\cdot X^\gamma, \;
\tau_{\mathcal{K}}\cdot T\gamma \cdot T\pi_{Q}(w))\\
& = -\omega_{\mathcal{K}}(\tau_{\mathcal{K}}\cdot T\gamma\cdot X^\gamma, \;
\tau_{\mathcal{K}}\cdot T\gamma \cdot T\pi_{Q}(w))\\
& = -\tau_{\mathcal{K}}\cdot
i_{\mathcal{M}}^* \cdot \omega(T\gamma \cdot T\pi_{Q}(X_{H}\cdot\gamma), \; T\gamma \cdot T\pi_{Q}(w))\\
& = -\tau_{\mathcal{K}}\cdot
i_{\mathcal{M}}^* \cdot\gamma^*\omega( T\pi_{Q}(X_{H}\cdot\gamma), \; T\pi_{Q}(w))\\
& = \tau_{\mathcal{K}}\cdot i_{\mathcal{M}}^* \cdot
\mathbf{d}\gamma(T\pi_{Q}( X_{H}\cdot\gamma ), \; T\pi_{Q}(w))=0.
\end{align*}
Because the distributional two-form $\omega_{\mathcal{K}}$ is non-degenerate,
the left side of (3.2) equals zero, only when
$\gamma$ satisfies the equation
$T\gamma\cdot X^\gamma= X_{\mathcal{K}}\cdot \gamma .$ Thus,
if the one-form $\gamma: Q \rightarrow T^*Q $ is closed on $\mathcal{D}$ with respect to
$T\pi_Q: TT^* Q \rightarrow TQ, $ then $\gamma$ must be a solution of the Type I of Hamilton-Jacobi equation
$T\gamma\cdot X^\gamma= X_{\mathcal{K}}\cdot \gamma .$
\hskip 0.3cm $\blacksquare$\\

Next, for any symplectic map $\varepsilon: T^* Q \rightarrow T^* Q $,
we can prove the following Type II of
Hamilton-Jacobi theorem for the distributional Hamiltonian system.
For convenience, the
maps involved in the following theorem and its proof are shown in
Diagram-2.

\begin{center}
\hskip 0cm \xymatrix{& \mathcal{M} \ar[d]_{X_{\mathcal{K}}}
\ar[r]^{i_{\mathcal{M}}} & T^* Q \ar[d]_{X_{H\cdot \varepsilon}}
\ar[dr]^{X^\varepsilon} \ar[r]^{\pi_Q}
& Q \ar[r]^{\gamma} & T^*Q \ar[d]^{X_H} \\
& \mathcal{K}  & T(T^*Q) \ar[l]_{\tau_{\mathcal{K}}} & TQ
\ar[l]_{T\gamma} & T(T^* Q)\ar[l]_{T\pi_Q}}
\end{center}
$$\mbox{Diagram-2}$$

\begin{theo} (Type II of Hamilton-Jacobi Theorem for a Distributional Hamiltonian System)
For the nonholonomic Hamiltonian system $(T^*Q,\omega,\mathcal{D},H)$
with an associated distributional Hamiltonian system
$(\mathcal{K},\omega_{\mathcal {K}},H_{\mathcal{K}})$, assume that $\gamma: Q
\rightarrow T^*Q$ is an one-form on $Q$, and $\lambda=\gamma \cdot
\pi_{Q}: T^* Q \rightarrow T^* Q, $ and for any
symplectic map $\varepsilon: T^* Q \rightarrow T^* Q $, denote by
$X^\varepsilon = T\pi_{Q}\cdot X_H \cdot \varepsilon$, where $X_{H}$ is the dynamical
vector field of the corresponding unconstrained Hamiltonian system
$(T^*Q,\omega,H)$. Moreover, assume that $\textmd{Im}(\gamma)\subset
\mathcal{M}=\mathcal{F}L(\mathcal{D}), $ and $
\textmd{Im}(T\gamma)\subset \mathcal{K}. $ If $\varepsilon$ is a
solution of the equation
$\tau_{\mathcal{K}}\cdot T\varepsilon(X_{H\cdot\varepsilon})= T\lambda \cdot X_H\cdot\varepsilon,$
if and only if it is a solution of the equation
$T\gamma \cdot
X^\varepsilon= X_{\mathcal{K}} \cdot \varepsilon $. Here $
X_{H\cdot\varepsilon}$ is the Hamiltonian vector field of the function
$H\cdot\varepsilon: T^* Q\rightarrow \mathbb{R}, $ and $X_{\mathcal{K}}$
is the dynamical vector field of the distributional Hamiltonian system
$(\mathcal{K},\omega_{\mathcal {K}},H_{\mathcal{K}})$. The equation $T\gamma \cdot
X^\varepsilon= X_{\mathcal{K}} \cdot \varepsilon,$ is called the Type II of
Hamilton-Jacobi equation for the distributional Hamiltonian system
$(\mathcal{K},\omega_{\mathcal {K}},H_{\mathcal{K}})$.
\end{theo}

\noindent{\bf Proof: } In the same way, we note that
$\textmd{Im}(\gamma)\subset \mathcal{M}, $ and
$\textmd{Im}(T\gamma)\subset \mathcal{K}, $ in this case,
$\omega_{\mathcal{K}}\cdot
\tau_{\mathcal{K}}=\tau_{\mathcal{K}}\cdot \omega_{\mathcal{M}}=
\tau_{\mathcal{K}}\cdot i_{\mathcal{M}}^* \cdot \omega, $ along
$\textmd{Im}(T\gamma)$. Thus, using the non-degenerate
distributional two-form $\omega_{\mathcal{K}}$, from Lemma 3.4, if
we take that $v= \tau_{\mathcal{K}}\cdot X_H\cdot \varepsilon = X_{\mathcal{K}}\cdot
\varepsilon \in \mathcal{K} (\subset \mathcal{F}), $ and for any $w
\in \mathcal{F}, \; T\lambda(w)\neq 0, $ and
$\tau_{\mathcal{K}}\cdot w \neq 0, $ then we have that
\begin{align*}
& \omega_{\mathcal{K}}(T\gamma \cdot X^\varepsilon, \;
\tau_{\mathcal{K}}\cdot w)=
\omega_{\mathcal{K}}(\tau_{\mathcal{K}}\cdot T\gamma \cdot
X^\varepsilon, \; \tau_{\mathcal{K}}\cdot w)\\ & =
\tau_{\mathcal{K}}\cdot i_{\mathcal{M}}^* \cdot\omega(T\gamma \cdot
X^\varepsilon, \; w ) = \tau_{\mathcal{K}}\cdot
i_{\mathcal{M}}^* \cdot\omega(T(\gamma \cdot \pi_Q)\cdot X_H \cdot \varepsilon, \; w)\\
& =\tau_{\mathcal{K}}\cdot i_{\mathcal{M}}^* \cdot(\omega(X_H \cdot
\varepsilon, \; w-T(\gamma
\cdot \pi_Q)\cdot w)-\mathbf{d}\gamma(T\pi_{Q}(X_H\cdot \varepsilon), \; T\pi_{Q}(w)))\\
& = \tau_{\mathcal{K}}\cdot i_{\mathcal{M}}^* \cdot\omega(X_H \cdot
\varepsilon, \; w) - \tau_{\mathcal{K}}\cdot i_{\mathcal{M}}^* \cdot
\omega(X_H \cdot \varepsilon, \; T\lambda \cdot w)\\
& \;\;\;\;\;\;
-\tau_{\mathcal{K}}\cdot i_{\mathcal{M}}^* \cdot\mathbf{d}\gamma(T\pi_{Q}(X_H\cdot \varepsilon), \; T\pi_{Q}(w))\\
& = \omega_{\mathcal{K}}( \tau_{\mathcal{K}}\cdot X_H \cdot \varepsilon,
\; \tau_{\mathcal{K}}\cdot w) -
\omega_{\mathcal{K}}(\tau_{\mathcal{K}}\cdot X_H \cdot \varepsilon, \;
\tau_{\mathcal{K}}\cdot T\lambda \cdot w)\\
& \;\;\;\;\;\;
+\tau_{\mathcal{K}}\cdot i_{\mathcal{M}}^* \cdot \lambda^* \omega(X_H\cdot \varepsilon, \; w)\\
& = \omega_{\mathcal{K}}(X_{\mathcal{K}}\cdot \varepsilon, \;
\tau_{\mathcal{K}} \cdot w) -
\omega_{\mathcal{K}}(\tau_{\mathcal{K}}\cdot X_H \cdot \varepsilon,
\; T\lambda \cdot w)+ \omega_{\mathcal{K}}(T\lambda\cdot X_H\cdot \varepsilon,
\; T\lambda \cdot w),
\end{align*}
where we have used that
$ \tau_{\mathcal{K}}\cdot T\gamma= T\gamma, \; \tau_{\mathcal{K}}\cdot T\lambda= T\lambda, $ and
$\tau_{\mathcal{K}}\cdot X_H\cdot \varepsilon = X_{\mathcal{K}}\cdot
\varepsilon, $ since $\textmd{Im}(T\gamma)\subset \mathcal{K}. $ Note
that $\varepsilon: T^* Q \rightarrow T^* Q $ is symplectic, and $
X_H\cdot \varepsilon = T\varepsilon \cdot X_{H\cdot\varepsilon}, $ along
$\varepsilon$, and hence $\tau_{\mathcal{K}}\cdot X_H \cdot \varepsilon=
\tau_{\mathcal{K}}\cdot T\varepsilon \cdot X_{H \cdot \varepsilon}, $ along $\varepsilon$.
Then we have that
\begin{align*}
& \omega_{\mathcal{K}}(T\gamma \cdot X^\varepsilon, \;
\tau_{\mathcal{K}}\cdot w)-
\omega_{\mathcal{K}}(X_{\mathcal{K}}\cdot \varepsilon, \;
\tau_{\mathcal{K}} \cdot w) \nonumber \\
& = -
\omega_{\mathcal{K}}(\tau_{\mathcal{K}}\cdot X_H \cdot \varepsilon, \;
 T\lambda \cdot w)+ \omega_{\mathcal{K}}(T\lambda\cdot X_H\cdot \varepsilon,
\; T\lambda \cdot w)\\
&= \omega_{\mathcal{K}}(T\lambda\cdot X_H\cdot \varepsilon
-\tau_{\mathcal{K}}\cdot T\varepsilon \cdot X_{H \cdot \varepsilon},
\; T\lambda \cdot w).
\end{align*}
Because the distributional two-form
$\omega_{\mathcal{K}}$ is non-degenerate, it follows that the equation
$T\gamma\cdot X^\varepsilon= X_{\mathcal{K}}\cdot
\varepsilon ,$ is equivalent to the equation
$\tau_{\mathcal{K}}\cdot T\varepsilon \cdot X_{H\cdot\varepsilon} = T\lambda\cdot X_H\cdot \varepsilon $.
Thus, $\varepsilon$ is a solution of the equation
$\tau_{\mathcal{K}}\cdot T\varepsilon\cdot X_{H\cdot\varepsilon}= T\lambda \cdot X_H \cdot\varepsilon,$
if and only if it is a solution of
the Type II of Hamilton-Jacobi equation $T\gamma\cdot X^\varepsilon= X_{\mathcal{K}}\cdot
\varepsilon .$
\hskip 0.3cm $\blacksquare$

\begin{rema}
If the nonholonomic Hamiltonian system we considered has not any constrains, in this case,
the distributional Hamiltonian system is just the Hamiltonian system itself.
From the above Type I and Type II of Hamilton-Jacobi theorems, that is,
Theorem 3.5 and Theorem 3.6, we can get the Theorem 2.5
and Theorem 2.6 in Wang \cite{wa17}.
It shows that Theorem 3.5 and Theorem 3.6 can be regarded as an extension of two types of geometric
Hamilton-Jacobi theorem for Hamiltonian system given in \cite{wa17} to the
nonholonomic context. In particular, in this case,
if the one-form $\gamma$ is
given by a generating function of a symplectic map, then the
classical Hamilton-Jacobi equation $H(q, \gamma(q))=E,$(constant in $t$),
or equivalently, $\mathbf{d}(H\cdot \gamma)=0, $ as well as the Type I of
Hamilton-Jacobi equation $T\gamma \cdot X^\gamma= X_H \cdot\gamma, $ and the Type II of
Hamilton-Jacobi theorem, all of them hold, see Wang \cite{wa17}.
\end{rema}

\begin{rema}
It is worthy of note that the formulations of Type I
and Type II of Hamilton-Jacobi equation for a distributional
Hamiltonian system, given by Theorem 3.5 and Theorem 3.6, have
more extensive sense, because, in general, the one-form $\gamma$ is not
given by a generating function of a symplectic map. When $\gamma$ is a solution
of the classical Hamilton-Jacobi equation, that is, $X_H\cdot
\gamma=0, $ which is equivalent to the equation $\mathbf{d}(H\cdot
\gamma)=0, $ or $H(q, \gamma(q))=E, \; q\in Q, $ and $E$ is a
constant, in this case, $X_H^\gamma= T\pi_Q\cdot X_H\cdot \gamma=0,$
and hence from the Type I of Hamilton-Jacobi equation, we have that
$X_{\mathcal{K}}\cdot \gamma= T\gamma\cdot X_H^\gamma=0. $ Since the
classical Hamilton-Jacobi equation $X_H\cdot \gamma=0, $ shows that
the dynamical vector field of the corresponding unconstrained
Hamiltonian system $(T^*Q,\omega,H)$ is degenerate along $\gamma$,
then the equation $X_{\mathcal{K}}\cdot \gamma=0,$ shows that the
dynamical vector field of the distributional Hamiltonian system
$(\mathcal{K},\omega_{\mathcal {K}},H_{\mathcal{K}})$ is degenerate along
$\gamma$. The equation $X_{\mathcal{K}}\cdot \gamma=0$ is called the
classical Hamilton-Jacobi equation for the distributional
Hamiltonian system $(\mathcal{K},\omega_{\mathcal {K}},H_{\mathcal{K}}). $ In
addition, for a symplectic map $\varepsilon: T^* Q \rightarrow T^* Q
$, if $X_H\cdot \varepsilon=0, $ then from the Type II of
Hamilton-Jacobi equation, we have that $X_{\mathcal{K}}\cdot
\varepsilon= T\gamma\cdot X_H^\varepsilon=0. $ But, from the
equation $\tau_{\mathcal{K}}\cdot T\varepsilon\cdot
X_{H\cdot\varepsilon}= T\lambda \cdot X_H \cdot\varepsilon,$ we know
that $X_{\mathcal{K}}\cdot \varepsilon=0 $ is not equivalent to
$X_{H\cdot\varepsilon}=0.$
\end{rema}

\section{Hamilton-Jacobi Equations for a Reduced Distributional Hamiltonian System }

It is well-known that the reduction of nonholonomically constrained mechanical systems
is very important subject in geometric mechanics, and it is also
regarded as a useful tool for simplifying and studying
concrete nonholonomic systems, see Koiller \cite{ko92},
Bates and $\acute{S}$niatycki \cite{basn93},  Cantrijn et al.
\cite{calemama99, calemama98}, Cushman et al. \cite{cudusn10} and \cite{cukesnba95},
Cendra et al. \cite{cemara01}, Bloch et al. \cite{blkrmamu96}
and de Le\'{o}n and Rodrigues \cite{lero89} and so on.\\

In this section, for a nonholonomic Hamiltonian
system with symmetry $(T^*Q,G,\omega,\mathcal{D},H)$,
where $\omega$ is the canonical
symplectic form on $T^* Q$, and $\mathcal{D}\subset TQ$ is a
$\mathcal{D}$-completely and $\mathcal{D}$-regularly nonholonomic
constraint of the system, and $\mathcal{D}$ and $H$ are both
$G$-invariant, we first give the nonholonomic reduction of the system,
and a reduced distribution $\bar{\mathcal{K}}$, an associated
non-degenerate and nonholonomic reduced distributional two-form $\omega_{\bar{\mathcal{K}}}$
induced by the canonical symplectic form $\omega$,
and a nonholonomic reduced distributional Hamiltonian system.
Then we derive precisely
the geometric constraint conditions of the reduced distributional two-form
$\omega_{\bar{\mathcal{K}}}$ for the nonholonomic dynamical vector field,
that is, the two types of Hamilton-Jacobi equation for the
nonholonomic reduced distributional Hamiltonian system, which are an
extension of the above two types of Hamilton-Jacobi equation for the distributional
Hamiltonian system under nonholonomic reduction. In order to do
this, we need first to describe the nonholonomic reduction, and
analyze carefully the dynamics and structure of
the nonholonomic Hamiltonian system with symmetry,
following the results given by
Bates and $\acute{S}$niatycki in \cite{basn93}, see also, Cushman et
al. \cite{cudusn10} and \cite{cukesnba95}.\\

Now, we assume that the
5-tuple $(T^*Q,G,\omega,\mathcal{D},H)$ is a
$\mathcal{D}$-completely and $\mathcal{D}$-regularly nonholonomic
Hamiltonian system with symmetry, and the Lie group $G$ acts
smoothly on $Q$ by the left, and we also consider the natural lifted actions on $TQ$ and $T^* Q$,
and assume that the cotangent lifted action on $T^\ast Q$ is free, proper and
symplectic. The orbit space $T^* Q/ G$ is a smooth manifold and the
canonical projection $\pi_{/G}: T^* Q \rightarrow T^* Q /G $ is
a surjective submersion.\\

Assume that $H: T^*Q \rightarrow \mathbb{R}$ is a $G$-invariant
Hamiltonian, and that the $\mathcal{D}$-completely and
$\mathcal{D}$-regularly nonholonomic constraints $\mathcal{D}\subset
TQ$ is a $G$-invariant distribution, that is, the tangent of group action maps
$\mathcal{D}_q$ to $\mathcal{D}_{gq}$ for any
$q\in Q $. Note that the Legendre transformation $\mathcal{F}L: TQ
\rightarrow T^*Q$ is a fiber-preserving map, from \S2 we know that
$\mathcal{M}=\mathcal{F}L(\mathcal{D})\subset T^*Q$ is
$G$-invariant, and the quotient space
$\bar{\mathcal{M}}=\mathcal{M}/G$ of the $G$-orbit in $\mathcal{M}$
is a smooth manifold with projection $\pi_{/G}:
\mathcal{M}\rightarrow \bar{\mathcal{M}}( \subset T^* Q /G),$ which
is a surjective submersion.\\

Since $G$ is the symmetry group of the system, all intrinsically
defined vector fields and distributions push down to
$\bar{\mathcal{M}}$. In particular, the vector field $X_\mathcal{M}$
on $\mathcal{M}$ pushes down to a vector field
$X_{\bar{\mathcal{M}}}=T\pi_{/G}\cdot X_\mathcal{M}$, and the
distribution $\mathcal{K}$ pushes down to a distribution
$T\pi_{/G}\cdot \mathcal{K}$ on $\bar{\mathcal{M}}$, and the
Hamiltonian $H$ pushes down to $h_{\bar{\mathcal{M}}}$, such that
$h_{\bar{\mathcal{M}}}\cdot \pi_{/G}=
\tau_{\mathcal{M}}\cdot H$. However, $\omega_\mathcal{K}$ need not
push down to a distributional two-form defined on $T\pi_{/G}\cdot
\mathcal{K}$, despite of the fact that $\omega_\mathcal{K}$ is
$G$-invariant. This is because there may be infinitesimal symmetry
$\eta_{\mathcal{K}}$ that lies in $\mathcal{M}$, such that
$\mathbf{i}_{\eta_\mathcal{K}} \omega_\mathcal{K}\neq 0$. From Bates
and $\acute{S}$niatycki \cite{basn93}, we know that to eliminate
this difficulty, $\omega_\mathcal{K}$ is restricted to a
sub-distribution $\mathcal{U}$ of $\mathcal{K}$ defined by
$$\mathcal{U}=\{u\in\mathcal{K} \; | \; \omega_\mathcal{K}(u,v)
=0,\quad \forall \; v \in \mathcal{V}\cap \mathcal{K}\},$$ where
$\mathcal{V}$ is the distribution on $\mathcal{M}$ tangent to the
orbits of $G$ in $\mathcal{M}$ and it is spanned by the infinitesimal
symmetries. Clearly, $\mathcal{U}$ and $\mathcal{V}$ are both
$G$-invariant, project down to $\bar{\mathcal{M}}$ and
$T\pi_{/G}\cdot \mathcal{V}=0$, and define the distribution $\bar{\mathcal{K}}$ by
$\bar{\mathcal{K}}= T\pi_{/G}\cdot \mathcal{U}$. Moreover, we take
that $\omega_\mathcal{U}= \tau_{\mathcal{U}}\cdot
\omega_{\mathcal{M}}$ is the restriction of the induced symplectic form
$\omega_{\mathcal{M}}$ on $T^*\mathcal{M}$ fibrewise to the
distribution $\mathcal{U}$, where $\tau_{\mathcal{U}}$ is the
restriction map to distribution $\mathcal{U}$, and the
$\omega_{\mathcal{U}}$ pushes down to a non-degenerate
distributional two-form $\omega_{\bar{\mathcal{K}}}$ on
$\bar{\mathcal{K}}$, such that $\pi_{/G}^*
\omega_{\bar{\mathcal{K}}}= \omega_{\mathcal{U}}$. Because
$\omega_{\bar{\mathcal{K}}}$ is non-degenerate as a bilinear form on
each fibre of $\bar{\mathcal{K}}$, there exists a vector field
$X_{\bar{\mathcal{K}}}$ on $\bar{\mathcal{M}}$ which takes values in
the constraint distribution $\bar{\mathcal{K}}$, such that the
reduced distributional Hamiltonian equation holds, that is,
$\mathbf{i}_{X_{\bar{\mathcal{K}}}}\omega_{\bar{\mathcal{K}}}=\mathbf{d}h_{\bar{\mathcal{K}}}$,
where $\mathbf{d}h_{\bar{\mathcal{K}}}$ is the restriction of
$\mathbf{d}h_{\bar{\mathcal{M}}}$ to $\bar{\mathcal{K}}$ and
the function $h_{\bar{\mathcal{K}}}$ satisfies
$\mathbf{d}h_{\bar{\mathcal{K}}}= \tau_{\bar{\mathcal{K}}}\cdot \mathbf{d}h_{\bar{\mathcal{M}}}$,
and $h_{\bar{\mathcal{M}}}\cdot \pi_{/G}= H_{\mathcal{M}}$ and
$H_{\mathcal{M}}$ is the restriction of the Hamiltonian function $H$
to $\mathcal{M}$. In addition, the vector fields $X_{\mathcal{K}}$
and $X_{\bar{\mathcal{K}}}$ are $\pi_{/G}$-related. Thus, the
geometrical formulation of a nonholonomic reduced distributional
Hamiltonian system may be summarized as follows.

\begin{defi} (Nonholonomic Reduced Distributional Hamiltonian System)
Assume that the 5-tuple $(T^*Q,G,\omega,\mathcal{D},H)$ is a nonholonomic
Hamiltonian system with symmetry, where $\omega$ is the canonical
symplectic form on $T^* Q$, and $\mathcal{D}\subset TQ$ is a
$\mathcal{D}$-completely and $\mathcal{D}$-regularly nonholonomic
constraint of the system, and $\mathcal{D}$ and $H$ are both
$G$-invariant. If there exists a nonholonomic reduced distribution $\bar{\mathcal{K}}$,
an associated non-degenerate  and nonholonomic reduced
distributional two-form $\omega_{\bar{\mathcal{K}}}$
and a vector field $X_{\bar{\mathcal {K}}}$ on the reduced constraint
submanifold $\bar{\mathcal{M}}=\mathcal{M}/G, $ where
$\mathcal{M}=\mathcal{F}L(\mathcal{D})\subset T^*Q$, such that the
nonholonomic reduced distributional Hamiltonian equation
$\mathbf{i}_{X_{\bar{\mathcal{K}}}}\omega_{\bar{\mathcal{K}}} =
\mathbf{d}h_{\bar{\mathcal{K}}}$ holds,
where $\mathbf{d}h_{\bar{\mathcal{K}}}$ is the restriction of
$\mathbf{d}h_{\bar{\mathcal{M}}}$ to $\bar{\mathcal{K}}$ and
the function $h_{\bar{\mathcal{K}}}$ satisfies
$\mathbf{d}h_{\bar{\mathcal{K}}}= \tau_{\bar{\mathcal{K}}}\cdot \mathbf{d}h_{\bar{\mathcal{M}}}$
and $h_{\bar{\mathcal{M}}}\cdot \pi_{/G}= H_{\mathcal{M}}$ as defined above.
Then the triple $(\bar{\mathcal{K}},\omega_{\bar{\mathcal {K}}},h_{\bar{\mathcal{K}}})$
is called a nonholonomic reduced distributional Hamiltonian system
of the nonholonomic Hamiltonian system with symmetry
$(T^*Q,G,\omega,\mathcal{D},H)$, and $X_{\bar{\mathcal {K}}}$ is
called a dynamical vector field of the
nonholonomic reduced distributional Hamiltonian system
$(\bar{\mathcal{K}},\omega_{\bar{\mathcal{K}}},h_{\bar{\mathcal{K}}})$. Under the above
circumstances, we refer to $(T^*Q,G,\omega,\mathcal{D},H)$ as a
nonholonomic reducible Hamiltonian system with an associated
nonholonomic reduced distributional Hamiltonian system
$(\bar{\mathcal{K}},\omega_{\bar{\mathcal{K}}},h_{\bar{\mathcal{K}}})$.
\end{defi}

Since the non-degenerate and nonholonomic reduced distributional two-form
$\omega_{\bar{\mathcal{K}}}$ is not symplectic,
and the nonholonomic reduced distributional Hamiltonian system
$(\bar{\mathcal{K}},\omega_{\bar{\mathcal{K}}},h_{\bar{\mathcal{K}}})$ is not yet a Hamiltonian system,
and has no yet generating function,
and hence we can not describe the Hamilton-Jacobi equation for the nonholonomic reduced
distributional Hamiltonian system just like as in Theorem 1.1.
But, for a given nonholonomic reducible Hamiltonian system
$(T^*Q,G,\omega,\mathcal{D},H)$ with an associated
nonholonomic reduced distributional Hamiltonian system
$(\bar{\mathcal{K}},\omega_{\bar{\mathcal {K}}},h_{\bar{\mathcal{K}}})$, by using Lemma 3.4,
we can derive precisely
the geometric constraint conditions of the nonholonomic reduced distributional two-form
$\omega_{\bar{\mathcal{K}}}$ for the dynamical vector field $X_{\bar{\mathcal {K}}}$,
that is, the two types of Hamilton-Jacobi equation for the
nonholonomic reduced distributional Hamiltonian system
$(\bar{\mathcal{K}},\omega_{\bar{\mathcal {K}}},h_{\bar{\mathcal{K}}})$.
At first, using the fact that the one-form $\gamma: Q
\rightarrow T^*Q $ is closed on $\mathcal{D}$ with respect to
$T\pi_Q: TT^* Q \rightarrow TQ, $
$\textmd{Im}(\gamma)\subset \mathcal{M}, $ and it is $G$-invariant,
as well as $ \textmd{Im}(T\gamma)\subset \mathcal{K}, $
we can prove the Type I of
Hamilton-Jacobi theorem for the nonholonomic reduced distributional
Hamiltonian system. For convenience, the maps involved in the
following theorem and its proof are shown in Diagram-3.
\begin{center}
\hskip 0cm \xymatrix{ & \mathcal{M} \ar[d]_{X_{\mathcal{K}}}
\ar[r]^{i_{\mathcal{M}}} & T^* Q \ar[d]_{X_{H}}
 \ar[r]^{\pi_Q}
  & Q \ar[d]_{X_H^\gamma} \ar[r]^{\gamma}
  & T^* Q \ar[d]_{X_H} \ar[r]^{\pi_{/G}} & T^* Q/G \ar[d]_{X_{h}}
  & \mathcal{\bar{M}} \ar[l]_{i_{\mathcal{\bar{M}}}} \ar[d]_{X_{\mathcal{\bar{K}}}}\\
  & \mathcal{K}
  & T(T^*Q) \ar[l]_{\tau_{\mathcal{K}}}
  & TQ \ar[l]_{T\gamma}
  & T(T^* Q) \ar[l]_{T\pi_Q} \ar[r]^{T\pi_{/G}}
  & T(T^* Q/G) \ar[r]^{\tau_{\mathcal{\bar{K}}}} & \mathcal{\bar{K}} }
\end{center}
$$\mbox{Diagram-3}$$

\begin{theo} (Type I of Hamilton-Jacobi Theorem for a Nonholonomic
Reduced Distributional Hamiltonian System)
For a given nonholonomic reducible Hamiltonian system
$(T^*Q,G,\omega,\mathcal{D},H)$ with an associated
nonholonomic reduced distributional Hamiltonian system
$(\bar{\mathcal{K}},\omega_{\bar{\mathcal{K}}},h_{\bar{\mathcal{K}}})$, assume that
$\gamma: Q \rightarrow T^*Q$ is an one-form on $Q$, and
$X_H^\gamma = T\pi_{Q}\cdot X_H \cdot \gamma$, where $X_{H}$ is
the dynamical vector field of the corresponding unconstrained
Hamiltonian system with symmetry $(T^*Q,G,\omega,H)$. Moreover,
assume that $\textmd{Im}(\gamma)\subset \mathcal{M}, $ and it is
$G$-invariant, $ \textmd{Im}(T\gamma)\subset \mathcal{K}, $ and
$\bar{\gamma}=\pi_{/G}(\gamma): Q \rightarrow T^* Q/G .$ If the
one-form $\gamma: Q \rightarrow T^*Q $ is closed on $\mathcal{D}$ with respect to
$T\pi_Q: TT^* Q \rightarrow TQ, $ then $\bar{\gamma}$ is a solution
of the equation $T\bar{\gamma}\cdot X_H^ \gamma =
X_{\bar{\mathcal{K}}}\cdot \bar{\gamma}. $
Here $X_{\bar{\mathcal{K}}}$ is the dynamical vector
field of the nonholonomic reduced distributional Hamiltonian system
$(\bar{\mathcal{K}},\omega_{\bar{\mathcal{K}}},h_{\bar{\mathcal{K}}})$. The equation $
T\bar{\gamma}\cdot X_H^ \gamma = X_{\bar{\mathcal{K}}}\cdot
\bar{\gamma},$ is called the Type I of Hamilton-Jacobi equation for
the nonholonomic reduced distributional Hamiltonian system
$(\bar{\mathcal{K}},\omega_{\bar{\mathcal{K}}},h_{\bar{\mathcal{K}}})$.
\end{theo}

\noindent{\bf Proof: } At first, from Theorem 3.5, we know that
$\gamma$ is a solution of the Hamilton-Jacobi equation
$T\gamma\cdot X_H^\gamma= X_{\mathcal{K}}\cdot \gamma .$ Next, we note that
$\textmd{Im}(\gamma)\subset \mathcal{M}, $ and it is $G$-invariant,
$ \textmd{Im}(T\gamma)\subset \mathcal{K}, $ and hence
$\textmd{Im}(T\bar{\gamma})\subset \bar{\mathcal{K}}, $ in this case,
$\pi^*_{/G}\cdot\omega_{\bar{\mathcal{K}}}\cdot\tau_{\bar{\mathcal{K}}}= \tau_{\mathcal{U}}\cdot
\omega_{\mathcal{M}}= \tau_{\mathcal{U}}\cdot i_{\mathcal{M}}^*\cdot
\omega, $ along $\textmd{Im}(T\bar{\gamma})$. Thus, using the
non-degenerate distributional two-form $\omega_{\bar{\mathcal{K}}}$,
from Lemma 3.4, if we take that
$v= \tau_{\bar{\mathcal{K}}}\cdot T\pi_{/G}\cdot X_{H}\cdot \gamma =X_{\bar{\mathcal{K}}} \cdot
\bar{\gamma} \in \bar{\mathcal{K}},
$ and for any $w \in \mathcal{F}, \; T\lambda(w)\neq 0, $ and
$\tau_{\bar{\mathcal{K}}}\cdot T\pi_{/G}\cdot w \neq 0, $ then we have that
\begin{align*}
& \omega_{\bar{\mathcal{K}}}(T\bar{\gamma} \cdot X_H^\gamma, \;
\tau_{\bar{\mathcal{K}}}\cdot T\pi_{/G} \cdot w)
= \omega_{\bar{\mathcal{K}}}(\tau_{\bar{\mathcal{K}}}\cdot T(\pi_{/G} \cdot
\gamma) \cdot X_H^\gamma, \; \tau_{\bar{\mathcal{K}}}\cdot T\pi_{/G} \cdot w )\\
& = \pi^*_{/G}
\cdot \omega_{\bar{\mathcal{K}}}\cdot\tau_{\bar{\mathcal{K}}}(T\gamma \cdot X_H^\gamma, \; w) =
\tau_{\mathcal{U}}\cdot i^*_{\mathcal{M}} \cdot\omega(T\gamma \cdot
X_H^\gamma, \; w)\\
& = \tau_{\mathcal{U}}\cdot i^*_{\mathcal{M}} \cdot
\omega(T(\gamma
\cdot \pi_Q)\cdot X_H \cdot \gamma, \; w) \\
& = \tau_{\mathcal{U}}\cdot i^*_{\mathcal{M}} \cdot
(\omega(X_H \cdot \gamma, \; w-T(\gamma \cdot \pi_Q)\cdot w)- \mathbf{d}\gamma(T\pi_{Q}(X_H\cdot \gamma), \; T\pi_{Q}(w)))\\
& = \tau_{\mathcal{U}}\cdot i^*_{\mathcal{M}} \cdot \omega(X_H \cdot
\gamma, \; w) - \tau_{\mathcal{U}}\cdot i^*_{\mathcal{M}} \cdot
\omega(X_H \cdot \gamma, \; T(\gamma \cdot \pi_Q) \cdot w)\\
& \;\;\;\;\;\;
- \tau_{\mathcal{U}}\cdot i^*_{\mathcal{M}} \cdot\mathbf{d}\gamma(T\pi_{Q}(X_H\cdot \gamma), \; T\pi_{Q}(w))\\
& =\pi^*_{/G}\cdot \omega_{\bar{\mathcal{K}}}\cdot\tau_{\bar{\mathcal{K}}}(X_H \cdot \gamma, \;
w) - \pi^*_{/G}\cdot \omega_{\bar{\mathcal{K}}}\cdot\tau_{\bar{\mathcal{K}}}(X_H \cdot \gamma,
\; T(\gamma \cdot \pi_Q) \cdot w)\\
& \;\;\;\;\;\; - \tau_{\mathcal{U}}\cdot i^*_{\mathcal{M}}
\cdot\mathbf{d}\gamma(T\pi_{Q}(X_H\cdot \gamma), \; T\pi_{Q}(w))\\
& = \omega_{\bar{\mathcal{K}}}(\tau_{\bar{\mathcal{K}}}\cdot T\pi_{/G}(X_H \cdot \gamma), \;
\tau_{\bar{\mathcal{K}}}\cdot T\pi_{/G} \cdot w) - \omega_{\bar{\mathcal{K}}}(\tau_{\bar{\mathcal{K}}}\cdot T\pi_{/G}(X_H \cdot
\gamma), \; \tau_{\bar{\mathcal{K}}}\cdot T(\pi_{/G} \cdot\gamma) \cdot T\pi_{Q}(w))\\ & \;\;\;\;\;\; -
\tau_{\mathcal{U}}\cdot i^*_{\mathcal{M}}
\cdot\mathbf{d}\gamma(T\pi_{Q}(X_H\cdot \gamma), \; T\pi_{Q}(w))\\
& = \omega_{\bar{\mathcal{K}}}(\tau_{\bar{\mathcal{K}}}\cdot T\pi_{/G}(X_H)\cdot
\pi_{/G}(\gamma), \; \tau_{\bar{\mathcal{K}}}\cdot T\pi_{/G} \cdot w) -
\omega_{\bar{\mathcal{K}}}(\tau_{\bar{\mathcal{K}}}\cdot T\pi_{/G}(X_H)\cdot \pi_{/G}(\gamma), \;
\tau_{\bar{\mathcal{K}}}\cdot T\bar{\gamma} \cdot T\pi_{Q}(w))\\ & \;\;\;\;\;\; - \tau_{\mathcal{U}}\cdot
i^*_{\mathcal{M}}
\cdot\mathbf{d}\gamma(T\pi_{Q}(X_H\cdot \gamma), \; T\pi_{Q}(w))\\
& = \omega_{\bar{\mathcal{K}}}(X_{\bar{\mathcal{K}}} \cdot
\bar{\gamma}, \; \tau_{\bar{\mathcal{K}}}\cdot T\pi_{/G} \cdot w)-
\omega_{\bar{\mathcal{K}}}(X_{\bar{\mathcal{K}}} \cdot
\bar{\gamma}, \; T\bar{\gamma} \cdot T\pi_{Q}(w))- \tau_{\mathcal{U}}\cdot
i^*_{\mathcal{M}} \cdot\mathbf{d}\gamma(T\pi_{Q}(X_H\cdot \gamma),
\; T\pi_{Q}(w)),
\end{align*}
where we have used that $\tau_{\bar{\mathcal{K}}}\cdot T\pi_{/G}(X_H)\cdot \bar{\gamma}
=\tau_{\bar{\mathcal{K}}}\cdot X_{h_{\bar{\mathcal{K}}}}\cdot \bar{\gamma}=
X_{\bar{\mathcal{K}}}\cdot \bar{\gamma}, $ and $\tau_{\bar{\mathcal{K}}}\cdot T\bar{\gamma}=T\bar{\gamma}, $ since
$\textmd{Im}(T\bar{\gamma})\subset \bar{\mathcal{K}}. $
If the one-form $\gamma: Q \rightarrow T^*Q $ is closed on $\mathcal{D}$ with respect to
$T\pi_Q: TT^* Q \rightarrow TQ, $ then we have that
$\mathbf{d}\gamma(T\pi_{Q}(X_H\cdot \gamma), \; T\pi_{Q}(w))=0, $
since $X_{H}\cdot \gamma, \; w \in \mathcal{F},$ and
$T\pi_{Q}(X_H\cdot \gamma), \; T\pi_{Q}(w) \in \mathcal{D}, $ and hence
$$
\tau_{\mathcal{U}}\cdot
i_{\mathcal{M}}^* \cdot\mathbf{d}\gamma(T\pi_{Q}(X_H\cdot \gamma),
\; T\pi_{Q}(w))=0,
$$
and
\begin{equation}
\omega_{\bar{\mathcal{K}}}(T\bar{\gamma} \cdot X_H^\gamma, \;
\tau_{\bar{\mathcal{K}}}\cdot T\pi_{/G} \cdot w)- \omega_{\bar{\mathcal{K}}}(X_{\bar{\mathcal{K}}} \cdot
\bar{\gamma}, \; \tau_{\bar{\mathcal{K}}}\cdot T\pi_{/G} \cdot w)
= -\omega_{\bar{\mathcal{K}}}(X_{\bar{\mathcal{K}}} \cdot
\bar{\gamma}, \; T\bar{\gamma} \cdot T\pi_{Q}(w)).
\label{4.1} \end{equation}
If $\bar{\gamma}$ satisfies the equation $
T\bar{\gamma}\cdot X_H^ \gamma = X_{\bar{\mathcal{K}}}\cdot
\bar{\gamma} ,$
from Lemma 3.4(i) we know that the right side of (4.1) becomes
\begin{align*}
 -\omega_{\bar{\mathcal{K}}}(X_{\bar{\mathcal{K}}} \cdot
\bar{\gamma}, \; T\bar{\gamma} \cdot T\pi_{Q}(w))
& = -\omega_{\bar{\mathcal{K}}}\cdot\tau_{\bar{\mathcal{K}}}(T\bar{\gamma}\cdot X_H^\gamma, \; T\bar{\gamma} \cdot T\pi_{Q}(w))\\
& = -\bar{\gamma}^*\omega_{\bar{\mathcal{K}}}\cdot\tau_{\bar{\mathcal{K}}}(T\pi_{Q} \cdot X_{H} \cdot
\gamma, \; T\pi_{Q}(w))\\
& = - \gamma^* \cdot \pi^*_{/G}\cdot \omega_{\bar{\mathcal{K}}}\cdot\tau_{\bar{\mathcal{K}}}(T\pi_{Q} \cdot X_{H} \cdot
\gamma, \; T\pi_{Q}(w))\\
& = - \gamma^* \cdot \tau_{\mathcal{U}}\cdot
i_{\mathcal{M}}^* \cdot \omega(T\pi_{Q}(X_{H}\cdot\gamma), \; T\pi_{Q}(w))\\
& = -\tau_{\mathcal{U}}\cdot
i_{\mathcal{M}}^* \cdot\gamma^*\omega( T\pi_{Q}(X_{H}\cdot\gamma), \; T\pi_{Q}(w))\\
& = \tau_{\mathcal{U}}\cdot i_{\mathcal{M}}^* \cdot
\mathbf{d}\gamma(T\pi_{Q}( X_{H}\cdot\gamma ), \; T\pi_{Q}(w))=0,
\end{align*}
where $\gamma^*\cdot \tau_{\mathcal{U}}\cdot i^*_{\mathcal{M}}
\cdot \omega= \tau_{\mathcal{U}}\cdot i^*_{\mathcal{M}}
\cdot\gamma^*\cdot \omega, $ because $\textmd{Im}(\gamma)\subset
\mathcal{M}. $
But, since the nonholonomic reduced distributional two-form
$\omega_{\bar{\mathcal{K}}}$ is non-degenerate,
the left side of (4.1) equals zero, only when
$\bar{\gamma}$ satisfies the equation $
T\bar{\gamma}\cdot X_H^ \gamma = X_{\bar{\mathcal{K}}}\cdot
\bar{\gamma} .$ Thus,
if the one-form $\gamma: Q \rightarrow T^*Q $ is closed on $\mathcal{D}$ with respect to
$T\pi_Q: TT^* Q \rightarrow TQ, $ then $\bar{\gamma}$ must be a solution of the Type I of Hamilton-Jacobi equation
$T\bar{\gamma}\cdot X_H^ \gamma = X_{\bar{\mathcal{K}}}\cdot
\bar{\gamma}. $
\hskip 0.3cm $\blacksquare$\\

Next, for any $G$-invariant symplectic map $\varepsilon: T^* Q \rightarrow T^* Q $, we can prove
the following Type II of
Hamilton-Jacobi theorem for the nonholonomic reduced distributional Hamiltonian system.
For convenience, the maps involved in the following theorem and its
proof are shown in Diagram-4.
\begin{center}
\hskip 0cm \xymatrix{ & \mathcal{M} \ar[d]_{X_{\mathcal{K}}}
\ar[r]^{i_{\mathcal{M}}} & T^* Q \ar[d]_{X_{H\cdot \varepsilon}}
\ar[dr]^{X_H^\varepsilon} \ar[r]^{\pi_Q}
  & Q \ar[r]^{\gamma} & T^* Q \ar[d]_{X_H} \ar[r]^{\pi_{/G}} & T^* Q/G \ar[d]_{X_{h}}
  & \mathcal{\bar{M}} \ar[l]_{i_{\mathcal{\bar{M}}}} \ar[d]_{X_{\mathcal{\bar{K}}}}\\
  & \mathcal{K}
  & T(T^*Q) \ar[l]_{\tau_{\mathcal{K}}}
  & TQ \ar[l]_{T\gamma}
  & T(T^* Q) \ar[l]_{T\pi_Q} \ar[r]^{T\pi_{/G}}
  & T(T^* Q/G) \ar[r]^{\tau_{\mathcal{\bar{K}}}} & \mathcal{\bar{K}} }
\end{center}
$$\mbox{Diagram-4}$$

\begin{theo} (Type II of Hamilton-Jacobi Theorem for a Nonholonomic
Reduced Distributional Hamiltonian System)
For a given nonholonomic reducible Hamiltonian system
$(T^*Q,G,\omega,\mathcal{D},H)$ with an associated
nonholonomic reduced distributional Hamiltonian system
$(\bar{\mathcal{K}},\omega_{\bar{\mathcal{K}}},h_{\bar{\mathcal{K}}})$, assume that
$\gamma: Q \rightarrow T^*Q$ is an one-form on $Q$, and $\lambda=
\gamma \cdot \pi_{Q}: T^* Q \rightarrow T^* Q, $ and for any $G$-invariant
symplectic map $\varepsilon: T^* Q \rightarrow T^* Q $, denote by
$X_H^\varepsilon = T\pi_{Q}\cdot X_H \cdot \varepsilon$, where $X_{H}$ is
the dynamical vector field of the corresponding unconstrained
Hamiltonian system with symmetry $(T^*Q,G,\omega,H)$. Moreover,
assume that $\textmd{Im}(\gamma)\subset \mathcal{M}, $ and it is
$G$-invariant, $ \textmd{Im}(T\gamma)\subset \mathcal{K}, $ and
$\bar{\gamma}=\pi_{/G}(\gamma): Q \rightarrow T^* Q/G $, and
$\bar{\lambda}=\pi_{/G}(\lambda): T^* Q \rightarrow T^* Q/G, $ and
$\bar{\varepsilon}=\pi_{/G}(\varepsilon): T^* Q \rightarrow T^* Q/G. $ Then
$\varepsilon$ and $\bar{\varepsilon}$ satisfy the equation
$\tau_{\bar{\mathcal{K}}}\cdot T\bar{\varepsilon}\cdot X_{h_{\bar{\mathcal{K}}}\cdot
\bar{\varepsilon}}= T\bar{\lambda} \cdot X_H\cdot \varepsilon, $ if and only if they satisfy the
equation $T\bar{\gamma}\cdot X_H^ \varepsilon =
X_{\bar{\mathcal{K}}}\cdot \bar{\varepsilon}. $ Here
$ X_{h_{\bar{\mathcal{K}}} \cdot\bar{\varepsilon}}$ is the Hamiltonian
vector field of the function $h_{\bar{\mathcal{K}}}\cdot \bar{\varepsilon}: T^* Q\rightarrow
\mathbb{R}, $ and $X_{\bar{\mathcal{K}}}$ is the dynamical vector
field of the nonholonomic reduced distributional Hamiltonian system
$(\bar{\mathcal{K}},\omega_{\bar{\mathcal{K}}},h_{\bar{\mathcal{K}}})$. The equation $
T\bar{\gamma}\cdot X_H^\varepsilon = X_{\bar{\mathcal{K}}}\cdot
\bar{\varepsilon},$ is called the Type II of Hamilton-Jacobi equation for the
nonholonomic reduced distributional Hamiltonian system
$(\bar{\mathcal{K}},\omega_{\bar{\mathcal{K}}},h_{\bar{\mathcal{K}}})$.
\end{theo}

\noindent{\bf Proof: } In the same way, we note that
$\textmd{Im}(\gamma)\subset \mathcal{M}, $ and it is $G$-invariant,
$ \textmd{Im}(T\gamma)\subset \mathcal{K}, $ and hence
$\textmd{Im}(T\bar{\gamma})\subset \bar{\mathcal{K}}, $ in this case,
$\pi^*_{/G}\cdot\omega_{\bar{\mathcal{K}}}\cdot \tau_{\bar{\mathcal{K}}}= \tau_{\mathcal{U}}\cdot
\omega_{\mathcal{M}}= \tau_{\mathcal{U}}\cdot i_{\mathcal{M}}^*\cdot
\omega, $ along $\textmd{Im}(T\bar{\gamma})$. Thus, using the
non-degenerate and nonholonomic reduced distributional two-form $\omega_{\bar{\mathcal{K}}}$,
from Lemma 3.4, if we take that $v=\tau_{\bar{\mathcal{K}}}\cdot T\pi_{/G} \cdot X_{H}\cdot \varepsilon
=X_{\bar{\mathcal{K}}}\cdot\bar{\varepsilon} \in \bar{\mathcal{K}},
$ and for any $w \in \mathcal{F}, \; T\lambda(w)\neq 0, $ and
$\tau_{\bar{\mathcal{K}}}\cdot T\pi_{/G}\cdot w \neq 0, $ then we have that
\begin{align*}
& \omega_{\bar{\mathcal{K}}}(T\bar{\gamma} \cdot X_H^\varepsilon, \;
\tau_{\bar{\mathcal{K}}}\cdot T\pi_{/G} \cdot w) = \omega_{\bar{\mathcal{K}}}(\tau_{\bar{\mathcal{K}}}\cdot T(\pi_{/G} \cdot
\gamma) \cdot X_H^\varepsilon, \; \tau_{\bar{\mathcal{K}}}\cdot T\pi_{/G} \cdot w )\\
& = \pi^*_{/G}\cdot \omega_{\bar{\mathcal{K}}}
\cdot\tau_{\bar{\mathcal{K}}}(T\gamma \cdot X_H^\varepsilon, \; w) =
\tau_{\mathcal{U}}\cdot i^*_{\mathcal{M}} \cdot\omega(T\gamma \cdot
X_H^\varepsilon, \; w)\\
& = \tau_{\mathcal{U}}\cdot i^*_{\mathcal{M}} \cdot
\omega(T(\gamma
\cdot \pi_Q)\cdot X_H \cdot \varepsilon, \; w) \\
& = \tau_{\mathcal{U}}\cdot i^*_{\mathcal{M}} \cdot
(\omega(X_H \cdot \varepsilon, \; w-T(\gamma \cdot \pi_Q)\cdot w)
- \mathbf{d}\gamma(T\pi_{Q}(X_H\cdot \varepsilon), \; T\pi_{Q}(w)))\\
& = \tau_{\mathcal{U}}\cdot i^*_{\mathcal{M}} \cdot \omega(X_H \cdot
\varepsilon, \; w) - \tau_{\mathcal{U}}\cdot i^*_{\mathcal{M}} \cdot
\omega(X_H \cdot \varepsilon, \; T\lambda \cdot w)
- \tau_{\mathcal{U}}\cdot i^*_{\mathcal{M}} \cdot\mathbf{d}\gamma(T\pi_{Q}(X_H\cdot \varepsilon), \; T\pi_{Q}(w))\\
& =\pi^*_{/G}\cdot \omega_{\bar{\mathcal{K}}}\cdot \tau_{\bar{\mathcal{K}}}(X_H \cdot \varepsilon, \;
w) - \pi^*_{/G}\cdot \omega_{\bar{\mathcal{K}}}\cdot \tau_{\bar{\mathcal{K}}}(X_H \cdot \varepsilon,
\; T\lambda \cdot w)+ \tau_{\mathcal{U}}\cdot i^*_{\mathcal{M}}
\cdot \lambda^* \omega(X_H\cdot \varepsilon, \; w)\\
& = \omega_{\bar{\mathcal{K}}}(\tau_{\bar{\mathcal{K}}}\cdot T\pi_{/G}(X_H \cdot \varepsilon), \;
\tau_{\bar{\mathcal{K}}}\cdot T\pi_{/G} \cdot w) - \omega_{\bar{\mathcal{K}}}(\tau_{\bar{\mathcal{K}}}\cdot T\pi_{/G}(X_H \cdot
\varepsilon), \; \tau_{\bar{\mathcal{K}}}\cdot T(\pi_{/G} \cdot\lambda) \cdot w)\\
& \;\;\;\;\;\; +\pi^*_{/G}\cdot \omega_{\bar{\mathcal{K}}}\cdot \tau_{\bar{\mathcal{K}}}(T\lambda\cdot X_H \cdot \varepsilon, \;
T\lambda \cdot w)\\
& = \omega_{\bar{\mathcal{K}}}(\tau_{\bar{\mathcal{K}}}\cdot T\pi_{/G}(X_H)\cdot
\pi_{/G}(\varepsilon), \; \tau_{\bar{\mathcal{K}}}\cdot T\pi_{/G} \cdot w) -
\omega_{\bar{\mathcal{K}}}(\tau_{\bar{\mathcal{K}}}\cdot T\pi_{/G}(X_H)\cdot \pi_{/G}(\varepsilon), \;
\tau_{\bar{\mathcal{K}}}\cdot T\bar{\lambda} \cdot w)\\
& \;\;\;\;\;\; +\omega_{\bar{\mathcal{K}}}(\tau_{\bar{\mathcal{K}}}\cdot T\pi_{/G}\cdot T\lambda\cdot X_H \cdot \varepsilon,
\; \tau_{\bar{\mathcal{K}}}\cdot T\pi_{/G}\cdot T\lambda \cdot w)\\
& = \omega_{\bar{\mathcal{K}}}(X_{\bar{\mathcal{K}}} \cdot
\bar{\varepsilon}, \; \tau_{\bar{\mathcal{K}}}\cdot T\pi_{/G} \cdot w)-
\omega_{\bar{\mathcal{K}}}(\tau_{\bar{\mathcal{K}}} \cdot X_{h_{\bar{\mathcal{K}}}}\cdot
\bar{\varepsilon}, \; T\bar{\lambda} \cdot w)+ \omega_{\bar{\mathcal{K}}}(T\bar{\lambda}\cdot X_H \cdot \varepsilon,
\; T\bar{\lambda} \cdot w),
\end{align*}
where we have used that $\tau_{\bar{\mathcal{K}}}\cdot T\pi_{/G}(X_H)\cdot \bar{\varepsilon}
=\tau_{\bar{\mathcal{K}}}(X_{h_{\bar{\mathcal{K}}}})\cdot \bar{\varepsilon}=
X_{\bar{\mathcal{K}}}\cdot \bar{\varepsilon}, $
and $\tau_{\bar{\mathcal{K}}}\cdot T\pi_{/G}\cdot T\lambda=T\bar{\lambda}, $ since
$\textmd{Im}(T\bar{\gamma})\subset \bar{\mathcal{K}}. $ Note that
$\varepsilon: T^* Q \rightarrow T^* Q $ is symplectic, and
$\bar{\varepsilon}^*= \varepsilon^*\cdot \pi_{/G}^*: T^*(T^* Q)/G
\rightarrow T^*T^* Q$ is also symplectic along $\bar{\varepsilon}$, and
hence $X_{h_{\bar{\mathcal{K}}}}\cdot \bar{\varepsilon}
= T\bar{\varepsilon} \cdot X_{h_{\bar{\mathcal{K}}} \cdot
\bar{\varepsilon}}, $ along $\bar{\varepsilon}$, and hence
$\tau_{\bar{\mathcal{K}}}\cdot X_{h_{\bar{\mathcal{K}}}} \cdot\bar{\varepsilon}
= \tau_{\bar{\mathcal{K}}}\cdot T\bar{\varepsilon} \cdot X_{h_{\bar{\mathcal{K}}}
\cdot \bar{\varepsilon}}, $ along $\bar{\varepsilon}$.
Then we have that
\begin{align*}
& \omega_{\bar{\mathcal{K}}}(T\bar{\gamma} \cdot X_H^\varepsilon, \;
\tau_{\bar{\mathcal{K}}}\cdot T\pi_{/G}\cdot w)-
\omega_{\bar{\mathcal{K}}}(X_{\bar{\mathcal{K}}}\cdot \bar{\varepsilon},
\; \tau_{\bar{\mathcal{K}}}\cdot T\pi_{/G} \cdot w) \nonumber \\
& = -\omega_{\bar{\mathcal{K}}}(\tau_{\bar{\mathcal{K}}} \cdot X_{h_{\bar{\mathcal{K}}}} \cdot
\bar{\varepsilon}, \; T\bar{\lambda} \cdot w)+ \omega_{\bar{\mathcal{K}}}(T\bar{\lambda}\cdot X_H \cdot \varepsilon,
\; T\bar{\lambda} \cdot w)\\
& = \omega_{\bar{\mathcal{K}}}(T\bar{\lambda}\cdot X_H \cdot \varepsilon- \tau_{\bar{\mathcal{K}}}\cdot T\bar{\varepsilon} \cdot X_{h_{\bar{\mathcal{K}}}\cdot \bar{\varepsilon}},
\; T\bar{\lambda} \cdot w).
\end{align*}
Because the nonholonomic reduced distributional two-form
$\omega_{\bar{\mathcal{K}}}$ is non-degenerate, it follows that the equation
$T\bar{\gamma}\cdot X_H^\varepsilon = X_{\bar{\mathcal{K}}}\cdot
\bar{\varepsilon},$ is equivalent to the equation $T\bar{\lambda}\cdot X_H \cdot \varepsilon
= \tau_{\bar{\mathcal{K}}}\cdot T\bar{\varepsilon} \cdot X_{h_{\bar{\mathcal{K}}}
\cdot \bar{\varepsilon}}. $
Thus, $\varepsilon$ and $\bar{\varepsilon}$ satisfy the equation
$T\bar{\lambda}\cdot X_H \cdot \varepsilon
= \tau_{\bar{\mathcal{K}}}\cdot T\bar{\varepsilon} \cdot X_{h_{\bar{\mathcal{K}}}
\cdot \bar{\varepsilon}}, $ if and only if they satisfy
the Type II of Hamilton-Jacobi equation $T\bar{\gamma}\cdot X_H^\varepsilon = X_{\bar{\mathcal{K}}}\cdot
\bar{\varepsilon} .$
\hskip 0.3cm $\blacksquare$\\

For a given nonholonomic reducible Hamiltonian system
$(T^*Q,G,\omega,\mathcal{D},H)$ with an associated
nonholonomic reduced distributional Hamiltonian system
$(\bar{\mathcal{K}},\omega_{\bar{\mathcal{K}}},h_{\bar{\mathcal{K}}})$,
we know that the nonholonomic dynamical vector field
$X_{\mathcal{K}}$ and the nonholonomic reduced dynamical
vector field $X_{\bar{\mathcal{K}}}$ are $\pi_{/G}$-related,
that is, $X_{\bar{\mathcal{K}}}\cdot \pi_{/G}=T\pi_{/G}\cdot
X_{\mathcal{K}}. $ Then we can prove the following Theorem 4.4,
which states the relationship between the solutions of Type II of
Hamilton-Jacobi equations and nonholonomic reduction.

\begin{theo}
For a given nonholonomic reducible Hamiltonian system
$(T^*Q,G,\omega,\mathcal{D},H)$ with an associated
nonholonomic reduced distributional Hamiltonian system
$(\bar{\mathcal{K}},\omega_{\bar{\mathcal{K}}},h_{\bar{\mathcal{K}}})$, assume that
$\gamma: Q \rightarrow T^*Q$ is an one-form on $Q$, and
$\lambda=\gamma \cdot \pi_{Q}: T^* Q \rightarrow T^* Q, $
and $\varepsilon: T^* Q \rightarrow T^* Q $ is a $G$-invariant symplectic map.
Moreover, assume that $\textmd{Im}(\gamma)\subset
\mathcal{M}, $ and it is $G$-invariant, $
\textmd{Im}(T\gamma)\subset \mathcal{K}, $ and
$\bar{\gamma}=\pi_{/G}(\gamma): Q \rightarrow T^* Q/G $, and
$\bar{\lambda}=\pi_{/G}(\lambda): T^* Q \rightarrow T^* Q/G, $ and
$\bar{\varepsilon}=\pi_{/G}(\varepsilon): T^* Q \rightarrow T^* Q/G. $ Then $\varepsilon$
is a solution of the Type II of Hamilton-Jacobi equation, $T\gamma\cdot
X_H^\varepsilon= X_{\mathcal{K}}\cdot \varepsilon, $ for the distributional
Hamiltonian system $(\mathcal{K},\omega_{\mathcal{K}},H_{\mathcal{K}})$, if and
only if $\varepsilon$ and $\bar{\varepsilon}$ satisfy the Type II of
Hamilton-Jacobi equation $T\bar{\gamma}\cdot X_H^\varepsilon =
X_{\bar{\mathcal{K}}}\cdot \bar{\varepsilon}, $ for the nonholonomic reduced
distributional Hamiltonian system $ (\bar{\mathcal{K}},
\omega_{\bar{\mathcal{K}}}, h_{\bar{\mathcal{K}}} ). $
\end{theo}

\noindent{\bf Proof: } Note that
$\textmd{Im}(\gamma)\subset \mathcal{M},$ and
it is $G$-invariant, $\textmd{Im}(T\gamma)\subset \mathcal{K}, $
and hence $\textmd{Im}(T\bar{\gamma})\subset \bar{\mathcal{K}}, $ in
this case, $\pi^*_{/G}\cdot\omega_{\bar{\mathcal{K}}}\cdot \tau_{\bar{\mathcal{K}}}= \tau_{\mathcal{U}}\cdot
\omega_{\mathcal{M}}= \tau_{\mathcal{U}}\cdot i_{\mathcal{M}}^*\cdot
\omega, $ along $\textmd{Im}(T\bar{\gamma})$, and
$\tau_{\bar{\mathcal{K}}}\cdot T\bar{\gamma}= T\bar{\gamma},
\; \tau_{\bar{\mathcal{K}}}\cdot X_{\bar{\mathcal{K}}}= X_{\bar{\mathcal{K}}}. $
Since nonholonomic vector field $X_{\mathcal{K}}$ and the
vector field $X_{\bar{\mathcal{K}}}$ are $\pi_{/G}$-related,
that is, $X_{\bar{\mathcal{K}}}\cdot \pi_{/G}=T\pi_{/G}\cdot
X_{\mathcal{K}}, $ using the non-degenerate and nonholonomic reduced
distributional two-form $\omega_{\bar{\mathcal{K}}}$,
we have that
\begin{align*}
& \omega_{\bar{\mathcal{K}}}(T\bar{\gamma} \cdot X_H^\varepsilon
- X_{\bar{\mathcal{K}}}\cdot \bar{\varepsilon}, \; \tau_{\bar{\mathcal{K}}}\cdot T\pi_{/G}\cdot w)\\
& = \omega_{\bar{\mathcal{K}}}(T\bar{\gamma} \cdot X_H^\varepsilon, \;
\tau_{\bar{\mathcal{K}}}\cdot T\pi_{/G}\cdot w)-
\omega_{\bar{\mathcal{K}}}(X_{\bar{\mathcal{K}}}\cdot \bar{\varepsilon},
\; \tau_{\bar{\mathcal{K}}}\cdot T\pi_{/G} \cdot w) \\
& = \omega_{\bar{\mathcal{K}}}(\tau_{\bar{\mathcal{K}}}\cdot T\bar{\gamma}\cdot X_H^
\varepsilon, \; \tau_{\bar{\mathcal{K}}}\cdot T\pi_{/G}\cdot w)
-\omega_{\bar{\mathcal{K}}}(\tau_{\bar{\mathcal{K}}}\cdot X_{\bar{\mathcal{K}}}\cdot \pi_{/G}\cdot \varepsilon,
\; \tau_{\bar{\mathcal{K}}}\cdot T\pi_{/G}\cdot w)\\
& = \omega_{\bar{\mathcal{K}}}\cdot \tau_{\bar{\mathcal{K}}}(T\pi_{/G}\cdot T\gamma \cdot X_H^
\varepsilon, \; T\pi_{/G} \cdot w)
- \omega_{\bar{\mathcal{K}}}\cdot \tau_{\bar{\mathcal{K}}}(T\pi_{/G}\cdot
X_{\mathcal{K}}\cdot
\varepsilon, \; T\pi_{/G}\cdot w)\\
& = \pi^*_{/G}\cdot\omega_{\bar{\mathcal{K}}}\cdot \tau_{\bar{\mathcal{K}}}(T\gamma \cdot X_H^
\varepsilon, \; w)
- \pi^*_{/G}\cdot\omega_{\bar{\mathcal{K}}}\cdot \tau_{\bar{\mathcal{K}}}(X_{\mathcal{K}} \cdot
\varepsilon, \; w)\\
& = \tau_{\mathcal{U}}\cdot i_{\mathcal{M}}^* \cdot
\omega(T\gamma \cdot X_H^
\varepsilon, \; w)- \tau_{\mathcal{U}}\cdot i_{\mathcal{M}}^* \cdot
\omega(X_{\mathcal{K}} \cdot
\varepsilon, \; w).
\end{align*}
In the case we considered that $\tau_{\mathcal{U}}\cdot i_{\mathcal{M}}^* \cdot
\omega=\tau_{\mathcal{K}}\cdot i_{\mathcal{M}}^* \cdot
\omega= \omega_{\mathcal{K}}\cdot \tau_{\mathcal{K}}, $
and
$\tau_{\mathcal{K}}\cdot T\gamma =T\gamma, \; \tau_{\mathcal{K}} \cdot X_{\mathcal{K}}= X_{\mathcal{K}}$,
since $\textmd{Im}(\gamma)\subset
\mathcal{M}, $ and $\textmd{Im}(T\gamma)\subset \mathcal{K}. $
Thus, we have that
\begin{align*}
& \omega_{\bar{\mathcal{K}}}(T\bar{\gamma} \cdot X_H^\varepsilon
- X_{\bar{\mathcal{K}}}\cdot \bar{\varepsilon}, \; \tau_{\bar{\mathcal{K}}}\cdot T\pi_{/G}\cdot w)\\
& = \omega_{\mathcal{K}}\cdot \tau_{\mathcal{K}}(T\gamma \cdot X_H^
\varepsilon, \; w)- \omega_{\mathcal{K}}\cdot \tau_{\mathcal{K}}(X_{\mathcal{K}} \cdot \varepsilon, \; w)\\
& = \omega_{\mathcal{K}}(\tau_{\mathcal{K}} \cdot T\gamma \cdot X_H^
\varepsilon, \; \tau_{\mathcal{K}} \cdot w)
- \omega_{\mathcal{K}}(\tau_{\mathcal{K}} \cdot X_{\mathcal{K}}\cdot \varepsilon, \; \tau_{\mathcal{K}} \cdot w)\\
& = \omega_{\mathcal{K}}(T\gamma \cdot X_H^
\varepsilon- X_{\mathcal{K}}\cdot \varepsilon, \; \tau_{\mathcal{K}} \cdot w).
\end{align*}
Because the distributional two-form $\omega_{\mathcal{K}}$
and the nonholonomic reduced distributional
two-form $\omega_{\bar{\mathcal{K}}}$ are both non-degenerate,
it follows that the equation
$T\bar{\gamma}\cdot X_H^\varepsilon=
X_{\bar{\mathcal{K}}}\cdot \bar{\varepsilon}, $
is equivalent to the equation $T\gamma\cdot X_H^\varepsilon= X_{\mathcal{K}}\cdot \varepsilon. $
Thus, $\varepsilon$ is a solution of the Type II of Hamilton-Jacobi equation
$T\gamma\cdot X_H^\varepsilon= X_{\mathcal{K}}\cdot \varepsilon, $ for the distributional
Hamiltonian system $(\mathcal{K},\omega_{\mathcal {K}},H_{\mathcal{K}})$, if and only if
$\varepsilon$ and $\bar{\varepsilon} $ satisfy the Type II of Hamilton-Jacobi
equation $T\bar{\gamma}\cdot X_H^\varepsilon=
X_{\bar{\mathcal{K}}}\cdot \bar{\varepsilon}, $ for the
nonholonomic reduced distributional Hamiltonian system
$(\bar{\mathcal{K}},\omega_{\bar{\mathcal{K}}},h_{\bar{\mathcal{K}}})$.
\hskip 0.3cm
$\blacksquare$

\begin{rema}
It is worthy of noting that the
formulations of Type I and Type II of Hamilton-Jacobi equation for a
nonholonomic reduced distributional Hamiltonian system, given by
Theorem $4.2$ and Theorem $4.3$, have more extensive sense, because
in general, the one-form $\gamma$ is not
given by a generating function of a symplectic map.
When $\gamma$ is a solution of the classical
Hamilton-Jacobi equation, that is, $X_H\cdot \gamma=0, $ then
$X_H^\gamma= T\pi_Q\cdot X_H\cdot \gamma=0,$ and hence from the Type
I of Hamilton-Jacobi equation, we have that
$X_{\bar{\mathcal{K}}}\cdot \bar{\gamma}= T\bar{\gamma}\cdot
X_H^\gamma=0, $ which shows that the dynamical vector field
of the nonholonomic reduced distributional Hamiltonian system
$(\bar{\mathcal{K}},\omega_{\bar{\mathcal{K}}},h_{\bar{\mathcal{K}}})$ is degenerate
along $\bar{\gamma}$. The equation $X_{\bar{\mathcal{K}}}\cdot
\bar{\gamma}=0$ is called the classical Hamilton-Jacobi equation for
the nonholonomic reduced distributional Hamiltonian system
$(\bar{\mathcal{K}},\omega_{\bar{\mathcal{K}}},h_{\bar{\mathcal{K}}}). $ In addition,
for a symplectic map $\varepsilon: T^* Q \rightarrow T^* Q $, if
$X_H\cdot \varepsilon=0, $ then from the Type II of Hamilton-Jacobi
equation, we have that $X_{\bar{\mathcal{K}}}\cdot
\bar{\varepsilon}= T\bar{\gamma}\cdot X_H^\varepsilon=0. $ But, from
the equation $T\bar{\lambda}\cdot X_H \cdot \varepsilon =
\tau_{\bar{\mathcal{K}}}\cdot T\bar{\varepsilon} \cdot X_{h_{\bar{\mathcal{K}}} \cdot
\bar{\varepsilon}}, $ we know that the equation
$X_{\bar{\mathcal{K}}}\cdot \bar{\varepsilon}=0 $ is not equivalent
to the equation $X_{h_{\bar{\mathcal{K}}} \cdot \bar{\varepsilon}}=0.$
\end{rema}

\section{Nonholonomic Hamiltonian System with Symmetry and Momentum Map}

As it is well known that momentum map is a very important notion in modern study
of geometric mechanics, and it is a geometric generalization of the classical
linear and angular momentum. A fundamental fact about momentum map is that
if the Hamiltonian $H$ is invariant under the action of a Lie group $G$,
then the vector valued function $\mathbf{J}$ is a constant of the motion
for the dynamics of the Hamiltonian vector field $X_H$ associated to $H$,
that is, all momentum maps are conserved quantities. Moreover, momentum map
has infinitesimal equivariance, such that it plays an important role in the
study of reduction theory of Hamiltonian systems with symmetries, see Marsden
\cite{ma92} and Marsden et al.\cite{mamiorpera07, mamora90}.
Now, it is a natural problem what and how we could do, when the Hamiltonian system
we considered has nonholonomic constrains, and the Lie group $G$ is
not Abelian, and $G_\mu\neq G ,$ where $G_\mu$ is the isotropy subgroup
of coadjoint $G$-action at the point $\mu\in \mathfrak{g}^*$,
and hence the above procedure of nonholonomic reduction
given in $\S 4$ does not work or is not efficient enough.
In this section, we shall consider a nonholonomic Hamiltonian
system with symmetry and momentum map, and give two types of Hamilton-Jacobi
theorems of the nonholonomic point and orbit reduced distributional Hamiltonian
systems with respect to momentum map.\\

\subsection{Hamilton-Jacobi equations in the case compatible with Marsden-Weinstein reduction}

In this subsection, for a nonholonomic Hamiltonian
system with symmetry and momentum map
$(T^*Q,G,\omega,\mathbf{J},\mathcal{D},H)$,
where $\omega$ is the canonical
symplectic form on $T^* Q$, and $\mathcal{D}\subset TQ$ is a
$\mathcal{D}$-completely and $\mathcal{D}$-regularly nonholonomic
constraint of the system, and $\mathcal{D}$ and $H$ are both
$G$-invariant, we first give the $\mathbf{J}$-nonholonomic regular point reduction of the system
compatible with Marsden-Weinstein reduction,
and a $\mathbf{J}$-nonholonomic $R_p$-reduced distribution $\mathcal{K}_\mu$, an associated
non-degenerate and $\mathbf{J}$-nonholonomic $R_p$-reduced
distributional two-form $\omega_{\mathcal{K}_\mu}$, which is
induced by the canonical symplectic form $\omega$ on $T^*Q$,
and a $\mathbf{J}$-nonholonomic $R_p$-reduced distributional Hamiltonian system,
where the "regular point reduced" is simply written as $R_p$-reduced.
Then we derive precisely
the geometric constraint conditions of the $\mathbf{J}$-nonholonomic
$R_p$-reduced distributional two-form
$\omega_{\mathcal{K}_\mu}$ for the nonholonomic reducible dynamical vector field,
that is, the two types of Hamilton-Jacobi equation for the
$\mathbf{J}$-nonholonomic $R_p$-reduced distributional Hamiltonian system, which are an
extension of the above two types of Hamilton-Jacobi equation for the distributional
Hamiltonian system under $\mathbf{J}$-nonholonomic regular point reduction.\\

At first, we need to give
carefully a geometric formulation of the $\mathbf{J}$-nonholonomic
$R_p$-reduced distributional Hamiltonian system, by using momentum map and the
nonholonomic reduction compatible with Marsden-Weinstein reduction.
Now, we assume that the
6-tuple $(T^*Q,G,\omega,\mathbf{J},\mathcal{D},H)$ is a
$\mathcal{D}$-completely and $\mathcal{D}$-regularly nonholonomic
Hamiltonian system with symmetry and momentum map, and the Lie group
$G$, which may not be Abelian, acts smoothly by the left on $Q$,
its tangent lifted action on
$TQ$ and its cotangent lifted action on $T^\ast Q$,
and we assume that the action on $T^\ast Q$ is free,
proper and symplectic, and admits an
$\operatorname{Ad}^\ast$-equivariant momentum map $\mathbf{J}:
T^\ast Q\rightarrow \mathfrak{g}^\ast$, where $\mathfrak{g}$ is a
Lie algebra of $G$ and $\mathfrak{g}^\ast$ is the dual of
$\mathfrak{g}$. Let $\mu \in\mathfrak{g}^\ast$ be a regular value of
$\mathbf{J}$ and denote by $G_\mu$ the isotropy subgroup of the
coadjoint $G$-action at the point $\mu\in\mathfrak{g}^\ast$, which
is defined by $G_\mu=\{g\in G|\operatorname{Ad}_g^\ast \mu=\mu \}$.
Since $G_\mu (\subset G)$ acts freely and properly on $Q$ and on
$T^\ast Q$, then $G_\mu$ acts also freely and properly on
$\mathbf{J}^{-1}(\mu)$, so that the space $(T^\ast
Q)_\mu=\mathbf{J}^{-1}(\mu)/G_\mu$ is a symplectic manifold with
symplectic form $\omega_\mu$ uniquely characterized by the relation
\begin{equation}\pi_\mu^\ast \omega_\mu=i_\mu^\ast
\omega. \label{5.1}\end{equation} The map
$i_\mu:\mathbf{J}^{-1}(\mu)\rightarrow T^\ast Q$ is the inclusion
and $\pi_\mu:\mathbf{J}^{-1}(\mu)\rightarrow (T^\ast Q)_\mu$ is the
projection. The pair $((T^\ast Q)_\mu,\omega_\mu)$ is the
Marsden-Weinstein reduced space of $(T^\ast Q,\omega)$ at $\mu$.\\

Assume that $H: T^*Q \rightarrow \mathbb{R}$ is a $G$-invariant
Hamiltonian, and the $\mathcal{D}$-completely and
$\mathcal{D}$-regularly nonholonomic constraints $\mathcal{D}\subset
TQ$ is a $G$-invariant distribution. From \S2, we know that, by using the
Legendre transformation $\mathcal{F}L: TQ \rightarrow T^*Q $, we can
define the constraint submanifold
$\mathcal{M}=\mathcal{F}L(\mathcal{D})\subset T^*Q $ and the
distribution $\mathcal{F}$ which is the pre-image of the
nonholonomic constraints $\mathcal{D}$ for the map $T\pi_Q: TT^* Q
\rightarrow TQ$, that is, $\mathcal{F}=(T\pi_Q)^{-1}(\mathcal{D})$,
and $\mathcal{K}=\mathcal{F} \cap T\mathcal{M}$. Moreover, we can
also define the distributional two-form $\omega_\mathcal{K}$, a vector
field $X_\mathcal{K}$ and the function $H_\mathcal{K}$, such that
$\mathbf{i}_{X_\mathcal{K}}\omega_\mathcal{K}=\mathbf{d}H_\mathcal{K}$.
Since $\mathcal{D}\subset TQ$ is a $G$-invariant distribution, and the
Legendre transformation $\mathcal{F}L: TQ \rightarrow T^*Q$ is a
fiber-preserving map, then
$\mathcal{M}=\mathcal{F}L(\mathcal{D})\subset T^*Q$ is
$G$-invariant. For a regular value $\mu\in\mathfrak{g}^\ast$ of the
momentum map $\mathbf{J}:T^\ast Q\rightarrow \mathfrak{g}^\ast$,
we shall assume that the constraint submanifold $\mathcal{M}$
is clean intersection with $\mathbf{J}^{-1}(\mu)$, that is,
$\mathcal{M} \cap \mathbf{J}^{-1}(\mu)\neq \emptyset$.
Note that $\mathcal{M}$ is also $G_\mu (\subset G)$ action
invariant, and so is $\mathbf{J}^{-1}(\mu)$, because $\mathbf{J}$ is
$\operatorname{Ad}^\ast$-equivariant. It follows that the quotient
space $\mathcal{M}_\mu =(\mathcal{M}\cap \mathbf{J}^{-1}(\mu))
/G_\mu \subset (T^\ast Q)_\mu$ of the $G_\mu$-orbit in
$\mathcal{M}\cap \mathbf{J}^{-1}(\mu)$, is a smooth manifold with
projection $\pi_\mu: \mathcal{M}\cap \mathbf{J}^{-1}(\mu)
\rightarrow \mathcal{M}_\mu$ which is a surjective submersion.
Denote $i_{\mathcal{M}_\mu}: \mathcal{M}_\mu\rightarrow
(T^*Q)_\mu, $ and $\omega_{\mathcal{M}_\mu}= i_{\mathcal{M}_\mu}^*
\omega_\mu $, that is, the symplectic form
$\omega_{\mathcal{M}_\mu}$ is induced from the Marsden-Weinstein reduced symplectic
form $\omega_\mu$ on $(T^* Q)_\mu$, where $i_{\mathcal{M}_\mu}^*:
T^*(T^*Q)_\mu \rightarrow T^*\mathcal{M}_\mu. $ Moreover, the
distribution $\mathcal{F}$ pushes down to a distribution
$\mathcal{F}_\mu= T\pi_\mu\cdot \mathcal{F}$ on $(T^\ast Q)_\mu$,
and we define $\mathcal{K}_\mu=\mathcal{F}_\mu \cap
T\mathcal{M}_\mu$. Assume that $\omega_{\mathcal{K}_\mu}=
\tau_{\mathcal{K}_\mu}\cdot \omega_{\mathcal{M}_\mu}$ is the
restriction of the symplectic form $\omega_{\mathcal{M}_\mu}$ on
$T^*\mathcal{M}_\mu$ fibrewise to the distribution $\mathcal{K}_\mu$,
where $\tau_{\mathcal{K}_\mu}$ is the restriction map to distribution
$\mathcal{K}_\mu$. \\

From the above construction, we know that
$\omega_{\mathcal{K}_\mu}$ is non-degenerate, and is called
as a $\mathbf{J}$-nonholonomic $R_p$-reduced distributional two-form to
avoid any confusion. Because $\omega_{\mathcal{K}_\mu}$ is
non-degenerate as a bilinear form on each fibre of
$\mathcal{K}_\mu$, there exists a vector field $X_{\mathcal{K}_\mu}$
on $\mathcal{M}_\mu$, which takes values in the constraint
distribution $\mathcal{K}_\mu$, such that the $\mathbf{J}$-nonholonomic
$R_p$-reduced distributional Hamiltonian equation holds, that is,
$\mathbf{i}_{X_{\mathcal{K}_\mu}}\omega_{\mathcal{K}_\mu}=\mathbf{d}h_{\mathcal{K}_\mu}$,
if the admissibility condition $\mathrm{dim}\mathcal{M}_\mu=
\mathrm{rank}\mathcal{F}_\mu$ and the compatibility condition
$T\mathcal{M}_\mu\cap \mathcal{F}_\mu^\bot= \{0\}$ hold, where
$\mathcal{F}_\mu^\bot$ denotes the symplectic orthogonal of
$\mathcal{F}_\mu$ with respect to the $\mathbf{J}$-nonholonomic $R_p$-reduced symplectic form
$\omega_\mu$, and $\mathbf{d}h_{\mathcal{K}_\mu}$ is the restriction
of $\mathbf{d}h_{\mathcal{M}_\mu}$ to $\mathcal{K}_\mu$,
and the function $h_{\mathcal {K}_\mu}$ satisfies
$\mathbf{d}h_{\mathcal{K}_\mu}= \tau_{\mathcal{K}_\mu}\cdot \mathbf{d}h_{\mathcal{M}_\mu} $,
and $h_{\mathcal{M}_\mu}= \tau_{\mathcal{M}_\mu}\cdot h_\mu$ is the
restriction of $h_\mu$ to $\mathcal{M}_\mu$,
and $h_\mu$ is the Marsden-Weinstein point reduced
Hamiltonian function $h_\mu: (T^* Q)_\mu \rightarrow \mathbb{R}$ defined
by $h_\mu\cdot \pi_\mu= H\cdot i_\mu$. Thus, the geometrical formulation
of the $\mathbf{J}$-nonholonomic $R_p$-reduced distributional Hamiltonian
system may be summarized as follows.

\begin{defi} ($\mathbf{J}$-Nonholonomic $R_p$-reduced Distributional Hamiltonian System)
Assume that the 6-tuple $(T^*Q,G,\omega,\mathbf{J},\mathcal{D},H)$
is a nonholonomic Hamiltonian system with symmetry and momentum map,
where $\omega$ is the canonical symplectic form on $T^* Q$,
and $\mathcal{D}\subset TQ$ is a $\mathcal{D}$-completely and
$\mathcal{D}$-regularly nonholonomic constraint of the system, and
$\mathcal{D}$ and $H$ are both $G$-invariant. For a regular value
$\mu\in\mathfrak{g}^\ast$ of the momentum map $\mathbf{J}:T^\ast
Q\rightarrow \mathfrak{g}^\ast$, assume that there exists a
$\mathbf{J}$-nonholonomic $R_p$-reduced distribution
$\mathcal{K}_\mu$, an associated non-degenerate and $\mathbf{J}$-nonholonomic
$R_p$-reduced distributional two-form
$\omega_{\mathcal{K}_\mu}$ and a vector field $X_{\mathcal {K}_\mu}$
on the $\mathbf{J}$-nonholonomic $R_p$-reduced constraint submanifold
$\mathcal{M}_\mu=(\mathcal{M}\cap \mathbf{J}^{-1}(\mu)) /G_\mu, $
where $\mathcal{M}=\mathcal{F}L(\mathcal{D}),$ and $\mathcal{M}\cap
\mathbf{J}^{-1}(\mu)\neq \emptyset, $ and $G_\mu=\{g\in G \;| \;
\operatorname{Ad}_g^\ast \mu=\mu \}$, such that the
$\mathbf{J}$-nonholonomic $R_p$-reduced distributional Hamiltonian
equation $\mathbf{i}_{X_{\mathcal{K}_\mu}}\omega_{\mathcal{K}_\mu} =
\mathbf{d}h_{\mathcal{K}_\mu}$ holds, where
$\mathbf{d}h_{\mathcal{K}_\mu}$ is the restriction of
$\mathbf{d}h_{\mathcal{M}_\mu}$ to $\mathcal{K}_\mu$, and
and the function $h_{\mathcal {K}_\mu}$ is defined above. Then the triple
$(\mathcal{K}_\mu,\omega_{\mathcal {K}_\mu},h_{\mathcal {K}_\mu})$ is called a
$\mathbf{J}$-nonholonomic $R_p$-reduced distributional Hamiltonian system
of the system $(T^*Q,G,\omega,\mathbf{J},\mathcal{D},H)$, and $X_{\mathcal
{K}_\mu}$ is the dynamical vector field of the
$\mathbf{J}$-nonholonomic $R_p$-reduced distributional Hamiltonian system
$(\mathcal{K}_\mu,\omega_{\mathcal{K}_\mu},h_{\mathcal {K}_\mu})$.
Under the above circumstances, we refer to
$(T^*Q,G,\omega,\mathbf{J},\mathcal{D},H)$ as
a $\mathbf{J}$-nonholonomic point reducible Hamiltonian system
with an associated $\mathbf{J}$-nonholonomic $R_p$-reduced
distributional Hamiltonian system
$(\mathcal{K}_\mu,\omega_{\mathcal{K}_\mu},h_{\mathcal {K}_\mu})$.
\end{defi}

Since the non-degenerate and $\mathbf{J}$-nonholonomic $R_p$-reduced
distributional two-form $\omega_{\mathcal{K}_\mu}$ may not be symplectic,
and the $\mathbf{J}$-nonholonomic $R_p$-reduced distributional Hamiltonian system
$(\mathcal{K}_\mu,\omega_{\mathcal {K}_\mu},h_{\mathcal {K}_\mu})$
may not be yet a Hamiltonian system, and has no yet generating function,
and hence we can not describe the Hamilton-Jacobi equation for a
$\mathbf{J}$-nonholonomic $R_p$-reduced
distributional Hamiltonian system just like as in Theorem 1.1.
But, for a given $\mathbf{J}$-nonholonomic regular point reducible Hamiltonian system $(T^*Q,G,\omega,\mathbf{J},\mathcal{D},H)$ with an associated
$\mathbf{J}$-nonholonomic $R_p$-reduced distributional Hamiltonian system
$(\mathcal{K}_\mu,\omega_{\mathcal {K}_\mu},h_{\mathcal {K}_\mu})$, by using Lemma 3.4,
we can derive precisely
the geometric constraint conditions of the $\mathbf{J}$-nonholonomic $R_p$-reduced distributional two-form
$\omega_{\mathcal{K}_\mu}$ for the $\mathbf{J}$-nonholonomic regular point reducible dynamical vector field,
that is, the two types of Hamilton-Jacobi equation for the
$\mathbf{J}$-nonholonomic $R_p$-reduced distributional Hamiltonian system
$(\mathcal{K}_\mu,\omega_{\mathcal {K}_\mu},h_{\mathcal {K}_\mu})$.
At first, by using the fact that the one-form $\gamma: Q
\rightarrow T^*Q $ is closed on $\mathcal{D}$ with respect to
$T\pi_Q: TT^* Q \rightarrow TQ, $
and $\textmd{Im}(\gamma)\subset \mathcal{M} \cap
\mathbf{J}^{-1}(\mu), $ and it is $G_\mu$-invariant,
as well as $ \textmd{Im}(T\bar{\gamma}_\mu)\subset \mathcal{K}_\mu, $
we can prove the Type I of
Hamilton-Jacobi theorem for the $\mathbf{J}$-nonholonomic $R_p$-reduced distributional
Hamiltonian system. For convenience, the maps involved in the
following theorem and its proof are shown in Diagram-5.
\begin{center}
\hskip 0cm \xymatrix{ \mathbf{J}^{-1}(\mu) \ar[r]^{i_\mu}
& T^* Q  \ar[r]^{\pi_Q}
& Q \ar[d]_{X_H^\gamma} \ar[r]^{\gamma}
& T^*Q \ar[d]_{X_H} \ar[r]^{\pi_\mu} & (T^* Q)_\mu \ar[d]_{X_{h_\mu}}
& \mathcal{M}_\mu  \ar[l]_{i_{\mathcal{M}_\mu}} \ar[d]_{X_{\mathcal{K}_\mu}}\\
& T(T^*Q)  & TQ \ar[l]_{T\gamma} & T(T^*Q) \ar[l]^{T\pi_Q} \ar[r]_{T\pi_\mu}
& T(T^* Q)_\mu \ar[r]^{\tau_{\mathcal{K}_\mu}} & \mathcal{K}_\mu }
\end{center}
$$\mbox{Diagram-5}$$

\begin{theo} (Type I of Hamilton-Jacobi Theorem for a $\mathbf{J}$-Nonholonomic
$R_p$-reduced Distributional Hamiltonian System) For a given
$\mathbf{J}$-nonholonomic point reducible Hamiltonian system
$(T^*Q,G,\omega,\mathbf{J},\mathcal{D},H)$ with an associated
$\mathbf{J}$-nonholonomic $R_p$-reduced distributional Hamiltonian system
$(\mathcal{K}_\mu,\omega_{\mathcal{K}_\mu},h_{\mathcal {K}_\mu})$, assume that $\gamma:
Q \rightarrow T^*Q$ is an one-form on $Q$, and $X_H^\gamma =
T\pi_{Q}\cdot X_H \cdot \gamma, $ where $X_{H}$ is the dynamical
vector field of the corresponding unconstrained Hamiltonian system with
symmetry and momentum map $(T^*Q,G,\omega,\mathbf{J},H)$. Moreover,
assume that $\mu\in\mathfrak{g}^\ast$ is a regular value of the momentum
map $\mathbf{J}$, and $\textmd{Im}(\gamma)\subset \mathcal{M} \cap
\mathbf{J}^{-1}(\mu), $ and it is $G_\mu$-invariant, and
$\bar{\gamma}_\mu=\pi_\mu(\gamma): Q \rightarrow \mathcal{M}_\mu $,
and $ \textmd{Im}(T\bar{\gamma}_\mu)\subset \mathcal{K}_\mu. $
If the one-form $\gamma: Q \rightarrow T^*Q $ is closed on $\mathcal{D}$ with respect to
$T\pi_Q: TT^* Q \rightarrow TQ, $ then
$\bar{\gamma}_\mu$ is a solution of the
equation $T\bar{\gamma}_\mu\cdot
X_H^\gamma= X_{\mathcal{K}_\mu}\cdot \bar{\gamma}_\mu. $
Here $X_{\mathcal{K}_\mu}$ is the dynamical vector field of the reduced
system $(\mathcal{K}_\mu,\omega_{\mathcal{K}_\mu},h_{\mathcal {K}_\mu})$. The equation
$T\bar{\gamma}_\mu\cdot X_H^\gamma= X_{\mathcal{K}_\mu}\cdot
\bar{\gamma}_\mu,$ is called the Type I of Hamilton-Jacobi equation for the
$\mathbf{J}$-nonholonomic $R_p$-reduced distributional Hamiltonian system
$(\mathcal{K}_\mu,\omega_{\mathcal{K}_\mu},h_{\mathcal {K}_\mu})$.
\end{theo}

\noindent{\bf Proof: } At first, from Theorem 3.5, we know that
$\gamma$ is a solution of the Hamilton-Jacobi equation
$T\gamma\cdot X_H^\gamma= X_{\mathcal{K}}\cdot \gamma .$ Next, we note that
$\textmd{Im}(\gamma)\subset \mathcal{M} \cap \mathbf{J}^{-1}(\mu), $
and it is $G_\mu$-invariant, in this case,
$\pi_\mu^*\omega_\mu=
i_\mu^*\omega= \omega, $ along $\textmd{Im}(\gamma)$.
On the other hand, because
$\textmd{Im}(T\bar{\gamma}_\mu)\subset \mathcal{K}_\mu, $
then $\omega_{\mathcal{K}_\mu}\cdot
\tau_{\mathcal{K}_\mu}=\tau_{\mathcal{K}_\mu}\cdot
\omega_{\mathcal{M}_\mu}= \tau_{\mathcal{K}_\mu}\cdot
i_{\mathcal{M}_\mu}^* \cdot \omega_\mu, $ along
$\textmd{Im}(T\bar{\gamma}_\mu)$. Thus, using
the $\mathbf{J}$-nonholonomic $R_p$-reduced distributional two-form
$\omega_{\mathcal{K}_\mu}$, from Lemma 3.4, if we take that $v=
\tau_{\mathcal{K}_\mu}\cdot T\pi_\mu \cdot X_{H}\cdot \gamma= X_{\mathcal{K}_\mu}\cdot
\bar{\gamma}_\mu \in \mathcal{K}_\mu, $ and for any $w \in
\mathcal{F}, \; T\lambda(w)\neq 0, $ and
$\tau_{\mathcal{K}_\mu}\cdot T\pi_\mu \cdot w \neq 0, $ then we have that
\begin{align*}
& \omega_{\mathcal{K}_\mu}(T\bar{\gamma}_\mu \cdot X_H^\gamma, \;
\tau_{\mathcal{K}_\mu}\cdot T\pi_\mu \cdot w)=
\omega_{\mathcal{K}_\mu}(\tau_{\mathcal{K}_\mu}\cdot
T\bar{\gamma}_\mu \cdot X_H^\gamma, \; \tau_{\mathcal{K}_\mu}\cdot T\pi_\mu \cdot w)\\
& = \tau_{\mathcal{K}_\mu}\cdot \omega_{\mathcal{M}_\mu}(T(\pi_\mu
\cdot\gamma )\cdot X_H^\gamma, \; T\pi_\mu \cdot w ) =
\tau_{\mathcal{K}_\mu}\cdot i_{\mathcal{M}_\mu}^* \cdot
\omega_\mu(T\pi_\mu \cdot T\gamma \cdot X_H^\gamma, \; T\pi_\mu
\cdot w )\\
& = \tau_{\mathcal{K}_\mu}\cdot i_{\mathcal{M}_\mu}^*
\cdot \pi_\mu^*\omega_\mu(T\gamma \cdot T\pi_Q \cdot X_H \cdot
\gamma, \; w) = \tau_{\mathcal{K}_\mu}\cdot i_{\mathcal{M}_\mu}^*
\cdot \omega(T(\gamma \cdot \pi_Q) \cdot X_H \cdot \gamma, \; w)\\
& =\tau_{\mathcal{K}_\mu}\cdot i_{\mathcal{M}_\mu}^* \cdot
(\omega(X_H \cdot \gamma, \; w-T(\gamma
\cdot \pi_Q)\cdot w) -\mathbf{d}\gamma(T\pi_{Q}(X_H\cdot \gamma), \; T\pi_{Q}(w)))\\
& = \tau_{\mathcal{K}_\mu}\cdot i_{\mathcal{M}_\mu}^*
\cdot\pi_\mu^*\omega_\mu(X_H \cdot \gamma, \; w) -
\tau_{\mathcal{K}_\mu}\cdot i_{\mathcal{M}_\mu}^*
\cdot\pi_\mu^*\omega_\mu(X_H \cdot
\gamma, \; T(\gamma
\cdot \pi_Q) \cdot w)\\
& \;\;\;\;\;\; -\tau_{\mathcal{K}_\mu}\cdot i_{\mathcal{M}_\mu}^* \cdot\mathbf{d}\gamma(T\pi_{Q}(X_H\cdot \gamma), \; T\pi_{Q}(w))\\
& = \tau_{\mathcal{K}_\mu}\cdot i_{\mathcal{M}_\mu}^*
\cdot\omega_\mu(T\pi_\mu(X_H\cdot \gamma), \; T\pi_\mu \cdot w) -
\tau_{\mathcal{K}_\mu}\cdot i_{\mathcal{M}_\mu}^* \cdot
\omega_\mu(T\pi_\mu\cdot(X_H\cdot \gamma), \; T(\pi_\mu
\cdot\gamma) \cdot T\pi_{Q}(w))\\
& \;\;\;\;\;\; -\tau_{\mathcal{K}_\mu}\cdot
i_{\mathcal{M}_\mu}^* \cdot\mathbf{d}\gamma(T\pi_{Q}(X_H\cdot
\gamma), \; T\pi_{Q}(w))\\
& = \tau_{\mathcal{K}_\mu}\cdot i_{\mathcal{M}_\mu}^* \cdot
\omega_\mu(X_{h_{\mathcal {K}_\mu}} \cdot \bar{\gamma}_\mu, \; T\pi_\mu \cdot w)-
\tau_{\mathcal{K}_\mu}\cdot i_{\mathcal{M}_\mu}^* \cdot
\omega_\mu(X_{h_{\mathcal {K}_\mu}} \cdot \bar{\gamma}_\mu, \; T\bar{\gamma}_\mu \cdot T\pi_{Q}(w))\\
& \;\;\;\;\;\; -\tau_{\mathcal{K}_\mu}\cdot i_{\mathcal{M}_\mu}^* \cdot\mathbf{d}\gamma(T\pi_{Q}(X_H\cdot \gamma), \; T\pi_{Q}(w))\\
& = \omega_{\mathcal{K}_\mu}( \tau_{\mathcal{K}_\mu}\cdot X_{h_{\mathcal {K}_\mu}}
\cdot \bar{\gamma}_\mu, \; \tau_{\mathcal{K}_\mu}\cdot T\pi_\mu
\cdot w) - \omega_{\mathcal{K}_\mu}(\tau_{\mathcal{K}_\mu}\cdot
X_{h_{\mathcal {K}_\mu}} \cdot \bar{\gamma}_\mu, \;
\tau_{\mathcal{K}_\mu}\cdot T\bar{\gamma}_\mu \cdot T\pi_{Q}(w))\\
& \;\;\;\;\;\; -\tau_{\mathcal{K}_\mu}\cdot i_{\mathcal{M}_\mu}^* \cdot\mathbf{d}\gamma(T\pi_{Q}(X_H\cdot \gamma), \; T\pi_{Q}(w))\\
& = \omega_{\mathcal{K}_\mu}(X_{\mathcal{K}_\mu}\cdot
\bar{\gamma}_\mu, \; \tau_{\mathcal{K}_\mu} \cdot T\pi_\mu \cdot w)
- \omega_{\mathcal{K}_\mu}(X_{\mathcal{K}_\mu} \cdot
\bar{\gamma}_\mu, \; T\bar{\gamma}_\mu
\cdot T\pi_{Q}(w))\\
& \;\;\;\;\;\; -\tau_{\mathcal{K}_\mu}\cdot i_{\mathcal{M}_\mu}^*
\cdot\mathbf{d}\gamma(T\pi_{Q}(X_H\cdot \gamma), \; T\pi_{Q}(w)),
\end{align*}
where we have used that $ \tau_{\mathcal{K}_\mu}\cdot T\bar{\gamma}_\mu=
T\bar{\gamma}_\mu, $ and $\tau_{\mathcal{K}_\mu}\cdot X_{h_{\mathcal {K}_\mu}}\cdot
\bar{\gamma}_\mu = X_{\mathcal{K}_\mu}\cdot \bar{\gamma}_\mu, $
since $\textmd{Im}(T\bar{\gamma}_\mu)\subset \mathcal{K}_\mu. $
If the one-form $\gamma: Q \rightarrow T^*Q $ is closed on $\mathcal{D}$ with respect to
$T\pi_Q: TT^* Q \rightarrow TQ, $ then we have that
$\mathbf{d}\gamma(T\pi_{Q}(X_H\cdot \gamma), \; T\pi_{Q}(w))=0, $
since $X_{H}\cdot \gamma, \; w \in \mathcal{F},$ and
$T\pi_{Q}(X_H\cdot \gamma), \; T\pi_{Q}(w) \in \mathcal{D}, $ and hence
$$
\tau_{\mathcal{K}_\mu}\cdot
i_{\mathcal{M}_\mu}^* \cdot\mathbf{d}\gamma(T\pi_{Q}(X_H\cdot \gamma),
\; T\pi_{Q}(w))=0,
$$
and
\begin{align}
& \omega_{\mathcal{K}_\mu}(T\bar{\gamma}_\mu \cdot X_H^\gamma, \;
\tau_{\mathcal{K}_\mu}\cdot T\pi_\mu \cdot w)- \omega_{\mathcal{K}_\mu}(X_{\mathcal{K}_\mu}\cdot
\bar{\gamma}_\mu, \; \tau_{\mathcal{K}_\mu} \cdot T\pi_\mu \cdot w) \nonumber \\
& = -\omega_{\mathcal{K}_\mu}(X_{\mathcal{K}_\mu} \cdot
\bar{\gamma}_\mu, \; T\bar{\gamma}_\mu
\cdot T\pi_{Q}(w)).
\label {5.2}\end{align}
If $\bar{\gamma}_\mu$ satisfies the equation
$T\bar{\gamma}_\mu\cdot X_H^\gamma= X_{\mathcal{K}_\mu}\cdot
\bar{\gamma}_\mu ,$
from Lemma 3.4(i) we know that the right side of (5.2) becomes
\begin{align*}
& -\omega_{\mathcal{K}_\mu}(X_{\mathcal{K}_\mu} \cdot
\bar{\gamma}_\mu, \; \tau_{\mathcal{K}_\mu}\cdot T\bar{\gamma}_\mu
\cdot T\pi_{Q}(w))\\
& = -\omega_{\mathcal{K}_\mu}(T\bar{\gamma}_\mu\cdot X_H^\gamma,
\; T\bar{\gamma}_\mu \cdot T\pi_{Q}(w))\\
& = -\omega_{\mathcal{K}_\mu}(\tau_{\mathcal{K}_\mu}T\bar{\gamma}_\mu\cdot X_H^\gamma,
\; \tau_{\mathcal{K}_\mu}\cdot T\bar{\gamma}_\mu \cdot T\pi_{Q}(w))\\
& = -\tau_{\mathcal{K}_\mu}\cdot i_{\mathcal{M}_\mu}^* \cdot
\omega_\mu(T\bar{\gamma}_\mu\cdot X_H^\gamma,
\; T\bar{\gamma}_\mu \cdot T\pi_{Q}(w))\\
& = -\tau_{\mathcal{K}_\mu}\cdot i_{\mathcal{M}_\mu}^* \cdot \bar{\gamma}_\mu^* \cdot
\omega_\mu( T\pi_{Q} \cdot X_{H} \cdot
\gamma, \; T\pi_{Q}(w))\\
& = -\tau_{\mathcal{K}_\mu}\cdot i_{\mathcal{M}_\mu}^* \cdot
 \gamma^* \cdot \pi^*_\mu\cdot \omega_\mu(T\pi_{Q} \cdot X_{H} \cdot
\gamma, \; T\pi_{Q}(w))\\
& = -\tau_{\mathcal{K}_\mu}\cdot
i_{\mathcal{M}_\mu}^* \cdot\gamma^*\omega( T\pi_{Q}(X_{H}\cdot\gamma), \; T\pi_{Q}(w))\\
& = \tau_{\mathcal{K}_\mu}\cdot i_{\mathcal{M}_\mu}^* \cdot
\mathbf{d}\gamma(T\pi_{Q}( X_{H}\cdot\gamma ), \; T\pi_{Q}(w))=0.
\end{align*}
But, because the $\mathbf{J}$-nonholonomic $R_p$-reduced distributional two-form
$\omega_{\mathcal{K}_\mu}$ is non-degenerate,
the left side of (5.2) equals zero, only when
$\bar{\gamma}_\mu$ satisfies the equation 
$T\bar{\gamma}_\mu\cdot X_H^\gamma= X_{\mathcal{K}_\mu}\cdot
\bar{\gamma}_\mu .$ Thus,
if the one-form $\gamma: Q \rightarrow T^*Q $ is closed on $\mathcal{D}$ with respect to
$T\pi_Q: TT^* Q \rightarrow TQ, $ then $\bar{\gamma}_\mu$ must be a solution
of the Type I of Hamilton-Jacobi equation
$T\bar{\gamma}_\mu\cdot X_H^\gamma= X_{\mathcal{K}_\mu}\cdot
\bar{\gamma}_\mu. $
\hskip 0.3cm $\blacksquare$\\

Next, for any $G_\mu$-invariant symplectic map $\varepsilon: T^* Q \rightarrow T^* Q $,
we can prove the following Type II of
Hamilton-Jacobi theorem for the $\mathbf{J}$-nonholonomic $R_p$-reduced distributional
Hamiltonian system.
For convenience, the maps involved in the following
theorem and its proof are shown in Diagram-6.
\begin{center}
\hskip 0cm \xymatrix{ \mathbf{J}^{-1}(\mu) \ar[r]^{i_\mu} & T^* Q
\ar[d]_{X_{H\cdot \varepsilon}} \ar[dr]^{X_H^\varepsilon} \ar[r]^{\pi_Q}
& Q \ar[r]^{\gamma} & T^*Q \ar[d]_{X_H} \ar[dr]^{X_{h_\mu} \cdot\bar{\varepsilon}} \ar[r]^{\pi_\mu}
& (T^* Q)_\mu \ar[d]^{X_{h_\mu}}
& \mathcal{M}_\mu  \ar[l]_{i_{\mathcal{M}_\mu}} \ar[d]_{X_{\mathcal{K}_\mu}}\\
& T(T^*Q)  & TQ \ar[l]^{T\gamma} & T(T^*Q) \ar[l]^{T\pi_Q} \ar[r]_{T\pi_\mu}
& T(T^* Q)_\mu \ar[r]^{\tau_{\mathcal{K}_\mu}} & \mathcal{K}_\mu }
\end{center}
$$\mbox{Diagram-6}$$

\begin{theo} (Type II of Hamilton-Jacobi Theorem for a $\mathbf{J}$-Nonholonomic
$R_p$-reduced Distributional Hamiltonian System) For a given
$\mathbf{J}$-nonholonomic point reducible Hamiltonian system
$(T^*Q,G,\omega,\mathbf{J},\mathcal{D},H)$ with an associated
$\mathbf{J}$-nonholonomic $R_p$-reduced distributional Hamiltonian system
$(\mathcal{K}_\mu,\omega_{\mathcal{K}_\mu},h_{\mathcal {K}_\mu})$, assume that $\gamma:
Q \rightarrow T^*Q$ is an one-form on $Q$, and $\lambda=\gamma \cdot
\pi_{Q}: T^* Q \rightarrow T^* Q, $ and for any
symplectic map $\varepsilon:T^* Q \rightarrow T^* Q, $
denote $X_H^\varepsilon =
T\pi_{Q}\cdot X_H \cdot \varepsilon$, where $X_{H}$ is the dynamical
vector field of the corresponding unconstrained Hamiltonian system with
symmetry and momentum map $(T^*Q,G,\omega,\mathbf{J},H)$. Moreover,
assume that $\mu\in\mathfrak{g}^\ast$ is a regular value of the momentum
map $\mathbf{J}$, and $\textmd{Im}(\gamma)\subset \mathcal{M} \cap
\mathbf{J}^{-1}(\mu), $ and that it is $G_\mu$-invariant,
and $\varepsilon$ is $G_\mu$-invariant and
$\varepsilon(\mathbf{J}^{-1}(\mu)) \subset \mathbf{J}^{-1}(\mu). $ Denote
$\bar{\gamma}_\mu=\pi_\mu(\gamma): Q \rightarrow \mathcal{M}_\mu $,
and $ \textmd{Im}(T\bar{\gamma}_\mu)\subset \mathcal{K}_\mu, $ and
$\bar{\lambda}_\mu=\pi_\mu(\lambda): \mathbf{J}^{-1}(\mu) (\subset T^* Q) \rightarrow
\mathcal{M}_\mu, $ and
$\bar{\varepsilon}_\mu=\pi_\mu(\varepsilon): \mathbf{J}^{-1}(\mu) (\subset T^* Q) \rightarrow
\mathcal{M}_\mu. $ Then $\varepsilon$ and $\bar{\varepsilon}_\mu$ satisfy the
equation $\tau_{\mathcal{K}_\mu} \cdot T\bar{\varepsilon}(X_{h_{\mathcal {K}_\mu}\cdot \bar{\varepsilon}_\mu})
= T\bar{\lambda}_\mu \cdot X_H \cdot\varepsilon, $ if and only if
they satisfy the equation $T\bar{\gamma}_\mu\cdot
X_H^\varepsilon= X_{\mathcal{K}_\mu}\cdot \bar{\varepsilon}_\mu. $
Here $X_{h_{\mathcal {K}_\mu} \cdot\bar{\varepsilon}_\mu}$ is the Hamiltonian vector field of the
function $h_{\mathcal {K}_\mu}\cdot \bar{\varepsilon}_\mu: T^* Q\rightarrow \mathbb{R}, $
and $X_{\mathcal{K}_\mu}$ is the dynamical vector field of the reduced
system $(\mathcal{K}_\mu,\omega_{\mathcal{K}_\mu},h_{\mathcal {K}_\mu})$. The equation
$T\bar{\gamma}_\mu\cdot X_H^\varepsilon= X_{\mathcal{K}_\mu}\cdot
\bar{\varepsilon}_\mu,$ is called the Type II of Hamilton-Jacobi equation for the
$\mathbf{J}$-nonholonomic $R_p$-reduced distributional Hamiltonian system
$(\mathcal{K}_\mu,\omega_{\mathcal{K}_\mu},h_{\mathcal {K}_\mu})$.
\end{theo}

\noindent{\bf Proof: } At first, we note that
$\textmd{Im}(\gamma)\subset \mathcal{M} \cap \mathbf{J}^{-1}(\mu), $
and it is $G_\mu$-invariant, in this case,
$\pi_\mu^*\omega_\mu=
i_\mu^*\omega= \omega, $ along $\textmd{Im}(\gamma)$.
On the other hand, because
$\textmd{Im}(T\bar{\gamma}_\mu)\subset \mathcal{K}_\mu, $
then $\omega_{\mathcal{K}_\mu}\cdot
\tau_{\mathcal{K}_\mu}=\tau_{\mathcal{K}_\mu}\cdot
\omega_{\mathcal{M}_\mu}= \tau_{\mathcal{K}_\mu}\cdot
i_{\mathcal{M}_\mu}^* \cdot \omega_\mu, $ along
$\textmd{Im}(T\bar{\gamma}_\mu)$. Thus, using
the $\mathbf{J}$-nonholonomic $R_p$-reduced distributional two-form
$\omega_{\mathcal{K}_\mu}$, from Lemma 3.4, if we take that $v=
\tau_{\mathcal{K}_\mu}\cdot T\pi_\mu \cdot X_{H}\cdot \varepsilon= X_{\mathcal{K}_\mu}\cdot
\bar{\varepsilon}_\mu \in \mathcal{K}_\mu, $ and for any $w \in
\mathcal{F}, \; T\lambda(w)\neq 0, $ and
$\tau_{\mathcal{K}_\mu}\cdot T\pi_\mu \cdot w \neq 0, $ then we have that
\begin{align*}
& \omega_{\mathcal{K}_\mu}(T\bar{\gamma}_\mu \cdot X_H^\varepsilon, \;
\tau_{\mathcal{K}_\mu}\cdot T\pi_\mu \cdot w)=
\omega_{\mathcal{K}_\mu}(\tau_{\mathcal{K}_\mu}\cdot
T\bar{\gamma}_\mu \cdot X_H^\varepsilon, \; \tau_{\mathcal{K}_\mu}\cdot T\pi_\mu \cdot w)\\
& = \tau_{\mathcal{K}_\mu}\cdot \omega_{\mathcal{M}_\mu}(T(\pi_\mu
\cdot\gamma )\cdot X_H^\varepsilon, \; T\pi_\mu \cdot w ) =
\tau_{\mathcal{K}_\mu}\cdot i_{\mathcal{M}_\mu}^* \cdot
\omega_\mu(T\pi_\mu \cdot T\gamma \cdot X_H^\varepsilon, \; T\pi_\mu
\cdot w )\\
& = \tau_{\mathcal{K}_\mu}\cdot i_{\mathcal{M}_\mu}^*
\cdot \pi_\mu^*\omega_\mu(T\gamma \cdot T\pi_Q \cdot X_H \cdot
\varepsilon, \; w) = \tau_{\mathcal{K}_\mu}\cdot i_{\mathcal{M}_\mu}^*
\cdot \omega(T(\gamma \cdot \pi_Q) \cdot X_H \cdot \varepsilon, \; w)\\
& =\tau_{\mathcal{K}_\mu}\cdot i_{\mathcal{M}_\mu}^* \cdot
(\omega(X_H \cdot \varepsilon, \; w-T(\gamma
\cdot \pi_Q)\cdot w) -\mathbf{d}\gamma(T\pi_{Q}(X_H\cdot \varepsilon), \; T\pi_{Q}(w)))\\
& = \tau_{\mathcal{K}_\mu}\cdot i_{\mathcal{M}_\mu}^* \cdot
\omega(X_H \cdot \varepsilon, \; w) - \tau_{\mathcal{K}_\mu}\cdot
i_{\mathcal{M}_\mu}^* \cdot \omega(X_H \cdot
\varepsilon, \; T\lambda \cdot w)\\
& \;\;\;\;\;\; -\tau_{\mathcal{K}_\mu}\cdot i_{\mathcal{M}_\mu}^*
\cdot\mathbf{d}\gamma(T\pi_{Q}(X_H\cdot \varepsilon), \; T\pi_{Q}(w))\\
& = \tau_{\mathcal{K}_\mu}\cdot i_{\mathcal{M}_\mu}^*
\cdot\pi_\mu^*\omega_\mu(X_H \cdot \varepsilon, \; w) -
\tau_{\mathcal{K}_\mu}\cdot i_{\mathcal{M}_\mu}^*
\cdot\pi_\mu^*\omega_\mu(X_H \cdot
\varepsilon, \; T\lambda \cdot w)\\
& \;\;\;\;\;\; +\tau_{\mathcal{K}_\mu}\cdot i_{\mathcal{M}_\mu}^* \cdot \lambda^* \omega(X_H\cdot \varepsilon, \; w)\\
& = \tau_{\mathcal{K}_\mu}\cdot i_{\mathcal{M}_\mu}^*
\cdot\omega_\mu(T\pi_\mu(X_H\cdot \varepsilon), \; T\pi_\mu \cdot w) -
\tau_{\mathcal{K}_\mu}\cdot i_{\mathcal{M}_\mu}^* \cdot
\omega_\mu(T\pi_\mu\cdot(X_H\cdot \varepsilon), \; T(\pi_\mu
\cdot\lambda) \cdot w)\\
& \;\;\;\;\;\; +\tau_{\mathcal{K}_\mu}\cdot
i_{\mathcal{M}_\mu}^* \cdot
\pi_\mu^*\omega_\mu(T\lambda\cdot X_H \cdot
\varepsilon, \; T\lambda \cdot w)\\
& = \tau_{\mathcal{K}_\mu}\cdot
i_{\mathcal{M}_\mu}^* \cdot\omega_\mu(T\pi_\mu(X_H)\cdot
\pi_\mu(\varepsilon), \; T\pi_\mu \cdot w) - \tau_{\mathcal{K}_\mu}\cdot
i_{\mathcal{M}_\mu}^* \cdot \omega_\mu(T\pi_\mu(X_H)\cdot
\pi_\mu(\varepsilon), \; T\bar{\lambda}_\mu \cdot w)\\
& \;\;\;\;\;\; +\tau_{\mathcal{K}_\mu}\cdot
i_{\mathcal{M}_\mu}^* \cdot
\omega_\mu(T\pi_\mu\cdot T\lambda\cdot X_H \cdot
\varepsilon, \; T\pi_\mu\cdot T\lambda \cdot w)\\
& = \tau_{\mathcal{K}_\mu}\cdot i_{\mathcal{M}_\mu}^* \cdot
\omega_\mu(X_{h_{\mathcal {K}_\mu}} \cdot \bar{\varepsilon}_\mu, \; T\pi_\mu \cdot w)-
\tau_{\mathcal{K}_\mu}\cdot i_{\mathcal{M}_\mu}^* \cdot
\omega_\mu(X_{h_{\mathcal {K}_\mu}} \cdot \bar{\varepsilon}_\mu, \; T\bar{\lambda}_\mu \cdot w)\\
& \;\;\;\;\;\; +\tau_{\mathcal{K}_\mu}\cdot
i_{\mathcal{M}_\mu}^* \cdot
\omega_\mu(T\bar{\lambda}_\mu \cdot X_H \cdot
\varepsilon, \; T\bar{\lambda}_\mu \cdot w)\\
& = \omega_{\mathcal{K}_\mu}( \tau_{\mathcal{K}_\mu}\cdot X_{h_{\mathcal {K}_\mu}}
\cdot \bar{\varepsilon}_\mu, \; \tau_{\mathcal{K}_\mu}\cdot T\pi_\mu
\cdot w) - \omega_{\mathcal{K}_\mu}(\tau_{\mathcal{K}_\mu}\cdot
X_{h_{\mathcal {K}_\mu}} \cdot \bar{\varepsilon}_\mu, \;
\tau_{\mathcal{K}_\mu}\cdot T\bar{\lambda}_\mu \cdot w)\\
& \;\;\;\;\;\; +\omega_{\mathcal{K}_\mu}(\tau_{\mathcal{K}_\mu}\cdot T\bar{\lambda}_\mu \cdot X_H \cdot
\varepsilon, \; \tau_{\mathcal{K}_\mu}\cdot T\bar{\lambda}_\mu \cdot w)\\
& = \omega_{\mathcal{K}_\mu}(X_{\mathcal{K}_\mu}\cdot
\bar{\varepsilon}_\mu, \; \tau_{\mathcal{K}_\mu} \cdot T\pi_\mu \cdot w)
- \omega_{\mathcal{K}_\mu}(\tau_{\mathcal{K}_\mu}\cdot X_{h_{\mathcal {K}_\mu}} \cdot
\bar{\varepsilon}_\mu, \; T\bar{\lambda}_\mu \cdot w)\\
& \;\;\;\;\;\; +\omega_{\mathcal{K}_\mu}(T\bar{\lambda}_\mu \cdot X_H \cdot
\varepsilon, \; T\bar{\lambda}_\mu \cdot w),
\end{align*}
where we have used that $ \tau_{\mathcal{K}_\mu}\cdot T\bar{\gamma}_\mu=
T\bar{\gamma}_\mu, $ and $\tau_{\mathcal{K}_\mu}\cdot X_{h_{\mathcal {K}_\mu}}\cdot
\bar{\varepsilon}_\mu = X_{\mathcal{K}_\mu}\cdot \bar{\varepsilon}_\mu, $
since $\textmd{Im}(T\bar{\gamma}_\mu)\subset \mathcal{K}_\mu. $ Note
that $\varepsilon: T^* Q \rightarrow T^* Q $ is symplectic, and
$\pi_\mu^*\omega_\mu= i_\mu^*\omega= \omega, $ along
$\textmd{Im}(\gamma)$, and hence $\bar{\varepsilon}_\mu=
\pi_\mu(\varepsilon): T^* Q \rightarrow (T^* Q)_\mu $ is also symplectic
along $\textmd{Im}(\gamma)$, and hence $X_{h_{\mathcal {K}_\mu}}\cdot
\bar{\varepsilon}_\mu= T\bar{\varepsilon}_\mu \cdot X_{h_{\mathcal {K}_\mu} \cdot
\bar{\varepsilon}_\mu}, $ along $\bar{\varepsilon}_\mu$, and hence
$\tau_{\mathcal{K}_\mu}\cdot X_{h_{\mathcal {K}_\mu}} \cdot
\bar{\varepsilon}_\mu= \tau_{\mathcal{K}_\mu}\cdot T\bar{\lambda}_\mu
\cdot X_{h_{\mathcal {K}_\mu} \cdot \bar{\varepsilon}_\mu}, $ along $\bar{\varepsilon}_\mu$, because
$\textmd{Im}(T\bar{\gamma}_\mu)\subset \mathcal{K}_\mu. $ Then we
have that
\begin{align*}
& \omega_{\mathcal{K}_\mu}(T\bar{\gamma}_\mu \cdot X_H^\varepsilon, \;
\tau_{\mathcal{K}_\mu}\cdot T\pi_\mu \cdot w)-
\omega_{\mathcal{K}_\mu}(X_{\mathcal{K}_\mu}\cdot \bar{\varepsilon}_\mu,
\; \tau_{\mathcal{K}_\mu} \cdot T\pi_\mu \cdot w) \nonumber \\
& = - \omega_{\mathcal{K}_\mu}(\tau_{\mathcal{K}_\mu}\cdot X_{h_{\mathcal {K}_\mu}} \cdot
\bar{\varepsilon}_\mu, \; T\bar{\lambda}_\mu \cdot w)
+\omega_{\mathcal{K}_\mu}(T\bar{\lambda}_\mu \cdot X_H \cdot
\varepsilon, \; T\bar{\lambda}_\mu \cdot w)\\
& = \omega_{\mathcal{K}_\mu}(T\bar{\lambda}_\mu \cdot X_H \cdot
\varepsilon- \tau_{\mathcal{K}_\mu}\cdot T\bar{\varepsilon}_\mu
\cdot X_{h_{\mathcal {K}_\mu} \cdot \bar{\varepsilon}_\mu}, \; T\bar{\lambda}_\mu \cdot w).
\end{align*}
Because the $\mathbf{J}$-nonholonomic $R_p$-reduced distributional two-form
$\omega_{\mathcal{K}_\mu}$ is non-degenerate, it follows that the equation
$T\bar{\gamma}_\mu\cdot X_H^\varepsilon = X_{\mathcal{K}_\mu}\cdot
\bar{\varepsilon}_\mu, $ is equivalent to the equation $T\bar{\lambda}_\mu\cdot X_H \cdot \varepsilon
= \tau_{\mathcal{K}_\mu}\cdot T\bar{\varepsilon}_\mu \cdot X_{h_{\mathcal {K}_\mu}
\cdot \bar{\varepsilon}_\mu}. $
Thus, $\varepsilon$ and $\bar{\varepsilon}_\mu$ satisfy the equation
$T\bar{\lambda}_\mu\cdot X_H \cdot \varepsilon
= \tau_{\mathcal{K}_\mu}\cdot T\bar{\varepsilon}_\mu \cdot X_{h_{\mathcal {K}_\mu}
\cdot \bar{\varepsilon}_\mu}, $ if and only if they satisfy
the Type II of Hamilton-Jacobi equation $T\bar{\gamma}_\mu\cdot X_H^\varepsilon = X_{\mathcal{K}_\mu}\cdot
\bar{\varepsilon}_\mu .$
\hskip 0.3cm $\blacksquare$\\

\begin{rema}
If a $\mathbf{J}$-nonholonomic point reducible Hamiltonian system we considered has not any constrains,
in this case, the $\mathbf{J}$-nonholonomic $R_p$-reduced distributional
Hamiltonian system is just the Marsden-Weinstein reduced Hamiltonian system itself.
From the above Type I and Type II of Hamilton-Jacobi theorems, that is,
Theorem 5.2 and Theorem 5.3, we can get the Theorem 3.3
and Theorem 3.4 in Wang \cite{wa17}.
It shows that Theorem 5.2 and Theorem 5.3 can be regarded as an extension of two types of
Hamilton-Jacobi theorem for the Marsden-Weinstein reduced Hamiltonian system
given in \cite{wa17} to the nonholonomic context.
\end{rema}

\begin{rema}
It is worthy of note that the formulations of Type I and Type II of Hamilton-Jacobi
equation for a $\mathbf{J}$-nonholonomic $R_p$-reduced distributional
Hamiltonian system, given by Theorem 5.2 and Theorem 5.3, have
more extensive sense, because, in general, the one-form $\gamma$ is not
given by a generating function of a symplectic map. When $\gamma$ is a solution
of the classical Hamilton-Jacobi equation, that is, $X_H\cdot
\gamma=0, $ then $X_H^\gamma= T\pi_Q\cdot X_H\cdot \gamma=0,$ and
hence from the Type I of Hamilton-Jacobi equation, we have that
$X_{\mathcal{K}_\mu}\cdot \bar{\gamma}_\mu= T\bar{\gamma}_\mu \cdot
X_H^\gamma=0, $ which shows that the
dynamical vector field of the $\mathbf{J}$-nonholonomic $R_p$-reduced
distributional Hamiltonian system
$(\mathcal{K}_\mu,\omega_{\mathcal{K}_\mu},h_{\mathcal {K}_\mu})$ is degenerate along
$\bar{\gamma}_\mu$. The equation $X_{\mathcal{K}_\mu}\cdot
\bar{\gamma}_\mu=0$ is called the classical Hamilton-Jacobi equation
for the $\mathbf{J}$-nonholonomic $R_p$-reduced distributional Hamiltonian
system $(\mathcal{K}_\mu,\omega_{\mathcal{K}_\mu},h_{\mathcal {K}_\mu})$. In addition,
for a symplectic map $\varepsilon: T^* Q \rightarrow T^* Q $, if
$X_H\cdot \varepsilon=0, $ then from the Type II of Hamilton-Jacobi
equation, we have that $X_{\mathcal{K}_\mu}\cdot
\bar{\varepsilon}_\mu= T\bar{\gamma}_\mu\cdot X_H^\varepsilon=0. $
But, from the equation $T\bar{\lambda}_\mu\cdot X_H \cdot
\varepsilon = \tau_{\mathcal{K}_\mu}\cdot T\bar{\varepsilon}_\mu
\cdot X_{h_{\mathcal {K}_\mu} \cdot \bar{\varepsilon}_\mu}, $ we know that the equation
$X_{\mathcal{K}_\mu}\cdot \bar{\varepsilon}_\mu=0 $ is not
equivalent to the equation $X_{h_{\mathcal {K}_\mu} \cdot \bar{\varepsilon}_\mu}=0.$
\end{rema}

For a given $\mathbf{J}$-nonholonomic regular point reducible Hamiltonian system
$(T^*Q,G,\omega,\mathbf{J},\mathcal{D},H)$ with an associated $\mathbf{J}$-
nonholonomic $R_p$-reduced distributional Hamiltonian system
$(\mathcal{K}_\mu,\omega_{\mathcal{K}_\mu},h_{\mathcal {K}_\mu})$, we know that the
Hamiltonian vector field $X_{H}$ and the $\mathbf{J}$-
nonholonomic $R_p$-reduced dynamical vector field $X_{h_{\mathcal {K}_\mu}}$ are $\pi_\mu$-related,
that is, $X_{h_{\mathcal {K}_\mu}}\cdot \pi_\mu=T\pi_\mu\cdot X_{H}\cdot i_\mu.$ Then
we can prove the following Theorem 5.6, which states the
relationship between the solutions of Type II of Hamilton-Jacobi equations and
$\mathbf{J}$-nonholonomic regular point reduction.

\begin{theo}
For a given $\mathbf{J}$-nonholonomic regular point reducible Hamiltonian system
$(T^*Q,G,\omega, \mathbf{J}, \\ \mathcal{D},H)$ with an associated
$\mathbf{J}$-nonholonomic $R_p$-reduced distributional Hamiltonian system
$(\mathcal{K}_\mu,\omega_{\mathcal{K}_\mu},h_{\mathcal {K}_\mu})$, assume that $\gamma:
Q \rightarrow T^*Q$ is an one-form on $Q$, and $\lambda=\gamma \cdot
\pi_{Q}: T^* Q \rightarrow T^* Q, $ and
$\varepsilon: T^* Q \rightarrow T^* Q $ is a symplectic map. Moreover, assume
that $\mu \in \mathfrak{g}^\ast $ is a regular value of the momentum
map $\mathbf{J}$, and $\textmd{Im}(\gamma)\subset \mathcal{M} \cap
\mathbf{J}^{-1}(\mu), $ and it is $G_\mu$-invariant,
and $\varepsilon$ is $G_\mu$-invariant and
$\varepsilon(\mathbf{J}^{-1}(\mu)) \subset \mathbf{J}^{-1}(\mu). $
Denote $\bar{\gamma}_\mu=\pi_\mu(\gamma): Q \rightarrow \mathcal{M}_\mu $,
and $ \textmd{Im}(T\bar{\gamma}_\mu)\subset \mathcal{K}_\mu, $ and
$\bar{\lambda}_\mu=\pi_\mu(\lambda): \mathbf{J}^{-1}(\mu) (\subset T^*Q) \rightarrow (T^* Q)_\mu $,
and $\bar{\varepsilon}_\mu=\pi_\mu(\varepsilon): \mathbf{J}^{-1}(\mu) (\subset T^*Q) \rightarrow (T^* Q)_\mu $.
Then $\varepsilon$
is a solution of the Type II of Hamilton-Jacobi equation $T\gamma\cdot
X_H^\varepsilon= X_{\mathcal{K}}\cdot \varepsilon, $ for the distributional
Hamiltonian system $(\mathcal{K},\omega_{\mathcal {K}},H_{\mathcal {K}})$, if and
only if $\varepsilon$ and $\bar{\varepsilon}_\mu $ satisfy the Type II of Hamilton-Jacobi
equation $T\bar{\gamma}_\mu\cdot X_H^\varepsilon=
X_{\mathcal{K}_\mu}\cdot \bar{\varepsilon}_\mu, $ for the
$\mathbf{J}$-nonholonomic $R_p$-reduced distributional Hamiltonian system
$(\mathcal{K}_\mu,\omega_{\mathcal{K}_\mu},h_{\mathcal {K}_\mu})$.
\end{theo}

\noindent{\bf Proof: } Note that
$\textmd{Im}(\gamma)\subset \mathcal{M} \cap \mathbf{J}^{-1}(\mu), $
and it is $G_\mu$-invariant,
as well as $\textmd{Im}(T\bar{\gamma}_\mu)\subset \mathcal{K}_\mu, $
in this case, $\omega_{\mathcal{K}_\mu}\cdot
\tau_{\mathcal{K}_\mu}= \tau_{\mathcal{K}_\mu}\cdot
\omega_{\mathcal{M}_\mu}= \tau_{\mathcal{K}_\mu}\cdot
i_{\mathcal{M}_\mu}^* \cdot \omega_\mu, $ along
$\textmd{Im}(T\bar{\gamma}_\mu)$,
and
$\pi_\mu^*\omega_\mu=
i_\mu^*\omega= \omega, $ along $\textmd{Im}(\gamma)$,
and $\tau_{\mathcal{K}_\mu}\cdot T\bar{\gamma}_\mu
= T\bar{\gamma}_\mu, $ and $\tau_{\mathcal{K}_\mu}\cdot X_{h_{\mathcal{K}_\mu}} = X_{\mathcal{K}_\mu}. $
Since the Hamiltonian vector field $X_{H}$
and the $\mathbf{J}$-nonholonomic $R_p$-reduced dynamical vector field
$X_{h_{\mathcal{K}_\mu}}$ are $\pi_\mu$-related, that is, $X_{h_{\mathcal{K}_\mu}}\cdot \pi_\mu=
T\pi_\mu\cdot X_{H}\cdot i_\mu, $
using the $\mathbf{J}$-nonholonomic $R_p$-reduced distributional two-form
$\omega_{\mathcal{K}_\mu}$, we have that
\begin{align*}
& \omega_{\mathcal{K}_\mu}(T\bar{\gamma}_\mu \cdot X_H^\varepsilon- X_{\mathcal{K}_\mu}\cdot \bar{\varepsilon}_\mu, \;
\tau_{\mathcal{K}_\mu}\cdot T\pi_\mu \cdot w) \\
& = \omega_{\mathcal{K}_\mu}(T\bar{\gamma}_\mu \cdot X_H^\varepsilon, \;
\tau_{\mathcal{K}_\mu}\cdot T\pi_\mu \cdot w)-
\omega_{\mathcal{K}_\mu}(X_{\mathcal{K}_\mu}\cdot \bar{\varepsilon}_\mu,
\; \tau_{\mathcal{K}_\mu} \cdot T\pi_\mu \cdot w) \\
& = \omega_{\mathcal{K}_\mu}(\tau_{\mathcal{K}_\mu}\cdot T\bar{\gamma}_\mu \cdot X_H^\varepsilon, \;
\tau_{\mathcal{K}_\mu}\cdot T\pi_\mu \cdot w)-
\omega_{\mathcal{K}_\mu}(\tau_{\mathcal{K}_\mu}\cdot X_{h_{\mathcal{K}_\mu}} \cdot \pi_\mu \cdot
\varepsilon, \; \tau_{\mathcal{K}_\mu} \cdot T\pi_\mu \cdot w) \\
& = \omega_{\mathcal{K}_\mu}\cdot \tau_{\mathcal{K}_\mu}(T\pi_\mu \cdot T\gamma \cdot X_H^\varepsilon, \; T\pi_\mu \cdot w)
-\omega_{\mathcal{K}_\mu}\cdot \tau_{\mathcal{K}_\mu}(T\pi_\mu \cdot X_H\cdot \varepsilon, \; T\pi_\mu \cdot w)\\
& = \tau_{\mathcal{K}_\mu}\cdot
i_{\mathcal{M}_\mu}^* \cdot \omega_\mu(T\pi_\mu \cdot T\gamma \cdot X_H^\varepsilon, \; T\pi_\mu \cdot w)
-\tau_{\mathcal{K}_\mu}\cdot
i_{\mathcal{M}_\mu}^* \cdot \omega_\mu(T\pi_\mu \cdot X_H\cdot \varepsilon, \; T\pi_\mu \cdot w)\\
& = \tau_{\mathcal{K}_\mu}\cdot
i_{\mathcal{M}_\mu}^* \cdot \pi_\mu^*\omega_\mu(T\gamma \cdot X_H^\varepsilon, \; w)
-\tau_{\mathcal{K}_\mu}\cdot
i_{\mathcal{M}_\mu}^* \cdot \pi_\mu^*\omega_\mu(X_H\cdot \varepsilon, \; w)\\
& = \tau_{\mathcal{K}_\mu}\cdot
i_{\mathcal{M}_\mu}^* \cdot \omega(T\gamma \cdot X_H^\varepsilon, \; w)
-\tau_{\mathcal{K}_\mu}\cdot
i_{\mathcal{M}_\mu}^* \cdot \omega(X_H\cdot \varepsilon, \; w).
\end{align*}
In the case we considered that $\tau_{\mathcal{K}_\mu}\cdot
i_{\mathcal{M}_\mu}^* \cdot \omega=\tau_{\mathcal{K}}\cdot i_{\mathcal{M}}^* \cdot
\omega= \omega_{\mathcal{K}}\cdot \tau_{\mathcal{K}}, $
and
$\tau_{\mathcal{K}}\cdot T\gamma =T\gamma, \; \tau_{\mathcal{K}} \cdot X_H= X_{\mathcal{K}}$,
since $\textmd{Im}(\gamma)\subset
\mathcal{M}, $ and $\textmd{Im}(T\gamma)\subset \mathcal{K}. $
Thus, we have that
\begin{align*}
& \omega_{\mathcal{K}_\mu}(T\bar{\gamma}_\mu \cdot X_H^\varepsilon- X_{\mathcal{K}_\mu}\cdot \bar{\varepsilon}_\mu, \;
\tau_{\mathcal{K}_\mu}\cdot T\pi_\mu \cdot w) \\
& = \omega_{\mathcal{K}}\cdot \tau_{\mathcal{K}}(T\gamma \cdot X_H^\varepsilon, \; w)
-\omega_{\mathcal{K}}\cdot \tau_{\mathcal{K}}(X_H\cdot \varepsilon, \; w)\\
& = \omega_{\mathcal{K}}(\tau_{\mathcal{K}}\cdot T\gamma \cdot X_H^\varepsilon, \; \tau_{\mathcal{K}}\cdot w)
-\omega_{\mathcal{K}}(\tau_{\mathcal{K}}\cdot X_H\cdot \varepsilon, \; \tau_{\mathcal{K}}\cdot w)\\
& = \omega_{\mathcal{K}}(T\gamma \cdot X_H^\varepsilon, \; \tau_{\mathcal{K}}\cdot w)
-\omega_{\mathcal{K}}(X_{\mathcal{K}}\cdot \varepsilon, \; \tau_{\mathcal{K}}\cdot w)\\
& = \omega_{\mathcal{K}}(T\gamma \cdot X_H^\varepsilon- X_{\mathcal{K}}\cdot \varepsilon, \; \tau_{\mathcal{K}}\cdot w).
\end{align*}
Because the distributional two-form $\omega_{\mathcal{K}}$ and
the $\mathbf{J}$-nonholonomic $R_p$-reduced distributional
two-form $\omega_{\mathcal{K}_\mu}$ are both non-degenerate,
it follows that the equation
$T\bar{\gamma}_\mu\cdot X_H^\varepsilon=
X_{\mathcal{K}_\mu}\cdot \bar{\varepsilon}_\mu, $ is equivalent to the equation
$T\gamma\cdot X_H^\varepsilon= X_{\mathcal{K}}\cdot \varepsilon$. Thus,
$\varepsilon$ is a solution of the Type II of Hamilton-Jacobi equation
$T\gamma\cdot X_H^\varepsilon= X_{\mathcal{K}}\cdot \varepsilon, $ for the distributional
Hamiltonian system $(\mathcal{K},\omega_{\mathcal {K}},H_{\mathcal {K}})$, if and only if
$\varepsilon$ and $\bar{\varepsilon}_\mu $ satisfy the Type II of Hamilton-Jacobi
equation $T\bar{\gamma}_\mu\cdot X_H^\varepsilon=
X_{\mathcal{K}_\mu}\cdot \bar{\varepsilon}_\mu, $ for the
$\mathbf{J}$-nonholonomic $R_p$-reduced distributional Hamiltonian system
$(\mathcal{K}_\mu,\omega_{\mathcal{K}_\mu},h_{\mathcal{K}_\mu})$.  \hskip 0.3cm
$\blacksquare$

\begin{rema}
If $(T^\ast Q, \omega)$ is a connected symplectic manifold, and
$\mathbf{J}:T^\ast Q\rightarrow \mathfrak{g}^\ast$ is a
non-equivariant momentum map with a non-equivariance group
one-cocycle $\sigma: G\rightarrow \mathfrak{g}^\ast$,
which is defined by $\sigma(g):=\mathbf{J}(g\cdot
z)-\operatorname{Ad}^\ast_{g^{-1}}\mathbf{J}(z)$, where $g\in G$ and
$z\in T^\ast Q$. Then we know that $\sigma$ produces a new affine
action $\Theta: G\times \mathfrak{g}^\ast \rightarrow
\mathfrak{g}^\ast $ defined by
$\Theta(g,\mu):=\operatorname{Ad}^\ast_{g^{-1}}\mu + \sigma(g)$,
where $\mu \in \mathfrak{g}^\ast$, with respect to which the given
momentum map $\mathbf{J}$ is equivariant. Assume that $G$ acts
freely and properly on $T^\ast Q$, and $\tilde{G}_\mu$ denotes the
isotropy subgroup of $\mu \in \mathfrak{g}^\ast$ relative to this
affine action $\Theta$ and $\mu$ is a regular value of $\mathbf{J}$.
Then the quotient space $(T^\ast
Q)_\mu=\mathbf{J}^{-1}(\mu)/\tilde{G}_\mu$ is also a symplectic
manifold with symplectic form $\omega_\mu$ uniquely characterized by
$(5.1)$, see Ortega and Ratiu \cite{orra04}. In this case,
we can also define the $\mathbf{J}$-nonholonomic regular point reducible
Hamiltonian system $(T^*Q,G,\omega,\mathbf{J},\mathcal{D},H)$ with an associated
$\mathbf{J}$-nonholonomic $R_p$-reduced distributional Hamiltonian system
$(\mathcal{K}_\mu,\omega_{\mathcal{K}_\mu},h_{\mathcal{K}_\mu})$,
and prove the Type I and Type II of
the Hamilton-Jacobi theorem for the $\mathbf{J}$-nonholonomic $R_p$-reduced
distributional Hamiltonian system
$(\mathcal{K}_\mu,\omega_{\mathcal{K}_\mu}, h_{\mathcal{K}_\mu})$ by using the above similar way,
in which the $\mathbf{J}$-nonholonomic $R_p$-reduced space
$(\mathcal{K}_\mu,\omega_{\mathcal{K}_\mu})$
is determined by the affine action and $\mathbf{J}$-nonholonomic regular point reduction.
\end{rema}

\subsection{Hamilton-Jacobi equations in the case compatible with regular orbit reduction}

In this subsection, for a nonholonomic Hamiltonian
system with symmetry and momentum map
$(T^*Q,G,\omega,\mathbf{J},\mathcal{D},H)$,
where $\omega$ is the canonical
symplectic form on $T^* Q$, and $\mathcal{D}\subset TQ$ is a
$\mathcal{D}$-completely and $\mathcal{D}$-regularly nonholonomic
constraint of the system, and $\mathcal{D}$ and $H$ are both
$G$-invariant, we first give the $\mathbf{J}$-nonholonomic regular orbit reduction
of the system compatible with regular orbit reduction,
and a $\mathbf{J}$-nonholonomic $R_o$-reduced distribution $\mathcal{K}_{\mathcal{O}_\mu}$,
an associated non-degenerate and $\mathbf{J}$-nonholonomic $R_o$-reduced
distributional two-form $\omega_{\mathcal{K}_{\mathcal{O}_\mu}}$,
which is induced by the canonical symplectic form $\omega$ on $T^*Q$,
and a $\mathbf{J}$-nonholonomic $R_o$-reduced distributional Hamiltonian system,
where the "regular orbit reduced" is simply written as $R_o$-reduced.
Then we derive precisely
the geometric constraint conditions of the $\mathbf{J}$-nonholonomic $R_o$-reduced distributional two-form
$\omega_{\mathcal{K}_{\mathcal{O}_\mu}}$ for the nonholonomic reducible dynamical vector field,
that is, the two types of Hamilton-Jacobi equation for the
$\mathbf{J}$-nonholonomic $R_o$-reduced distributional Hamiltonian system, which are an
extension of the above two types of Hamilton-Jacobi equation for the distributional
Hamiltonian system under $\mathbf{J}$-nonholonomic regular orbit reduction.\\

At first, we need to give
carefully a geometric formulation of the $\mathbf{J}$-nonholonomic
$R_o$-reduced distributional Hamiltonian system, by using momentum map and the
nonholonomic reduction compatible with regular orbit reduction.
Now, we assume that the
6-tuple $(T^*Q,G,\omega,\mathbf{J},\mathcal{D},H)$ is a
$\mathcal{D}$-completely and $\mathcal{D}$-regularly nonholonomic
Hamiltonian system with symmetry and momentum map, and the Lie group
$G$ acts smoothly by the left on $Q$. For the cotangent lifted
left action $\Phi^{T^\ast}:G\times T^\ast Q\rightarrow T^\ast Q$, assume that
it is free, proper and symplectic, and the action admits an
$\operatorname{Ad}^\ast$-equivariant momentum map $\mathbf{J}:T^\ast
Q\rightarrow \mathfrak{g}^\ast$. Let $\mu\in \mathfrak{g}^\ast$ be a
regular value of the momentum map $\mathbf{J}$ and
$\mathcal{O}_\mu=G\cdot \mu\subset \mathfrak{g}^\ast$ be the
$G$-orbit of the coadjoint $G$-action through the point $\mu$. Since
$G$ acts freely, properly and symplectically on $T^\ast Q$, then the
quotient space $(T^\ast Q)_{\mathcal{O}_\mu}=
\mathbf{J}^{-1}(\mathcal{O}_\mu)/G$ is a regular quotient symplectic
manifold with the symplectic form $\omega_{\mathcal{O}_\mu}$
uniquely characterized by the relation
\begin{equation}i_{\mathcal{O}_\mu}^\ast \omega=\pi_{\mathcal{O}_{\mu}}^\ast
\omega_{\mathcal{O}
_\mu}+\mathbf{J}_{\mathcal{O}_\mu}^\ast\omega_{\mathcal{O}_\mu}^+,
\label{4.1}\end{equation}
where $\mathbf{J}_{\mathcal{O}_\mu}$ is
the restriction of the momentum map $\mathbf{J}$ to
$\mathbf{J}^{-1}(\mathcal{O}_\mu)$, that is,
$\mathbf{J}_{\mathcal{O}_\mu}=\mathbf{J}\cdot i_{\mathcal{O}_\mu}$
and $\omega_{\mathcal{O}_\mu}^+$ is the $+$-symplectic structure on
the orbit $\mathcal{O}_\mu$ given by
\begin{equation}\omega_{\mathcal{O}_\mu}^
+(\nu)(\xi_{\mathfrak{g}^\ast}(\nu),\eta_{\mathfrak{g}^\ast}(\nu))
=<\nu,[\xi,\eta]>,\;\; \forall\;\nu\in\mathcal{O}_\mu, \;
\xi,\eta\in \mathfrak{g}.
\label{3.3}\end{equation}
The maps
$i_{\mathcal{O}_\mu}:\mathbf{J}^{-1}(\mathcal{O}_\mu)\rightarrow
T^\ast Q$ and
$\pi_{\mathcal{O}_\mu}:\mathbf{J}^{-1}(\mathcal{O}_\mu)\rightarrow
(T^\ast Q)_{\mathcal{O}_\mu}$ are natural injection and the
projection, respectively. The pair $((T^\ast
Q)_{\mathcal{O}_\mu},\omega_{\mathcal{O}_\mu})$ is called the regular
orbit reduced symplectic space of $(T^\ast Q,\omega)$ at $\mu$.\\

Let $H:T^\ast Q\rightarrow \mathbb{R}$ be a $G$-invariant
Hamiltonian, the flow $F_t$ of the Hamiltonian vector field $X_H$
leaves the connected components of
$\mathbf{J}^{-1}(\mathcal{O}_\mu)$ invariant and commutes with the
$G$-action, so it induces a flow $f_t^{\mathcal{O}_\mu}$ on $(T^\ast
Q)_{\mathcal{O}_\mu}$, defined by $f_t^{\mathcal{O}_\mu}\cdot
\pi_{\mathcal{O}_\mu}=\pi_{\mathcal{O}_\mu} \cdot F_t\cdot
i_{\mathcal{O}_\mu}$, and the vector field $X_{h_{\mathcal{O}_\mu}}$
generated by the flow $f_t^{\mathcal{O}_\mu}$ on $((T^\ast
Q)_{\mathcal{O}_\mu},\omega_{\mathcal{O}_\mu})$ is Hamiltonian with
the associated $R_o$-reduced Hamiltonian function
$h_{\mathcal{O}_\mu}:(T^\ast Q)_{\mathcal{O}_\mu}\rightarrow
\mathbb{R}$ defined by $h_{\mathcal{O}_\mu}\cdot
\pi_{\mathcal{O}_\mu}= H\cdot i_{\mathcal{O}_\mu}$, and the
Hamiltonian vector fields $X_H$ and $X_{h_{\mathcal{O}_\mu}}$ are
$\pi_{\mathcal{O}_\mu}$-related. \\

In the same way, by using the
Legendre transformation $\mathcal{F}L: TQ \rightarrow T^*Q $, we can
define the constraint submanifold
$\mathcal{M}=\mathcal{F}L(\mathcal{D})\subset T^*Q $ and the
distribution $\mathcal{F}=(T\pi_Q)^{-1}(\mathcal{D})$,
and $\mathcal{K}=\mathcal{F} \cap T\mathcal{M}$. Moreover, we can
also define the distributional two-form $\omega_\mathcal{K}$, a vector
field $X_\mathcal{K}$ and the function $H_\mathcal{K}$, such that
$\mathbf{i}_{X_\mathcal{K}}\omega_\mathcal{K}=\mathbf{d}H_\mathcal{K}$.
Since $\mathcal{D}\subset TQ$ is a $G$-invariant distribution, and the
Legendre transformation $\mathcal{F}L: TQ \rightarrow T^*Q$ is a
fiber-preserving map, then
$\mathcal{M}=\mathcal{F}L(\mathcal{D})\subset T^*Q$ is
$G$-invariant. For a regular value $\mu\in\mathfrak{g}^\ast$ of the
momentum map $\mathbf{J}$, $\mathcal{O}_\mu=G\cdot \mu\subset \mathfrak{g}^\ast$ is the
$G$-orbit of the coadjoint $G$-action through the point $\mu$,
we shall assume that the constraint submanifold $\mathcal{M}$
is clean intersection with $\mathbf{J}^{-1}(\mathcal{O}_\mu)$, that is,
$\mathcal{M} \cap \mathbf{J}^{-1}(\mathcal{O}_\mu)\neq \emptyset$.
It follows that the quotient
space $\mathcal{M}_{\mathcal{O}_\mu} =(\mathcal{M}\cap \mathbf{J}^{-1}(\mathcal{O}_\mu))
/G \subset (T^\ast Q)_{\mathcal{O}_\mu}$ of the $G$-orbit in
$\mathcal{M}\cap \mathbf{J}^{-1}(\mathcal{O}_\mu)$, is a smooth manifold with
projection $\pi_{\mathcal{O}_\mu}: \mathcal{M}\cap \mathbf{J}^{-1}(\mathcal{O}_\mu)
\rightarrow \mathcal{M}_{\mathcal{O}_\mu}$ which is a surjective submersion.
Denote $i_{\mathcal{M}_{\mathcal{O}_\mu}}: \mathcal{M}_{\mathcal{O}_\mu}\rightarrow
(T^*Q)_{\mathcal{O}_\mu}, $ and $\omega_{\mathcal{M}_{\mathcal{O}_\mu}}= i_{\mathcal{M}_{\mathcal{O}_\mu}}^*
\omega_{\mathcal{O}_\mu} $, that is, the symplectic form
$\omega_{\mathcal{M}_{\mathcal{O}_\mu}}$ is induced from the $R_o$-reduced symplectic
form $\omega_{\mathcal{O}_\mu}$ on $(T^* Q)_{\mathcal{O}_\mu}$,
where $i_{\mathcal{M}_{\mathcal{O}_\mu}}^*:
T^*(T^*Q)_{\mathcal{O}_\mu} \rightarrow T^*\mathcal{M}_{\mathcal{O}_\mu}.$
Moreover, the distribution $\mathcal{F}$ pushes down to a distribution
$\mathcal{F}_{\mathcal{O}_\mu}= T\pi_{\mathcal{O}_\mu}\cdot \mathcal{F}$
on $(T^\ast Q)_{\mathcal{O}_\mu}$,
and we define $\mathcal{K}_{\mathcal{O}_\mu}=\mathcal{F}_{\mathcal{O}_\mu} \cap
T\mathcal{M}_{\mathcal{O}_\mu}$. Assume that $\omega_{\mathcal{K}_{\mathcal{O}_\mu}}=
\tau_{\mathcal{K}_{\mathcal{O}_\mu}}\cdot \omega_{\mathcal{M}_{\mathcal{O}_\mu}}$ is the
restriction of the symplectic form $\omega_{\mathcal{M}_{\mathcal{O}_\mu}}$ on
$T^*\mathcal{M}_{\mathcal{O}_\mu}$ fibrewise to
the distribution $\mathcal{K}_{\mathcal{O}_\mu}$,
where $\tau_{\mathcal{K}_{\mathcal{O}_\mu}}$ is the restriction map to distribution
$\mathcal{K}_{\mathcal{O}_\mu}$. \\

From the above construction, we know that
$\omega_{\mathcal{K}_{\mathcal{O}_\mu}}$ is non-degenerate, and is called
as a $\mathbf{J}$-nonholonomic $R_o$-reduced distributional two-form to
avoid any confusion. Because $\omega_{\mathcal{K}_{\mathcal{O}_\mu}}$ is
non-degenerate as a bilinear form on each fibre of
$\mathcal{K}_{\mathcal{O}_\mu}$, there exists a vector field $X_{\mathcal{K}_{\mathcal{O}_\mu}}$
on $\mathcal{M}_{\mathcal{O}_\mu}$, which takes values in the constraint
distribution $\mathcal{K}_{\mathcal{O}_\mu}$,
such that the $\mathbf{J}$-nonholonomic $R_o$-reduced distributional
Hamiltonian equation holds, that is,
$\mathbf{i}_{X_{\mathcal{K}_{\mathcal{O}_\mu}}}\omega_{\mathcal{K}_{\mathcal{O}_\mu}}
=\mathbf{d}h_{\mathcal{K}_{\mathcal{O}_\mu}}$,
if the admissibility condition $\mathrm{dim}\mathcal{M}_{\mathcal{O}_\mu}=
\mathrm{rank}\mathcal{F}_{\mathcal{O}_\mu}$ and the compatibility condition
$T\mathcal{M}_{\mathcal{O}_\mu}\cap \mathcal{F}_{\mathcal{O}_\mu}^\bot= \{0\}$ hold, where
$\mathcal{F}_{\mathcal{O}_\mu}^\bot$ denotes the symplectic orthogonal of
$\mathcal{F}_{\mathcal{O}_\mu}$ with respect to the $R_o$-reduced symplectic form
$\omega_{\mathcal{O}_\mu}$, and $\mathbf{d}h_{\mathcal{K}_{\mathcal{O}_\mu}}$ is the restriction
of $\mathbf{d}h_{\mathcal{M}_{\mathcal{O}_\mu}}$ to $\mathcal{K}_{\mathcal{O}_\mu}$, and
the function $h_{\mathcal{K}_{\mathcal{O}_\mu}}$ satisfies
$\mathbf{d}h_{\mathcal{K}_{\mathcal{O}_\mu}}= \tau_{\mathcal{K}_{\mathcal{O}_\mu}}\cdot \mathbf{d}h_{\mathcal{M}_{\mathcal{O}_\mu}} $,
and $h_{\mathcal{M}_{\mathcal{O}_\mu}}
= \tau_{\mathcal{M}_{\mathcal{O}_\mu}}\cdot h_{\mathcal{O}_\mu}$ is the
restriction of $h_{\mathcal{O}_\mu}$ to $\mathcal{M}_{\mathcal{O}_\mu}$,
and $h_{\mathcal{O}_\mu}$ is the $R_o$-reduced
Hamiltonian function $h_{\mathcal{O}_\mu}: (T^* Q)_{\mathcal{O}_\mu} \rightarrow \mathbb{R}$ defined
by $h_{\mathcal{O}_\mu}\cdot \pi_{\mathcal{O}_\mu}= H\cdot i_{\mathcal{O}_\mu}$.
Thus, the geometrical formulation
of the $\mathbf{J}$-nonholonomic $R_o$-reduced distributional Hamiltonian
system may be summarized as follows.

\begin{defi} ($\mathbf{J}$-Nonholonomic $R_o$-reduced Distributional Hamiltonian System)
Assume that the 6-tuple $(T^*Q,G,\omega,\mathbf{J},\mathcal{D},H)$
is a nonholonomic Hamiltonian system with symmetry and momentum map,
where $\omega$ is the canonical symplectic form on $T^* Q$,
and $\mathcal{D}\subset TQ$ is a $\mathcal{D}$-completely and
$\mathcal{D}$-regularly nonholonomic constraint of the system, and
$\mathcal{D}$ and $H$ are both $G$-invariant. For a regular value
$\mu\in\mathfrak{g}^\ast$ of the momentum map $\mathbf{J}$,
$\mathcal{O}_\mu=G\cdot \mu\subset \mathfrak{g}^\ast$ is the
$G$-orbit of the coadjoint $G$-action through the point $\mu$,
assume that there exists a
$\mathbf{J}$-nonholonomic $R_o$-reduced distribution
$\mathcal{K}_{\mathcal{O}_\mu}$, an associated non-degenerate
and $\mathbf{J}$-nonholonomic $R_o$-reduced distributional two-form
$\omega_{\mathcal{K}_{\mathcal{O}_\mu}}$
and a vector field $X_{\mathcal {K}_{\mathcal{O}_\mu}}$
on the $\mathbf{J}$-nonholonomic $R_o$-reduced constraint submanifold
$\mathcal{M}_{\mathcal{O}_\mu}=(\mathcal{M}\cap \mathbf{J}^{-1}(\mathcal{O}_\mu)) /G, $
where $\mathcal{M}=\mathcal{F}L(\mathcal{D}),$ and $\mathcal{M}\cap
\mathbf{J}^{-1}({\mathcal{O}_\mu})\neq \emptyset, $ such that the
$\mathbf{J}$-nonholonomic $R_o$-reduced distributional Hamiltonian
equation $\mathbf{i}_{X_{\mathcal{K}_{\mathcal{O}_\mu}}}\omega_{\mathcal{K}_{\mathcal{O}_\mu}} =
\mathbf{d}h_{\mathcal{K}_{\mathcal{O}_\mu}}$ holds, where
$\mathbf{d}h_{\mathcal{K}_{\mathcal{O}_\mu}}$ is the restriction of
$\mathbf{d}h_{\mathcal{M}_{\mathcal{O}_\mu}}$ to $\mathcal{K}_{\mathcal{O}_\mu}$,
and the function $h_{\mathcal{K}_{\mathcal{O}_\mu}}$ is defined above.
Then the triple
$(\mathcal{K}_{\mathcal{O}_\mu},\omega_{\mathcal {K}_{\mathcal{O}_\mu}},
h_{\mathcal {K}_{\mathcal{O}_\mu}})$ is called a
$\mathbf{J}$-nonholonomic $R_o$-reduced distributional Hamiltonian system
of the system $(T^*Q,G,\omega,\mathbf{J},\mathcal{D},H)$, and $X_{\mathcal
{K}_{\mathcal{O}_\mu}}$ is the dynamical vector field of the
$\mathbf{J}$-nonholonomic $R_o$-reduced distributional Hamiltonian system
$(\mathcal{K}_{\mathcal{O}_\mu},\omega_{{\mathcal{K}_{\mathcal{O}_\mu}}},
h_{{\mathcal{K}_{\mathcal{O}_\mu}}})$.
Under the above circumstances, we refer to
$(T^*Q,G,\omega,\mathbf{J},\mathcal{D},H)$ as
a $\mathbf{J}$-nonholonomic regular orbit reducible Hamiltonian system
with an associated $\mathbf{J}$-nonholonomic $R_o$-reduced
distributional Hamiltonian system
$(\mathcal{K}_{\mathcal{O}_\mu},\omega_{\mathcal{K}_{\mathcal{O}_\mu}},
h_{\mathcal{K}_{\mathcal{O}_\mu}})$.
\end{defi}

Since the non-degenerate and $\mathbf{J}$-nonholonomic $R_o$-reduced
distributional two-form $\omega_{\mathcal{K}_{\mathcal{O}_\mu}}$ may not be symplectic,
and the $\mathbf{J}$-nonholonomic $R_o$-reduced distributional Hamiltonian system \\
$(\mathcal{K}_{\mathcal{O}_\mu},\omega_{\mathcal {K}_{\mathcal{O}_\mu}},
h_{\mathcal{K}_{\mathcal{O}_\mu}})$
may not be yet a Hamiltonian system, and has no yet generating function,
and hence we can not describe the Hamilton-Jacobi equation for a
$\mathbf{J}$-nonholonomic $R_o$-reduced
distributional Hamiltonian system just like as in Theorem 1.1.
But, for a given $\mathbf{J}$-nonholonomic regular orbit reducible Hamiltonian system $(T^*Q,G,\omega,\mathbf{J},\mathcal{D},H)$ with an associated
$\mathbf{J}$-nonholonomic $R_o$-reduced distributional Hamiltonian system
$(\mathcal{K}_{\mathcal{O}_\mu},\omega_{\mathcal {K}_{\mathcal{O}_\mu}},h_{\mathcal{K}_{\mathcal{O}_\mu}})$,
by using Lemma 3.4, we can derive precisely
the geometric constraint conditions of the $\mathbf{J}$-nonholonomic $R_o$-reduced distributional two-form
$\omega_{\mathcal{K}_{\mathcal{O}_\mu}}$ for the nonholonomic reducible dynamical vector field
that is, the two types of Hamilton-Jacobi equation for the
$\mathbf{J}$-nonholonomic $R_o$-reduced distributional Hamiltonian system
$(\mathcal{K}_{\mathcal{O}_\mu},\omega_{\mathcal
{K}_{\mathcal{O}_\mu}},h_{\mathcal{K}_{\mathcal{O}_\mu}})$.
At first, using the fact that the one-form $\gamma: Q
\rightarrow T^*Q $ is closed on $\mathcal{D}$ with respect to
$T\pi_Q: TT^* Q \rightarrow TQ, $
and $\textmd{Im}(\gamma)\subset \mathcal{M} \cap
\mathbf{J}^{-1}(\mathcal{O}_\mu), $ and it is $G$-invariant,
as well as $ \textmd{Im}(T\bar{\gamma}_{\mathcal{O}_\mu})\subset \mathcal{K}_{\mathcal{O}_\mu},$
we can prove the Type I of
Hamilton-Jacobi theorem for the $\mathbf{J}$-nonholonomic $R_o$-reduced distributional
Hamiltonian system. For convenience, the maps involved in the
following theorem and its proof are shown in Diagram-7.
\begin{center}
\hskip 0cm \xymatrix{ \mathbf{J}^{-1}(\mathcal{O}_\mu) \ar[r]^{i_{\mathcal{O}_\mu}}
& T^* Q  \ar[r]^{\pi_Q}
& Q \ar[d]_{X_H^\gamma} \ar[r]^{\gamma}
& T^*Q \ar[d]_{X_H} \ar[r]^{\pi_{\mathcal{O}_\mu}} & (T^* Q)_{\mathcal{O}_\mu} \ar[d]_{X_{h_{\mathcal{O}_\mu}}}
& \mathcal{M}_{\mathcal{O}_\mu}  \ar[l]_{i_{\mathcal{M}_{\mathcal{O}_\mu}}} \ar[d]_{X_{\mathcal{K}_{\mathcal{O}_\mu}}}\\
& T(T^*Q)  & TQ \ar[l]_{T\gamma} & T(T^*Q) \ar[l]^{T\pi_Q} \ar[r]_{T\pi_{\mathcal{O}_\mu}}
& T(T^* Q)_{\mathcal{O}_\mu} \ar[r]^{\tau_{\mathcal{K}_{\mathcal{O}_\mu}}} & \mathcal{K}_{\mathcal{O}_\mu} }
\end{center}
$$\mbox{Diagram-7}$$

\begin{theo} (Type I of Hamilton-Jacobi Theorem for a $\mathbf{J}$-Nonholonomic
$R_o$-reduced Distributional Hamiltonian System) For a given
$\mathbf{J}$-nonholonomic regular orbit reducible Hamiltonian system
$(T^*Q,G,\omega,\mathbf{J},\mathcal{D},H)$ with an associated
$\mathbf{J}$-nonholonomic $R_o$-reduced distributional Hamiltonian system
$(\mathcal{K}_{\mathcal{O}_\mu},\omega_{\mathcal{K}_{\mathcal{O}_\mu}},
h_{\mathcal{K}_{\mathcal{O}_\mu}})$, assume that $\gamma:
Q \rightarrow T^*Q$ is an one-form on $Q$, and $X_H^\gamma =
T\pi_{Q}\cdot X_H \cdot \gamma, $ where $X_{H}$ is the Hamiltonian
vector field of the corresponding unconstrained Hamiltonian system with
symmetry and momentum map $(T^*Q,G,\omega,\mathbf{J},H)$. Moreover,
assume that $\mu\in\mathfrak{g}^\ast$ is a regular value of the momentum
map $\mathbf{J}$, and $\textmd{Im}(\gamma)\subset \mathcal{M} \cap
\mathbf{J}^{-1}(\mathcal{O}_\mu), $ and it is $G$-invariant, and
$\bar{\gamma}_{\mathcal{O}_\mu}
=\pi_{\mathcal{O}_\mu}(\gamma): Q \rightarrow \mathcal{M}_{\mathcal{O}_\mu} $,
and $ \textmd{Im}(T\bar{\gamma}_{\mathcal{O}_\mu})\subset \mathcal{K}_{\mathcal{O}_\mu}. $
If the one-form $\gamma: Q \rightarrow T^*Q $ is closed on $\mathcal{D}$ with respect to
$T\pi_Q: TT^* Q \rightarrow TQ, $ then
$\bar{\gamma}_{\mathcal{O}_\mu}$ is a solution of the
equation $T\bar{\gamma}_{\mathcal{O}_\mu}\cdot
X_H^\gamma= X_{\mathcal{K}_{\mathcal{O}_\mu}}\cdot \bar{\gamma}_{\mathcal{O}_\mu}. $
Here $X_{\mathcal{K}_{\mathcal{O}_\mu}}$ is the dynamical vector field
of the $\mathbf{J}$-nonholonomic $R_o$-reduced
system $(\mathcal{K}_{\mathcal{O}_\mu},
\omega_{\mathcal{K}_{\mathcal{O}_\mu}},h_{\mathcal{K}_{\mathcal{O}_\mu}})$. The equation
$T\bar{\gamma}_{\mathcal{O}_\mu}\cdot X_H^\gamma
= X_{\mathcal{K}_{\mathcal{O}_\mu}}\cdot
\bar{\gamma}_{\mathcal{O}_\mu},$ is called the Type I of Hamilton-Jacobi equation for the
$\mathbf{J}$-nonholonomic $R_o$-reduced distributional Hamiltonian system
$(\mathcal{K}_{\mathcal{O}_\mu},\omega_{\mathcal{K}_{\mathcal{O}_\mu}},h_{\mathcal{K}_{\mathcal{O}_\mu}})$.
\end{theo}

\noindent{\bf Proof: } At first, from Theorem 3.5, we know that
$\gamma$ is a solution of the Hamilton-Jacobi equation
$T\gamma\cdot X_H^\gamma= X_{\mathcal{K}}\cdot \gamma .$ Next, we note that
the $R_o$-reduced symplectic space
$(T^\ast Q)_{\mathcal{O}_\mu}= \mathbf{J}^{-1}(\mathcal{O}_\mu)/G
\cong \mathbf{J}^{-1}(\mu)/G \times \mathcal{O}_\mu, $ with the
symplectic form $\omega_{\mathcal{O}_\mu}$ uniquely characterized by
the relation $i_{\mathcal{O}_\mu}^\ast
\omega=\pi_{\mathcal{O}_{\mu}}^\ast \omega_{\mathcal{O}
_\mu}+\mathbf{J}_{\mathcal{O}_\mu}^\ast\omega_{\mathcal{O}_\mu}^+. $
Since
$\textmd{Im}(\gamma)\subset \mathcal{M} \cap \mathbf{J}^{-1}(\mathcal{O}_\mu), $
and it is $G$-invariant,
in this case for any $V\in TQ, $ and $w\in TT^*Q, $
we have that
$\mathbf{J}_{\mathcal{O}_\mu}^\ast\omega_{\mathcal{O}_\mu}^+(T\gamma
\cdot V, \; w)=0, $ and hence
$\pi_{\mathcal{O}_\mu}^*\omega_{\mathcal{O}_\mu}=
i_{\mathcal{O}_\mu}^*\omega= \omega, $ along $\textmd{Im}(\gamma)$.
On the other hand, because
$\textmd{Im}(T\bar{\gamma}_{\mathcal{O}_\mu})\subset \mathcal{K}_{\mathcal{O}_\mu}, $
then
$\omega_{\mathcal{K}_{\mathcal{O}_\mu}}\cdot
\tau_{\mathcal{K}_{\mathcal{O}_\mu}}=\tau_{\mathcal{K}_{\mathcal{O}_\mu}}\cdot
\omega_{\mathcal{M}_{\mathcal{O}_\mu}}= \tau_{\mathcal{K}_{\mathcal{O}_\mu}}\cdot
i_{\mathcal{M}_{\mathcal{O}_\mu}}^* \cdot \omega_{\mathcal{O}_\mu}, $ along
$\textmd{Im}(T\bar{\gamma}_{\mathcal{O}_\mu})$. Thus, using
the $\mathbf{J}$-nonholonomic $R_o$-reduced distributional two-form
$\omega_{\mathcal{K}_{\mathcal{O}_\mu}}$, from Lemma 3.4, if we take that $v=
\tau_{\mathcal{K}_{\mathcal{O}_\mu}}\cdot T\pi_{\mathcal{O}_\mu} \cdot X_{H}\cdot \gamma
= X_{\mathcal{K}_{\mathcal{O}_\mu}}\cdot
\bar{\gamma}_{\mathcal{O}_\mu} \in \mathcal{K}_{\mathcal{O}_\mu}, $ and for any $w \in
\mathcal{F}, \; T\lambda(w)\neq 0, $ and
$\tau_{\mathcal{K}_{\mathcal{O}_\mu}}\cdot T\pi_{\mathcal{O}_\mu} \cdot w \neq 0, $ then we have that
\begin{align*}
& \omega_{\mathcal{K}_{\mathcal{O}_\mu}}(T\bar{\gamma}_{\mathcal{O}_\mu} \cdot X_H^\gamma, \;
\tau_{\mathcal{K}_{\mathcal{O}_\mu}}\cdot T\pi_{\mathcal{O}_\mu} \cdot w)=
\omega_{\mathcal{K}_{\mathcal{O}_\mu}}(\tau_{\mathcal{K}_{\mathcal{O}_\mu}}\cdot
T\bar{\gamma}_{\mathcal{O}_\mu} \cdot X_H^\gamma, \; \tau_{\mathcal{K}_{\mathcal{O}_\mu}}\cdot T\pi_{\mathcal{O}_\mu} \cdot w)\\
& = \tau_{\mathcal{K}_{\mathcal{O}_\mu}}\cdot \omega_{\mathcal{M}_{\mathcal{O}_\mu}}(T(\pi_{\mathcal{O}_\mu}
\cdot\gamma )\cdot X_H^\gamma, \; T\pi_{\mathcal{O}_\mu} \cdot w ) =
\tau_{\mathcal{K}_{\mathcal{O}_\mu}}\cdot i_{\mathcal{M}_{\mathcal{O}_\mu}}^* \cdot
\omega_{\mathcal{O}_\mu}(T\pi_{\mathcal{O}_\mu} \cdot T\gamma \cdot X_H^\gamma, \; T\pi_{\mathcal{O}_\mu}
\cdot w )\\
& = \tau_{\mathcal{K}_{\mathcal{O}_\mu}}\cdot i_{\mathcal{M}_{\mathcal{O}_\mu}}^*
\cdot \pi_{\mathcal{O}_\mu}^*\omega_{\mathcal{O}_\mu}(T\gamma \cdot T\pi_Q \cdot X_H \cdot
\gamma, \; w) = \tau_{\mathcal{K}_{\mathcal{O}_\mu}}\cdot i_{\mathcal{M}_{\mathcal{O}_\mu}}^*
\cdot \omega(T(\gamma \cdot \pi_Q) \cdot X_H \cdot \gamma, \; w)\\
& =\tau_{\mathcal{K}_{\mathcal{O}_\mu}}\cdot i_{\mathcal{M}_{\mathcal{O}_\mu}}^* \cdot
(\omega(X_H \cdot \gamma, \; w-T(\gamma
\cdot \pi_Q)\cdot w) -\mathbf{d}\gamma(T\pi_{Q}(X_H\cdot \gamma), \; T\pi_{Q}(w)))\\
& = \tau_{\mathcal{K}_{\mathcal{O}_\mu}}\cdot i_{\mathcal{M}_{\mathcal{O}_\mu}}^*
\cdot\pi_{\mathcal{O}_\mu}^*\omega_{\mathcal{O}_\mu}(X_H \cdot \gamma, \; w) -
\tau_{\mathcal{K}_{\mathcal{O}_\mu}}\cdot i_{\mathcal{M}_{\mathcal{O}_\mu}}^*
\cdot\pi_{\mathcal{O}_\mu}^*\omega_{\mathcal{O}_\mu}(X_H \cdot
\gamma, \; T(\gamma
\cdot \pi_Q) \cdot w)\\
& \;\;\;\;\;\; -\tau_{\mathcal{K}_{\mathcal{O}_\mu}}\cdot i_{\mathcal{M}_{\mathcal{O}_\mu}}^* \cdot\mathbf{d}\gamma(T\pi_{Q}(X_H\cdot \gamma), \; T\pi_{Q}(w))\\
& = \tau_{\mathcal{K}_{\mathcal{O}_\mu}}\cdot i_{\mathcal{M}_{\mathcal{O}_\mu}}^*
\cdot\omega_{\mathcal{O}_\mu}(T\pi_{\mathcal{O}_\mu}(X_H\cdot \gamma), \; T\pi_{\mathcal{O}_\mu} \cdot w)\\
& \;\;\;\;\;\; -\tau_{\mathcal{K}_{\mathcal{O}_\mu}}\cdot i_{\mathcal{M}_{\mathcal{O}_\mu}}^* \cdot
\omega_{\mathcal{O}_\mu}(T\pi_{\mathcal{O}_\mu}\cdot(X_H\cdot \gamma), \; T(\pi_{\mathcal{O}_\mu}
\cdot\gamma) \cdot T\pi_{Q}(w))\\
& \;\;\;\;\;\; -\tau_{\mathcal{K}_{\mathcal{O}_\mu}}\cdot
i_{\mathcal{M}_{\mathcal{O}_\mu}}^* \cdot\mathbf{d}\gamma(T\pi_{Q}(X_H\cdot
\gamma), \; T\pi_{Q}(w))\\
& = \tau_{\mathcal{K}_{\mathcal{O}_\mu}}\cdot i_{\mathcal{M}_{\mathcal{O}_\mu}}^* \cdot
\omega_{\mathcal{O}_\mu}(X_{h_{\mathcal{K}_{\mathcal{O}_\mu}}} \cdot \bar{\gamma}_{\mathcal{O}_\mu}, \; T\pi_{\mathcal{O}_\mu} \cdot w)-
\tau_{\mathcal{K}_{\mathcal{O}_\mu}}\cdot i_{\mathcal{M}_{\mathcal{O}_\mu}}^* \cdot
\omega_{\mathcal{O}_\mu}(X_{h_{\mathcal{K}_{\mathcal{O}_\mu}}} \cdot \bar{\gamma}_{\mathcal{O}_\mu},
\; T\bar{\gamma}_{\mathcal{O}_\mu} \cdot T\pi_{Q}(w))\\
& \;\;\;\;\;\; -\tau_{\mathcal{K}_{\mathcal{O}_\mu}}\cdot i_{\mathcal{M}_{\mathcal{O}_\mu}}^* \cdot\mathbf{d}\gamma(T\pi_{Q}(X_H\cdot \gamma), \; T\pi_{Q}(w))\\
& = \omega_{\mathcal{K}_{\mathcal{O}_\mu}}( \tau_{\mathcal{K}_{\mathcal{O}_\mu}}\cdot X_{h_{\mathcal{K}_{\mathcal{O}_\mu}}}
\cdot \bar{\gamma}_{\mathcal{O}_\mu}, \; \tau_{\mathcal{K}_{\mathcal{O}_\mu}}\cdot T\pi_{\mathcal{O}_\mu}
\cdot w) - \omega_{\mathcal{K}_{\mathcal{O}_\mu}}(\tau_{\mathcal{K}_{\mathcal{O}_\mu}}\cdot
X_{h_{\mathcal{K}_{\mathcal{O}_\mu}}} \cdot \bar{\gamma}_{\mathcal{O}_\mu}, \;
\tau_{\mathcal{K}_{\mathcal{O}_\mu}}\cdot T\bar{\gamma}_{\mathcal{O}_\mu} \cdot T\pi_{Q}(w))\\
& \;\;\;\;\;\; -\tau_{\mathcal{K}_{\mathcal{O}_\mu}}\cdot i_{\mathcal{M}_{\mathcal{O}_\mu}}^* \cdot\mathbf{d}\gamma(T\pi_{Q}(X_H\cdot \gamma), \; T\pi_{Q}(w))\\
& = \omega_{\mathcal{K}_{\mathcal{O}_\mu}}(X_{\mathcal{K}_{\mathcal{O}_\mu}}\cdot
\bar{\gamma}_{\mathcal{O}_\mu}, \; \tau_{\mathcal{K}_{\mathcal{O}_\mu}} \cdot T\pi_{\mathcal{O}_\mu} \cdot w)
- \omega_{\mathcal{K}_{\mathcal{O}_\mu}}(X_{\mathcal{K}_{\mathcal{O}_\mu}} \cdot
\bar{\gamma}_{\mathcal{O}_\mu}, \; T\bar{\gamma}_{\mathcal{O}_\mu}
\cdot T\pi_{Q}(w))\\
& \;\;\;\;\;\; -\tau_{\mathcal{K}_{\mathcal{O}_\mu}}\cdot i_{\mathcal{M}_{\mathcal{O}_\mu}}^*
\cdot\mathbf{d}\gamma(T\pi_{Q}(X_H\cdot \gamma), \; T\pi_{Q}(w)),
\end{align*}
where we have used that $ \tau_{\mathcal{K}_{\mathcal{O}_\mu}}\cdot T\bar{\gamma}_{\mathcal{O}_\mu}=
T\bar{\gamma}_{\mathcal{O}_\mu}, $ and $\tau_{\mathcal{K}_{\mathcal{O}_\mu}}\cdot X_{h_{\mathcal{K}_{\mathcal{O}_\mu}}}\cdot
\bar{\gamma}_{\mathcal{O}_\mu}
= X_{\mathcal{K}_{\mathcal{O}_\mu}}\cdot \bar{\gamma}_{\mathcal{O}_\mu}, $
since $\textmd{Im}(T\bar{\gamma}_{\mathcal{O}_\mu})\subset \mathcal{K}_{\mathcal{O}_\mu}. $
If the one-form $\gamma: Q \rightarrow T^*Q $ is closed on $\mathcal{D}$ with respect to
$T\pi_Q: TT^* Q \rightarrow TQ, $ then we have that
$\mathbf{d}\gamma(T\pi_{Q}(X_H\cdot \gamma), \; T\pi_{Q}(w))=0, $
since $X_{H}\cdot \gamma, \; w \in \mathcal{F},$ and
$T\pi_{Q}(X_H\cdot \gamma), \; T\pi_{Q}(w) \in \mathcal{D}, $ and hence
$$
\tau_{\mathcal{K}_{\mathcal{O}_\mu}}\cdot
i_{\mathcal{M}_{\mathcal{O}_\mu}}^* \cdot\mathbf{d}\gamma(T\pi_{Q}(X_H\cdot \gamma),
\; T\pi_{Q}(w))=0,
$$
and
\begin{align}
& \omega_{\mathcal{K}_{\mathcal{O}_\mu}}(T\bar{\gamma}_{\mathcal{O}_\mu} \cdot X_H^\gamma, \;
\tau_{\mathcal{K}_{\mathcal{O}_\mu}}\cdot T\pi_{\mathcal{O}_\mu} \cdot w)
- \omega_{\mathcal{K}_{\mathcal{O}_\mu}}(X_{\mathcal{K}_{\mathcal{O}_\mu}}\cdot
\bar{\gamma}_{\mathcal{O}_\mu}, \; \tau_{\mathcal{K}_{\mathcal{O}_\mu}} \cdot T\pi_{\mathcal{O}_\mu} \cdot w) \nonumber \\
& = -\omega_{\mathcal{K}_{\mathcal{O}_\mu}}(X_{\mathcal{K}_{\mathcal{O}_\mu}} \cdot
\bar{\gamma}_{\mathcal{O}_\mu}, \; T\bar{\gamma}_{\mathcal{O}_\mu}
\cdot T\pi_{Q}(w)).
\end{align}
If $\bar{\gamma}_{\mathcal{O}_\mu}$ satisfies the equation
$T\bar{\gamma}_{\mathcal{O}_\mu}\cdot X_H^\gamma= X_{\mathcal{K}_{\mathcal{O}_\mu}}\cdot
\bar{\gamma}_{\mathcal{O}_\mu},$
from Lemma 3.4(i) we know that the right side of (5.5) becomes
\begin{align*}
& -\omega_{\mathcal{K}_{\mathcal{O}_\mu}}(X_{\mathcal{K}_{\mathcal{O}_\mu}} \cdot
\bar{\gamma}_{\mathcal{O}_\mu}, \; \tau_{\mathcal{K}_{\mathcal{O}_\mu}}\cdot T\bar{\gamma}_{\mathcal{O}_\mu}
\cdot T\pi_{Q}(w))\\
& = -\omega_{\mathcal{K}_{\mathcal{O}_\mu}}(T\bar{\gamma}_{\mathcal{O}_\mu}\cdot X_H^\gamma,
\; T\bar{\gamma}_{\mathcal{O}_\mu} \cdot T\pi_{Q}(w))\\
& = - \omega_{\mathcal{K}_{\mathcal{O}_\mu}}
(\tau_{\mathcal{K}_{\mathcal{O}_\mu}}T\bar{\gamma}_{\mathcal{O}_\mu}\cdot X_H^\gamma,
\; \tau_{\mathcal{K}_{\mathcal{O}_\mu}}\cdot T\bar{\gamma}_{\mathcal{O}_\mu} \cdot T\pi_{Q}(w))\\
& = -\tau_{\mathcal{K}_{\mathcal{O}_\mu}}\cdot i_{\mathcal{M}_{\mathcal{O}_\mu}}^* \cdot
\omega_{\mathcal{O}_\mu}(T\bar{\gamma}_{\mathcal{O}_\mu}\cdot X_H^\gamma,
\; T\bar{\gamma}_{\mathcal{O}_\mu} \cdot T\pi_{Q}(w))\\
& = -\tau_{\mathcal{K}_{\mathcal{O}_\mu}}\cdot i_{\mathcal{M}_{\mathcal{O}_\mu}}^* \cdot \bar{\gamma}_{\mathcal{O}_\mu}^* \cdot
\omega_{\mathcal{O}_\mu}( T\pi_{Q} \cdot X_{H} \cdot
\gamma, \; T\pi_{Q}(w))\\
& = -\tau_{\mathcal{K}_{\mathcal{O}_\mu}}\cdot i_{\mathcal{M}_{\mathcal{O}_\mu}}^* \cdot
 \gamma^* \cdot \pi^*_{\mathcal{O}_\mu}\cdot \omega_{\mathcal{O}_\mu}(T\pi_{Q} \cdot X_{H} \cdot
\gamma, \; T\pi_{Q}(w))\\
& = -\tau_{\mathcal{K}_{\mathcal{O}_\mu}}\cdot
i_{\mathcal{M}_{\mathcal{O}_\mu}}^* \cdot\gamma^*\omega( T\pi_{Q}(X_{H}\cdot\gamma), \; T\pi_{Q}(w))\\
& = \tau_{\mathcal{K}_{\mathcal{O}_\mu}}\cdot i_{\mathcal{M}_{\mathcal{O}_\mu}}^* \cdot
\mathbf{d}\gamma(T\pi_{Q}( X_{H}\cdot\gamma ), \; T\pi_{Q}(w))=0.
\end{align*}
But, because the $\mathbf{J}$-nonholonomic $R_o$-reduced distributional two-form
$\omega_{\mathcal{K}_{\mathcal{O}_\mu}}$ is non-degenerate,
the left side of (5.5) equals zero, only when
$\bar{\gamma}_{\mathcal{O}_\mu}$ satisfies the equation
$T\bar{\gamma}_{\mathcal{O}_\mu}\cdot X_H^\gamma= X_{\mathcal{K}_{\mathcal{O}_\mu}}\cdot
\bar{\gamma}_{\mathcal{O}_\mu}.$ Thus,
if the one-form $\gamma: Q \rightarrow T^*Q $ is closed on $\mathcal{D}$ with respect to
$T\pi_Q: TT^* Q \rightarrow TQ, $ then $\bar{\gamma}_{\mathcal{O}_\mu}$ must be a solution
of the Type I of Hamilton-Jacobi equation
$T\bar{\gamma}_{\mathcal{O}_\mu}\cdot X_H^\gamma= X_{\mathcal{K}_{\mathcal{O}_\mu}}\cdot
\bar{\gamma}_{\mathcal{O}_\mu}. $
\hskip 0.3cm $\blacksquare$\\

Next, for any $G$-invariant symplectic map $\varepsilon: T^* Q \rightarrow T^* Q $,
we can prove the following Type II of geometric
Hamilton-Jacobi theorem for the $\mathbf{J}$-nonholonomic $R_o$-reduced distributional
Hamiltonian system.
For convenience, the maps involved in the following
theorem and its proof are shown in Diagram-8.
\begin{center}
\hskip 0cm \xymatrix{ \mathbf{J}^{-1}(\mathcal{O}_\mu) \ar[r]^{i_{\mathcal{O}_\mu}} & T^* Q
\ar[d]_{X_{H\cdot \varepsilon}} \ar[dr]^{X_H^\varepsilon} \ar[r]^{\pi_Q}
& Q \ar[r]^{\gamma} & T^*Q \ar[d]_{X_H} \ar[dr]^{X_{h_{\mathcal{O}_\mu}}
\cdot\bar{\varepsilon}} \ar[r]^{\pi_{\mathcal{O}_\mu}}
& (T^* Q)_{\mathcal{O}_\mu} \ar[d]^{X_{h_{\mathcal{O}_\mu}}}
& \mathcal{M}_{\mathcal{O}_\mu} \ar[l]_{i_{\mathcal{M}_{\mathcal{O}_\mu}}}
\ar[d]_{X_{\mathcal{K}_{\mathcal{O}_\mu}}}\\
& T(T^*Q)  & TQ \ar[l]^{T\gamma} & T(T^*Q) \ar[l]^{T\pi_Q} \ar[r]_{T\pi_{\mathcal{O}_\mu}}
& T(T^* Q)_{\mathcal{O}_\mu} \ar[r]^{\tau_{\mathcal{K}_{\mathcal{O}_\mu}}} & \mathcal{K}_{\mathcal{O}_\mu} }
\end{center}
$$\mbox{Diagram-8}$$

\begin{theo} (Type II of Hamilton-Jacobi Theorem for a $\mathbf{J}$-Nonholonomic
$R_o$-reduced Distributional Hamiltonian System) For a given
$\mathbf{J}$-nonholonomic regular orbit reducible Hamiltonian system
$(T^*Q,G,\omega,\mathbf{J},\mathcal{D},H)$ with an associated
$\mathbf{J}$-nonholonomic $R_o$-reduced distributional Hamiltonian system
$(\mathcal{K}_{\mathcal{O}_\mu},\omega_{\mathcal{K}_{\mathcal{O}_\mu}},
h_{\mathcal{K}_{\mathcal{O}_\mu}})$, assume that $\gamma:
Q \rightarrow T^*Q$ is an one-form on $Q$, and $\lambda=\gamma \cdot
\pi_{Q}: T^* Q \rightarrow T^* Q, $ and for any
symplectic map $\varepsilon:T^* Q \rightarrow T^* Q, $
denote $X_H^\varepsilon =
T\pi_{Q}\cdot X_H \cdot \varepsilon$, where $X_{H}$ is the Hamiltonian
vector field of the corresponding unconstrained Hamiltonian system with
symmetry and momentum map $(T^*Q,G,\omega,\mathbf{J},H)$. Moreover,
assume that $\mu\in\mathfrak{g}^\ast$ is a regular value of the momentum
map $\mathbf{J}$, and $\textmd{Im}(\gamma)\subset \mathcal{M} \cap
\mathbf{J}^{-1}(\mathcal{O}_\mu), $ and it is $G$-invariant,
and $\varepsilon$ is also $G$-invariant and
$\varepsilon(\mathbf{J}^{-1}(\mathcal{O}_\mu)) \subset \mathbf{J}^{-1}(\mathcal{O}_\mu). $ Denote
$\bar{\gamma}_{\mathcal{O}_\mu}=\pi_{\mathcal{O}_\mu}(\gamma): Q \rightarrow \mathcal{M}_{\mathcal{O}_\mu} $,
and $ \textmd{Im}(T\bar{\gamma}_{\mathcal{O}_\mu})\subset \mathcal{K}_{\mathcal{O}_\mu}, $ and
$\bar{\lambda}_{\mathcal{O}_\mu}=\pi_{\mathcal{O}_\mu}(\lambda): \mathbf{J}^{-1}(\mathcal{O}_\mu)
(\subset T^* Q) \rightarrow \mathcal{M}_{\mathcal{O}_\mu}, $ and
$\bar{\varepsilon}_{\mathcal{O}_\mu}=\pi_{\mathcal{O}_\mu}(\varepsilon):
\mathbf{J}^{-1}(\mathcal{O}_\mu) (\subset T^* Q) \rightarrow
\mathcal{M}_{\mathcal{O}_\mu}. $ Then $\varepsilon$ and $\bar{\varepsilon}_{\mathcal{O}_\mu}$ satisfy the
equation $\tau_{\mathcal{K}_{\mathcal{O}_\mu}} \cdot T\bar{\varepsilon}(X_{h_{\mathcal{K}_{\mathcal{O}_\mu}}\cdot \bar{\varepsilon}_{\mathcal{O}_\mu}})
= T\bar{\lambda}_{\mathcal{O}_\mu} \cdot X_H \cdot\varepsilon, $ if and only if
they satisfy the equation $T\bar{\gamma}_{\mathcal{O}_\mu}\cdot
X_H^\varepsilon= X_{\mathcal{K}_{\mathcal{O}_\mu}}\cdot \bar{\varepsilon}_{\mathcal{O}_\mu}. $
Here $X_{h_{\mathcal{K}_{\mathcal{O}_\mu}} \cdot\bar{\varepsilon}_{\mathcal{O}_\mu}}$
is the Hamiltonian vector field of the
function $h_{\mathcal{K}_{\mathcal{O}_\mu}}\cdot \bar{\varepsilon}_{\mathcal{O}_\mu}: T^* Q\rightarrow \mathbb{R}, $
and $X_{\mathcal{K}_{\mathcal{O}_\mu}}$ is the dynamical vector field
of the $\mathbf{J}$-nonholonomic $R_o$-reduced system $(\mathcal{K}_{\mathcal{O}_\mu},\omega_{\mathcal{K}_{\mathcal{O}_\mu}},h_{\mathcal{K}_{\mathcal{O}_\mu}})$.
The equation
$T\bar{\gamma}_{\mathcal{O}_\mu}\cdot X_H^\varepsilon= X_{\mathcal{K}_{\mathcal{O}_\mu}}\cdot
\bar{\varepsilon}_{\mathcal{O}_\mu},$ is called the Type II of Hamilton-Jacobi equation for the
$\mathbf{J}$-nonholonomic $R_o$-reduced distributional Hamiltonian system
$(\mathcal{K}_{\mathcal{O}_\mu},\omega_{\mathcal{K}_{\mathcal{O}_\mu}},h_{\mathcal{K}_{\mathcal{O}_\mu}})$.
\end{theo}

\noindent{\bf Proof: } At first, we note that
the $R_o$-reduced symplectic space
$(T^\ast Q)_{\mathcal{O}_\mu}= \mathbf{J}^{-1}(\mathcal{O}_\mu)/G
\cong (\mathbf{J}^{-1}(\mu)/G) \times \mathcal{O}_\mu, $ with the $R_o$-reduced
symplectic form $\omega_{\mathcal{O}_\mu}$ uniquely characterized by
the relation $i_{\mathcal{O}_\mu}^\ast
\omega=\pi_{\mathcal{O}_{\mu}}^\ast \omega_{\mathcal{O}
_\mu}+\mathbf{J}_{\mathcal{O}_\mu}^\ast\omega_{\mathcal{O}_\mu}^+. $
Since
$\textmd{Im}(\gamma)\subset \mathcal{M} \cap \mathbf{J}^{-1}(\mathcal{O}_\mu), $
and it is $G$-invariant,
in this case for any $V\in TQ, $ and $w\in TT^*Q, $
we have that
$\mathbf{J}_{\mathcal{O}_\mu}^\ast\omega_{\mathcal{O}_\mu}^+(T\gamma
\cdot V, \; w)=0, $ and hence
$\pi_{\mathcal{O}_\mu}^*\omega_{\mathcal{O}_\mu}=
i_{\mathcal{O}_\mu}^*\omega= \omega, $ along $\textmd{Im}(\gamma)$.
On the other hand, because
$\textmd{Im}(T\bar{\gamma}_{\mathcal{O}_\mu})\subset \mathcal{K}_{\mathcal{O}_\mu}, $
then $\omega_{\mathcal{K}_{\mathcal{O}_\mu}}\cdot
\tau_{\mathcal{K}_{\mathcal{O}_\mu}}=\tau_{\mathcal{K}_{\mathcal{O}_\mu}}\cdot
\omega_{\mathcal{M}_{\mathcal{O}_\mu}}= \tau_{\mathcal{K}_{\mathcal{O}_\mu}}\cdot
i_{\mathcal{M}_{\mathcal{O}_\mu}}^* \cdot \omega_{\mathcal{O}_\mu}, $ along
$\textmd{Im}(T\bar{\gamma}_{\mathcal{O}_\mu})$. Thus, using
the $\mathbf{J}$-nonholonomic $R_o$-reduced distributional two-form
$\omega_{\mathcal{K}_{\mathcal{O}_\mu}}$, from Lemma 3.4, if we take that $v=
\tau_{\mathcal{K}_{\mathcal{O}_\mu}}\cdot T\pi_{\mathcal{O}_\mu} \cdot X_{H}\cdot \varepsilon= X_{\mathcal{K}_{\mathcal{O}_\mu}}\cdot
\bar{\varepsilon}_{\mathcal{O}_\mu} \in \mathcal{K}_{\mathcal{O}_\mu}, $ and for any $w \in
\mathcal{F}, \; T\lambda(w)\neq 0, $ and
$\tau_{\mathcal{K}_{\mathcal{O}_\mu}}\cdot T\pi_{\mathcal{O}_\mu} \cdot w \neq 0, $ then we have that
\begin{align*}
& \omega_{\mathcal{K}_{\mathcal{O}_\mu}}(T\bar{\gamma}_{\mathcal{O}_\mu} \cdot X_H^\varepsilon, \;
\tau_{\mathcal{K}_{\mathcal{O}_\mu}}\cdot T\pi_{\mathcal{O}_\mu} \cdot w)=
\omega_{\mathcal{K}_{\mathcal{O}_\mu}}(\tau_{\mathcal{K}_{\mathcal{O}_\mu}}\cdot
T\bar{\gamma}_{\mathcal{O}_\mu} \cdot X_H^\varepsilon, \;
\tau_{\mathcal{K}_{\mathcal{O}_\mu}}\cdot T\pi_{\mathcal{O}_\mu} \cdot w)\\
& = \tau_{\mathcal{K}_{\mathcal{O}_\mu}}\cdot \omega_{\mathcal{M}_{\mathcal{O}_\mu}}(T(\pi_{\mathcal{O}_\mu}
\cdot\gamma )\cdot X_H^\varepsilon, \; T\pi_{\mathcal{O}_\mu} \cdot w ) =
\tau_{\mathcal{K}_{\mathcal{O}_\mu}}\cdot i_{\mathcal{M}_{\mathcal{O}_\mu}}^* \cdot
\omega_{\mathcal{O}_\mu}(T\pi_{\mathcal{O}_\mu} \cdot T\gamma \cdot X_H^\varepsilon,
\; T\pi_{\mathcal{O}_\mu} \cdot w )\\
& = \tau_{\mathcal{K}_{\mathcal{O}_\mu}}\cdot i_{\mathcal{M}_{\mathcal{O}_\mu}}^*
\cdot \pi_{\mathcal{O}_\mu}^*\omega_{\mathcal{O}_\mu}(T\gamma \cdot T\pi_Q \cdot X_H \cdot
\varepsilon, \; w) = \tau_{\mathcal{K}_{\mathcal{O}_\mu}}\cdot i_{\mathcal{M}_{\mathcal{O}_\mu}}^*
\cdot \omega(T(\gamma \cdot \pi_Q) \cdot X_H \cdot \varepsilon, \; w)\\
& =\tau_{\mathcal{K}_{\mathcal{O}_\mu}}\cdot i_{\mathcal{M}_{\mathcal{O}_\mu}}^* \cdot
(\omega(X_H \cdot \varepsilon, \; w-T(\gamma
\cdot \pi_Q)\cdot w) -\mathbf{d}\gamma(T\pi_{Q}(X_H\cdot \varepsilon), \; T\pi_{Q}(w)))\\
& = \tau_{\mathcal{K}_{\mathcal{O}_\mu}}\cdot i_{\mathcal{M}_{\mathcal{O}_\mu}}^* \cdot
\omega(X_H \cdot \varepsilon, \; w) - \tau_{\mathcal{K}_{\mathcal{O}_\mu}}\cdot
i_{\mathcal{M}_{\mathcal{O}_\mu}}^* \cdot \omega(X_H \cdot
\varepsilon, \; T\lambda \cdot w)\\
& \;\;\;\;\;\; -\tau_{\mathcal{K}_{\mathcal{O}_\mu}}\cdot i_{\mathcal{M}_{\mathcal{O}_\mu}}^* \cdot\mathbf{d}\gamma(T\pi_{Q}(X_H\cdot \varepsilon), \; T\pi_{Q}(w))\\
& = \tau_{\mathcal{K}_{\mathcal{O}_\mu}}\cdot i_{\mathcal{M}_{\mathcal{O}_\mu}}^*
\cdot\pi_{\mathcal{O}_\mu}^*\omega_{\mathcal{O}_\mu}(X_H \cdot \varepsilon, \; w) -
\tau_{\mathcal{K}_{\mathcal{O}_\mu}}\cdot i_{\mathcal{M}_{\mathcal{O}_\mu}}^*
\cdot\pi_{\mathcal{O}_\mu}^*\omega_{\mathcal{O}_\mu}(X_H \cdot
\varepsilon, \; T\lambda \cdot w)\\
& \;\;\;\;\;\; +\tau_{\mathcal{K}_{\mathcal{O}_\mu}}\cdot i_{\mathcal{M}_{\mathcal{O}_\mu}}^*
\cdot \lambda^* \omega(X_H\cdot \varepsilon, \; w)\\
& = \tau_{\mathcal{K}_{\mathcal{O}_\mu}}\cdot i_{\mathcal{M}_{\mathcal{O}_\mu}}^*
\cdot\omega_{\mathcal{O}_\mu}(T\pi_{\mathcal{O}_\mu}(X_H\cdot \varepsilon), \; T\pi_{\mathcal{O}_\mu} \cdot w)\\
& \;\;\;\;\;\; -\tau_{\mathcal{K}_{\mathcal{O}_\mu}}\cdot i_{\mathcal{M}_{\mathcal{O}_\mu}}^* \cdot
\omega_{\mathcal{O}_\mu}(T\pi_{\mathcal{O}_\mu}\cdot(X_H\cdot \varepsilon), \; T(\pi_{\mathcal{O}_\mu}
\cdot\lambda) \cdot w)\\
& \;\;\;\;\;\; +\tau_{\mathcal{K}_{\mathcal{O}_\mu}}\cdot
i_{\mathcal{M}_{\mathcal{O}_\mu}}^* \cdot
\pi_{\mathcal{O}_\mu}^*\omega_{\mathcal{O}_\mu}(T\lambda\cdot X_H \cdot
\varepsilon, \; T\lambda \cdot w)\\
& = \tau_{\mathcal{K}_{\mathcal{O}_\mu}}\cdot
i_{\mathcal{M}_{\mathcal{O}_\mu}}^* \cdot\omega_{\mathcal{O}_\mu}(T\pi_{\mathcal{O}_\mu}(X_H)\cdot
\pi_{\mathcal{O}_\mu}(\varepsilon), \; T\pi_{\mathcal{O}_\mu} \cdot w)\\
& \;\;\;\;\;\; -\tau_{\mathcal{K}_{\mathcal{O}_\mu}}\cdot
i_{\mathcal{M}_{\mathcal{O}_\mu}}^* \cdot \omega_{\mathcal{O}_\mu}(T\pi_{\mathcal{O}_\mu}(X_H)\cdot
\pi_{\mathcal{O}_\mu}(\varepsilon), \; T\bar{\lambda}_{\mathcal{O}_\mu} \cdot w)\\
& \;\;\;\;\;\; +\tau_{\mathcal{K}_{\mathcal{O}_\mu}}\cdot
i_{\mathcal{M}_{\mathcal{O}_\mu}}^* \cdot
\omega_{\mathcal{O}_\mu}(T\pi_{\mathcal{O}_\mu}\cdot T\lambda\cdot X_H \cdot
\varepsilon, \; T\pi_{\mathcal{O}_\mu}\cdot T\lambda \cdot w)\\
& = \tau_{\mathcal{K}_{\mathcal{O}_\mu}}\cdot i_{\mathcal{M}_{\mathcal{O}_\mu}}^* \cdot
\omega_{\mathcal{O}_\mu}(X_{h_{\mathcal{K}_{\mathcal{O}_\mu}}} \cdot \bar{\varepsilon}_{\mathcal{O}_\mu}, \; T\pi_{\mathcal{O}_\mu} \cdot w)-
\tau_{\mathcal{K}_{\mathcal{O}_\mu}}\cdot i_{\mathcal{M}_{\mathcal{O}_\mu}}^* \cdot
\omega_{\mathcal{O}_\mu}(X_{h_{\mathcal{K}_{\mathcal{O}_\mu}}} \cdot \bar{\varepsilon}_{\mathcal{O}_\mu},
\; T\bar{\lambda}_{\mathcal{O}_\mu} \cdot w)\\
& \;\;\;\;\;\; +\tau_{\mathcal{K}_{\mathcal{O}_\mu}}\cdot
i_{\mathcal{M}_{\mathcal{O}_\mu}}^* \cdot
\omega_{\mathcal{O}_\mu}(T\bar{\lambda}_{\mathcal{O}_\mu} \cdot X_H \cdot
\varepsilon, \; T\bar{\lambda}_{\mathcal{O}_\mu} \cdot w)\\
& = \omega_{\mathcal{K}_{\mathcal{O}_\mu}}( \tau_{\mathcal{K}_{\mathcal{O}_\mu}}\cdot X_{h_{\mathcal{K}_{\mathcal{O}_\mu}}}
\cdot \bar{\varepsilon}_{\mathcal{O}_\mu}, \; \tau_{\mathcal{K}_{\mathcal{O}_\mu}}\cdot T\pi_{\mathcal{O}_\mu}
\cdot w) - \omega_{\mathcal{K}_{\mathcal{O}_\mu}}(\tau_{\mathcal{K}_{\mathcal{O}_\mu}}\cdot
X_{h_{\mathcal{K}_{\mathcal{O}_\mu}}} \cdot \bar{\varepsilon}_{\mathcal{O}_\mu}, \;
\tau_{\mathcal{K}_{\mathcal{O}_\mu}}\cdot T\bar{\lambda}_{\mathcal{O}_\mu} \cdot w)\\
& \;\;\;\;\;\; +\omega_{\mathcal{K}_{\mathcal{O}_\mu}}(\tau_{\mathcal{K}_{\mathcal{O}_\mu}}
\cdot T\bar{\lambda}_{\mathcal{O}_\mu} \cdot X_H \cdot
\varepsilon, \; \tau_{\mathcal{K}_{\mathcal{O}_\mu}}\cdot T\bar{\lambda}_{\mathcal{O}_\mu} \cdot w)\\
& = \omega_{\mathcal{K}_{\mathcal{O}_\mu}}(X_{\mathcal{K}_{\mathcal{O}_\mu}}\cdot
\bar{\varepsilon}_{\mathcal{O}_\mu}, \; \tau_{\mathcal{K}_{\mathcal{O}_\mu}} \cdot T\pi_{\mathcal{O}_\mu} \cdot w)
- \omega_{\mathcal{K}_{\mathcal{O}_\mu}}(\tau_{\mathcal{K}_{\mathcal{O}_\mu}}\cdot X_{h_{\mathcal{K}_{\mathcal{O}_\mu}}} \cdot
\bar{\varepsilon}_{\mathcal{O}_\mu}, \; T\bar{\lambda}_{\mathcal{O}_\mu} \cdot w)\\
& \;\;\;\;\;\; +\omega_{\mathcal{K}_{\mathcal{O}_\mu}}(T\bar{\lambda}_{\mathcal{O}_\mu} \cdot X_H \cdot
\varepsilon, \; T\bar{\lambda}_{\mathcal{O}_\mu} \cdot w),
\end{align*}
where we have used that $ \tau_{\mathcal{K}_{\mathcal{O}_\mu}}\cdot T\bar{\gamma}_{\mathcal{O}_\mu}=
T\bar{\gamma}_{\mathcal{O}_\mu}, $ and $\tau_{\mathcal{K}_{\mathcal{O}_\mu}}\cdot X_{h_{\mathcal{K}_{\mathcal{O}_\mu}}}\cdot
\bar{\varepsilon}_{\mathcal{O}_\mu} = X_{\mathcal{K}_{\mathcal{O}_\mu}}\cdot \bar{\varepsilon}_{\mathcal{O}_\mu}, $
since $\textmd{Im}(T\bar{\gamma}_{\mathcal{O}_\mu})\subset \mathcal{K}_{\mathcal{O}_\mu}. $ Note
that $\varepsilon: T^* Q \rightarrow T^* Q $ is symplectic, and
$\pi_{\mathcal{O}_\mu}^*\omega_{\mathcal{O}_\mu}= i_{\mathcal{O}_\mu}^*\omega= \omega, $ along
$\textmd{Im}(\gamma)$, and hence $\bar{\varepsilon}_{\mathcal{O}_\mu}=
\pi_{\mathcal{O}_\mu}(\varepsilon): T^* Q \rightarrow (T^* Q)_{\mathcal{O}_\mu} $ is also symplectic
along $\textmd{Im}(\gamma)$, and hence $X_{h_{\mathcal{K}_{\mathcal{O}_\mu}}}\cdot
\bar{\varepsilon}_{\mathcal{O}_\mu}= T\bar{\varepsilon}_{\mathcal{O}_\mu} \cdot X_{h_{\mathcal{K}_{\mathcal{O}_\mu}} \cdot
\bar{\varepsilon}_{\mathcal{O}_\mu}}, $ along $\bar{\varepsilon}_{\mathcal{O}_\mu}$, and hence
$\tau_{\mathcal{K}_{\mathcal{O}_\mu}}\cdot X_{h_{\mathcal{K}_{\mathcal{O}_\mu}}} \cdot
\bar{\varepsilon}_{\mathcal{O}_\mu}= \tau_{\mathcal{K}_{\mathcal{O}_\mu}}\cdot T\bar{\lambda}_{\mathcal{O}_\mu}
\cdot X_{h_{\mathcal{K}_{\mathcal{O}_\mu}}\cdot \bar{\varepsilon}_{\mathcal{O}_\mu}}, $
along $\bar{\varepsilon}_{\mathcal{O}_\mu}$, because
$\textmd{Im}(T\bar{\gamma}_{\mathcal{O}_\mu})\subset \mathcal{K}_{\mathcal{O}_\mu}. $ Then we
have that
\begin{align*}
& \omega_{\mathcal{K}_{\mathcal{O}_\mu}}(T\bar{\gamma}_{\mathcal{O}_\mu} \cdot X_H^\varepsilon, \;
\tau_{\mathcal{K}_{\mathcal{O}_\mu}}\cdot T\pi_{\mathcal{O}_\mu} \cdot w)-
\omega_{\mathcal{K}_{\mathcal{O}_\mu}}(X_{\mathcal{K}_{\mathcal{O}_\mu}}\cdot \bar{\varepsilon}_{\mathcal{O}_\mu},
\; \tau_{\mathcal{K}_{\mathcal{O}_\mu}} \cdot T\pi_{\mathcal{O}_\mu} \cdot w) \nonumber \\
& = - \omega_{\mathcal{K}_{\mathcal{O}_\mu}}(\tau_{\mathcal{K}_{\mathcal{O}_\mu}}\cdot X_{h_{\mathcal{K}_{\mathcal{O}_\mu}}} \cdot
\bar{\varepsilon}_{\mathcal{O}_\mu}, \; T\bar{\lambda}_{\mathcal{O}_\mu} \cdot w)
+\omega_{\mathcal{K}_{\mathcal{O}_\mu}}(T\bar{\lambda}_{\mathcal{O}_\mu} \cdot X_H \cdot
\varepsilon, \; T\bar{\lambda}_{\mathcal{O}_\mu} \cdot w)\\
& = \omega_{\mathcal{K}_{\mathcal{O}_\mu}}(T\bar{\lambda}_{\mathcal{O}_\mu} \cdot X_H \cdot
\varepsilon- \tau_{\mathcal{K}_{\mathcal{O}_\mu}}\cdot T\bar{\varepsilon}_{\mathcal{O}_\mu}
\cdot X_{h_{\mathcal{K}_{\mathcal{O}_\mu}} \cdot \bar{\varepsilon}_{\mathcal{O}_\mu}}, \; T\bar{\lambda}_{\mathcal{O}_\mu} \cdot w).
\end{align*}
Because the $\mathbf{J}$-nonholonomic $R_o$-reduced distributional two-form
$\omega_{\mathcal{K}_{\mathcal{O}_\mu}}$ is non-degenerate, it follows that the equation
$T\bar{\gamma}_{\mathcal{O}_\mu}\cdot X_H^\varepsilon = X_{\mathcal{K}_{\mathcal{O}_\mu}}\cdot
\bar{\varepsilon}_{\mathcal{O}_\mu}, $ is equivalent to the equation
$T\bar{\lambda}_{\mathcal{O}_\mu}\cdot X_H \cdot \varepsilon
= \tau_{\mathcal{K}_{\mathcal{O}_\mu}}\cdot T\bar{\varepsilon}_{\mathcal{O}_\mu} \cdot X_{h_{\mathcal{K}_{\mathcal{O}_\mu}}
\cdot \bar{\varepsilon}_{\mathcal{O}_\mu}}. $
Thus, $\varepsilon$ and $\bar{\varepsilon}_{\mathcal{O}_\mu}$ satisfy the equation
$T\bar{\lambda}_{\mathcal{O}_\mu}\cdot X_H \cdot \varepsilon
= \tau_{\mathcal{K}_{\mathcal{O}_\mu}}\cdot T\bar{\varepsilon}_{\mathcal{O}_\mu} \cdot X_{h_{\mathcal{K}_{\mathcal{O}_\mu}}
\cdot \bar{\varepsilon}_{\mathcal{O}_\mu}}, $ if and only if they satisfy
the Type II of Hamilton-Jacobi equation
$T\bar{\gamma}_{\mathcal{O}_\mu}\cdot X_H^\varepsilon = X_{\mathcal{K}_{\mathcal{O}_\mu}}\cdot
\bar{\varepsilon}_{\mathcal{O}_\mu}.$
\hskip 0.3cm $\blacksquare$\\

\begin{rema}
If a $\mathbf{J}$-nonholonomic regular orbit reducible Hamiltonian system
we considered has not any constrains,
in this case, the $\mathbf{J}$-nonholonomic $R_o$-reduced distributional
Hamiltonian system is just the regular orbit reduced Hamiltonian system itself.
From the above Type I and Type II of Hamilton-Jacobi theorems, that is,
Theorem 5.9 and Theorem 5.10, we can get the Theorem 4.3
and Theorem 4.4 in Wang \cite{wa17}.
It shows that Theorem 5.9 and Theorem 5.10 can be regarded as an extension of two types of
Hamilton-Jacobi theorem for the $R_o$-reduced Hamiltonian system
given in \cite{wa17} to the nonholonomic context.
\end{rema}

\begin{rema}
It is worthy of note that the formulations of Type I and Type II of Hamilton-Jacobi
equation for a $\mathbf{J}$-nonholonomic $R_o$-reduced distributional
Hamiltonian system, given by Theorem 5.9 and Theorem 5.10, have
more extensive sense, because, in general, the one-form $\gamma$ is not
given by a generating function of a symplectic map. When $\gamma$ is a solution
of the classical Hamilton-Jacobi equation, that is, $X_H\cdot
\gamma=0, $ then $X_H^\gamma= T\pi_Q\cdot X_H\cdot \gamma=0,$ and
hence from the Type I of Hamilton-Jacobi equation, we have that
$X_{\mathcal{K}_{\mathcal{O}_\mu}}\cdot \bar{\gamma}_{\mathcal{O}_\mu}
= T\bar{\gamma}_{\mathcal{O}_\mu} \cdot X_H^\gamma=0, $ which shows that the
dynamical vector field of the $\mathbf{J}$-nonholonomic $R_o$-reduced
distributional Hamiltonian system
$(\mathcal{K}_{\mathcal{O}_\mu},\omega_{\mathcal{K}_{\mathcal{O}_\mu}},h_{\mathcal{K}_{\mathcal{O}_\mu}})$
is degenerate along
$\bar{\gamma}_{\mathcal{O}_\mu}$. The equation $X_{\mathcal{K}_{\mathcal{O}_\mu}}\cdot
\bar{\gamma}_{\mathcal{O}_\mu}=0$ is called the classical Hamilton-Jacobi equation
for the $\mathbf{J}$-nonholonomic $R_o$-reduced distributional Hamiltonian system
$(\mathcal{K}_{\mathcal{O}_\mu},\omega_{\mathcal{K}_{\mathcal{O}_\mu}},h_{\mathcal{K}_{\mathcal{O}_\mu}})$.
In addition,
for a symplectic map $\varepsilon: T^* Q \rightarrow T^* Q $, if
$X_H\cdot \varepsilon=0, $ then from the Type II of Hamilton-Jacobi
equation, we have that $X_{\mathcal{K}_{\mathcal{O}_\mu}}\cdot
\bar{\varepsilon}_{\mathcal{O}_\mu}= T\bar{\gamma}_{\mathcal{O}_\mu}\cdot X_H^\varepsilon=0. $
But, from the equation $T\bar{\lambda}_{\mathcal{O}_\mu}\cdot X_H \cdot
\varepsilon = \tau_{\mathcal{K}_{\mathcal{O}_\mu}}\cdot T\bar{\varepsilon}_{\mathcal{O}_\mu}
\cdot X_{h_{\mathcal{K}_{\mathcal{O}_\mu}}\cdot \bar{\varepsilon}_{\mathcal{O}_\mu}}, $
we know that the equation
$X_{\mathcal{K}_{\mathcal{O}_\mu}}\cdot \bar{\varepsilon}_{\mathcal{O}_\mu}=0 $ is not
equivalent to the equation
$X_{h_{\mathcal{K}_{\mathcal{O}_\mu}}\cdot \bar{\varepsilon}_{\mathcal{O}_\mu}}=0.$
\end{rema}

For a given $\mathbf{J}$-nonholonomic regular orbit reducible Hamiltonian system
$(T^*Q,G,\omega,\mathbf{J},\mathcal{D},H)$ with an associated $\mathbf{J}$-
nonholonomic $R_o$-reduced distributional Hamiltonian system
$(\mathcal{K}_{\mathcal{O}_\mu},\omega_{\mathcal{K}_{\mathcal{O}_\mu}},
h_{\mathcal{K}_{\mathcal{O}_\mu}})$,
we know that the Hamiltonian vector field
$X_{H}$ and the $\mathbf{J}$-nonholonomic $R_o$-reduced dynamical vector field
$X_{h_{\mathcal{K}_{\mathcal{O}_\mu}}}$ are $\pi_{\mathcal{O}_\mu}$-related,
that is, $X_{h_{\mathcal{K}_{\mathcal{O}_\mu}}}\cdot \pi_{\mathcal{O}_\mu}
=T\pi_{\mathcal{O}_\mu}\cdot X_{H}\cdot i_{\mathcal{O}_\mu}.$ Then
we can prove the following Theorem 5.13, which states the
relationship between the solutions of Type II of Hamilton-Jacobi equations and
$\mathbf{J}$-nonholonomic regular orbit reduction.

\begin{theo}
For a given $\mathbf{J}$-nonholonomic regular orbit reducible Hamiltonian system
$(T^*Q, G,\\ \omega,\mathbf{J},\mathcal{D},H)$ with an associated
$\mathbf{J}$-nonholonomic $R_o$-reduced distributional Hamiltonian system \\
$(\mathcal{K}_{\mathcal{O}_\mu},\omega_{\mathcal{K}_{\mathcal{O}_\mu}},
h_{\mathcal{K}_{\mathcal{O}_\mu}})$, assume that $\gamma:
Q \rightarrow T^*Q$ is an one-form on $Q$, and $\lambda=\gamma \cdot
\pi_{Q}: T^* Q \rightarrow T^* Q, $ and
$\varepsilon: T^* Q \rightarrow T^* Q $ is a symplectic map. Moreover, assume
that $\mu \in \mathfrak{g}^\ast $ is a regular value of the momentum
map $\mathbf{J}$, and $\textmd{Im}(\gamma)\subset \mathcal{M} \cap
\mathbf{J}^{-1}(\mathcal{O}_\mu), $ and it is $G$-invariant,
and $\varepsilon$ is also $G$-invariant and
$\varepsilon(\mathbf{J}^{-1}(\mathcal{O}_\mu)) \subset \mathbf{J}^{-1}(\mathcal{O}_\mu). $
Denote $\bar{\gamma}_{\mathcal{O}_\mu}=\pi_{\mathcal{O}_\mu}(\gamma): Q \rightarrow \mathcal{M}_{\mathcal{O}_\mu} $,
and $ \textmd{Im}(T\bar{\gamma}_{\mathcal{O}_\mu})\subset \mathcal{K}_{\mathcal{O}_\mu}, $ and
$\bar{\lambda}_{\mathcal{O}_\mu}=\pi_{\mathcal{O}_\mu}(\lambda):
\mathbf{J}^{-1}(\mathcal{O}_\mu) (\subset T^*Q) \rightarrow (T^* Q)_{\mathcal{O}_\mu} $,
and $\bar{\varepsilon}_{\mathcal{O}_\mu}=\pi_{\mathcal{O}_\mu}(\varepsilon):
\mathbf{J}^{-1}(\mathcal{O}_\mu) (\subset T^*Q) \rightarrow (T^* Q)_{\mathcal{O}_\mu} $.
Then $\varepsilon$
is a solution of the Type II of Hamilton-Jacobi equation $T\gamma\cdot
X_H^\varepsilon= X_{\mathcal{K}}\cdot \varepsilon, $ for the distributional
Hamiltonian system $(\mathcal{K},\omega_{\mathcal {K}},H_{\mathcal {K}})$, if and
only if $\varepsilon$ and $\bar{\varepsilon}_{\mathcal{O}_\mu} $ satisfy the Type II of Hamilton-Jacobi
equation $T\bar{\gamma}_{\mathcal{O}_\mu}\cdot X_H^\varepsilon=
X_{\mathcal{K}_{\mathcal{O}_\mu}}\cdot \bar{\varepsilon}_{\mathcal{O}_\mu}, $ for the
$\mathbf{J}$-nonholonomic $R_o$-reduced distributional Hamiltonian system
$(\mathcal{K}_{\mathcal{O}_\mu},\omega_{\mathcal{K}_{\mathcal{O}_\mu}},h_{\mathcal{K}_{\mathcal{O}_\mu}})$.
\end{theo}

\noindent{\bf Proof: } Note that
$\textmd{Im}(\gamma)\subset \mathcal{M} \cap \mathbf{J}^{-1}(\mathcal{O}_\mu), $
and it is $G$-invariant, as well as
$\textmd{Im}(T\bar{\gamma}_{\mathcal{O}_\mu})\subset \mathcal{K}_{\mathcal{O}_\mu}, $ in
this case, $\omega_{\mathcal{K}_{\mathcal{O}_\mu}}\cdot
\tau_{\mathcal{K}_{\mathcal{O}_\mu}}= \tau_{\mathcal{K}_{\mathcal{O}_\mu}}\cdot
\omega_{\mathcal{M}_{\mathcal{O}_\mu}}= \tau_{\mathcal{K}_{\mathcal{O}_\mu}}\cdot
i_{\mathcal{M}_{\mathcal{O}_\mu}}^* \cdot \omega_{\mathcal{O}_\mu}, $ along
$\textmd{Im}(T\bar{\gamma}_{\mathcal{O}_\mu})$, and $\pi_{\mathcal{O}_\mu}^*\omega_{\mathcal{O}_\mu}=
i_{\mathcal{O}_\mu}^*\omega= \omega, $ along $\textmd{Im}(\gamma)$, and
$\tau_{\mathcal{K}_{\mathcal{O}_\mu}}\cdot T\bar{\gamma}_{\mathcal{O}_\mu}
= T\bar{\gamma}_{\mathcal{O}_\mu}, $ and
$\tau_{\mathcal{K}_{\mathcal{O}_\mu}}\cdot X_{h_{\mathcal{K}_{\mathcal{O}_\mu}}} = X_{\mathcal{K}_{\mathcal{O}_\mu}}. $
Since the dynamical vector fields $X_{H}$
and $X_{h_{\mathcal{K}_{\mathcal{O}_\mu}}}$ are $\pi_{\mathcal{O}_\mu}$-related, that is,
$X_{h_{\mathcal{K}_{\mathcal{O}_\mu}}}\cdot \pi_{\mathcal{O}_\mu}=
T\pi_{\mathcal{O}_\mu}\cdot X_{H}\cdot i_{\mathcal{O}_\mu}, $
using the $\mathbf{J}$-nonholonomic $R_o$-reduced distributional two-form
$\omega_{\mathcal{K}_{\mathcal{O}_\mu}}$, we have that
\begin{align*}
& \omega_{\mathcal{K}_{\mathcal{O}_\mu}}(T\bar{\gamma}_{\mathcal{O}_\mu}
\cdot X_H^\varepsilon- X_{\mathcal{K}_{\mathcal{O}_\mu}}\cdot \bar{\varepsilon}_{\mathcal{O}_\mu}, \;
\tau_{\mathcal{K}_{\mathcal{O}_\mu}}\cdot T\pi_{\mathcal{O}_\mu} \cdot w) \\
& = \omega_{\mathcal{K}_{\mathcal{O}_\mu}}(T\bar{\gamma}_{\mathcal{O}_\mu} \cdot X_H^\varepsilon, \;
\tau_{\mathcal{K}_{\mathcal{O}_\mu}}\cdot T\pi_{\mathcal{O}_\mu} \cdot w)-
\omega_{\mathcal{K}_{\mathcal{O}_\mu}}(X_{\mathcal{K}_{\mathcal{O}_\mu}}\cdot \bar{\varepsilon}_{\mathcal{O}_\mu},
\; \tau_{\mathcal{K}_{\mathcal{O}_\mu}} \cdot T\pi_{\mathcal{O}_\mu} \cdot w) \\
& = \omega_{\mathcal{K}_{\mathcal{O}_\mu}}(\tau_{\mathcal{K}_{\mathcal{O}_\mu}}
\cdot T\bar{\gamma}_{\mathcal{O}_\mu} \cdot X_H^\varepsilon, \;
\tau_{\mathcal{K}_{\mathcal{O}_\mu}}\cdot T\pi_{\mathcal{O}_\mu} \cdot w)-
\omega_{\mathcal{K}_{\mathcal{O}_\mu}}(\tau_{\mathcal{K}_{\mathcal{O}_\mu}}
\cdot X_{h_{\mathcal{K}_{\mathcal{O}_\mu}}} \cdot \pi_{\mathcal{O}_\mu} \cdot
\varepsilon, \; \tau_{\mathcal{K}_{\mathcal{O}_\mu}} \cdot T\pi_{\mathcal{O}_\mu} \cdot w) \\
& = \omega_{\mathcal{K}_{\mathcal{O}_\mu}}\cdot \tau_{\mathcal{K}_{\mathcal{O}_\mu}}
(T\pi_{\mathcal{O}_\mu} \cdot T\gamma \cdot X_H^\varepsilon, \; T\pi_{\mathcal{O}_\mu} \cdot w)
-\omega_{\mathcal{K}_{\mathcal{O}_\mu}}\cdot \tau_{\mathcal{K}_{\mathcal{O}_\mu}}
(T\pi_{\mathcal{O}_\mu} \cdot X_H\cdot \varepsilon, \; T\pi_{\mathcal{O}_\mu} \cdot w)\\
& = \tau_{\mathcal{K}_{\mathcal{O}_\mu}}\cdot
i_{\mathcal{M}_{\mathcal{O}_\mu}}^* \cdot \omega_{\mathcal{O}_\mu}
(T\pi_{\mathcal{O}_\mu} \cdot T\gamma \cdot X_H^\varepsilon, \; T\pi_{\mathcal{O}_\mu} \cdot w)
-\tau_{\mathcal{K}_{\mathcal{O}_\mu}}\cdot
i_{\mathcal{M}_{\mathcal{O}_\mu}}^* \cdot \omega_{\mathcal{O}_\mu}
(T\pi_{\mathcal{O}_\mu} \cdot X_H\cdot \varepsilon, \; T\pi_{\mathcal{O}_\mu} \cdot w)\\
& = \tau_{\mathcal{K}_{\mathcal{O}_\mu}}\cdot
i_{\mathcal{M}_{\mathcal{O}_\mu}}^* \cdot \pi_{\mathcal{O}_\mu}^*\omega_{\mathcal{O}_\mu}
(T\gamma \cdot X_H^\varepsilon, \; w)
-\tau_{\mathcal{K}_{\mathcal{O}_\mu}}\cdot
i_{\mathcal{M}_{\mathcal{O}_\mu}}^* \cdot \pi_{\mathcal{O}_\mu}^*\omega_{\mathcal{O}_\mu}(X_H\cdot \varepsilon, \; w)\\
& = \tau_{\mathcal{K}_{\mathcal{O}_\mu}}\cdot
i_{\mathcal{M}_{\mathcal{O}_\mu}}^* \cdot \omega(T\gamma \cdot X_H^\varepsilon, \; w)
-\tau_{\mathcal{K}_{\mathcal{O}_\mu}}\cdot
i_{\mathcal{M}_{\mathcal{O}_\mu}}^* \cdot \omega(X_H\cdot \varepsilon, \; w).
\end{align*}
In the case we considered that $\tau_{\mathcal{K}_{\mathcal{O}_\mu}}\cdot
i_{\mathcal{M}_{\mathcal{O}_\mu}}^* \cdot \omega=\tau_{\mathcal{K}}\cdot i_{\mathcal{M}}^* \cdot
\omega= \omega_{\mathcal{K}}\cdot \tau_{\mathcal{K}}, $
and
$\tau_{\mathcal{K}}\cdot T\gamma =T\gamma, \; \tau_{\mathcal{K}} \cdot X_H= X_{\mathcal{K}}$,
since $\textmd{Im}(\gamma)\subset
\mathcal{M}, $ and $\textmd{Im}(T\gamma)\subset \mathcal{K}. $
Thus, we have that
\begin{align*}
& \omega_{\mathcal{K}_{\mathcal{O}_\mu}}(T\bar{\gamma}_{\mathcal{O}_\mu}
\cdot X_H^\varepsilon- X_{\mathcal{K}_{\mathcal{O}_\mu}}\cdot \bar{\varepsilon}_{\mathcal{O}_\mu}, \;
\tau_{\mathcal{K}_{\mathcal{O}_\mu}}\cdot T\pi_{\mathcal{O}_\mu} \cdot w) \\
& = \omega_{\mathcal{K}}\cdot \tau_{\mathcal{K}}(T\gamma \cdot X_H^\varepsilon, \; w)
-\omega_{\mathcal{K}}\cdot \tau_{\mathcal{K}}(X_H\cdot \varepsilon, \; w)\\
& = \omega_{\mathcal{K}}(\tau_{\mathcal{K}}\cdot T\gamma \cdot X_H^\varepsilon, \; \tau_{\mathcal{K}}\cdot w)
-\omega_{\mathcal{K}}(\tau_{\mathcal{K}}\cdot X_H\cdot \varepsilon, \; \tau_{\mathcal{K}}\cdot w)\\
& = \omega_{\mathcal{K}}(T\gamma \cdot X_H^\varepsilon, \; \tau_{\mathcal{K}}\cdot w)
-\omega_{\mathcal{K}}(X_{\mathcal{K}}\cdot \varepsilon, \; \tau_{\mathcal{K}}\cdot w)\\
& = \omega_{\mathcal{K}}(T\gamma \cdot X_H^\varepsilon- X_{\mathcal{K}}\cdot \varepsilon, \; \tau_{\mathcal{K}}\cdot w).
\end{align*}
Because the distributional two-form $\omega_{\mathcal{K}}$ and
the $\mathbf{J}$-nonholonomic $R_o$-reduced distributional
two-form $\omega_{\mathcal{K}_{\mathcal{O}_\mu}}$ are both non-degenerate,
it follows that the equation
$T\bar{\gamma}_{\mathcal{O}_\mu}\cdot X_H^\varepsilon=
X_{\mathcal{K}_{\mathcal{O}_\mu}}\cdot \bar{\varepsilon}_{\mathcal{O}_\mu}, $ is equivalent to the equation
$T\gamma\cdot X_H^\varepsilon= X_{\mathcal{K}}\cdot \varepsilon$. Thus,
$\varepsilon$ is a solution of the Type II of Hamilton-Jacobi equation
$T\gamma\cdot X_H^\varepsilon= X_{\mathcal{K}}\cdot \varepsilon, $ for the distributional
Hamiltonian system $(\mathcal{K},\omega_{\mathcal {K}},H_{\mathcal {K}})$, if and only if
$\varepsilon$ and $\bar{\varepsilon}_{\mathcal{O}_\mu} $ satisfy the Type II of Hamilton-Jacobi
equation $T\bar{\gamma}_{\mathcal{O}_\mu}\cdot X_H^\varepsilon=
X_{\mathcal{K}_{\mathcal{O}_\mu}}\cdot \bar{\varepsilon}_{\mathcal{O}_\mu}, $ for the
$\mathbf{J}$-nonholonomic $R_o$-reduced distributional Hamiltonian system
$(\mathcal{K}_{\mathcal{O}_\mu},\omega_{\mathcal{K}_{\mathcal{O}_\mu}},h_{\mathcal{K}_{\mathcal{O}_\mu}})$.  \hskip 0.3cm
$\blacksquare$

\begin{rema}
If $(T^\ast Q, \omega)$ is a connected symplectic manifold, and
$\mathbf{J}:T^\ast Q\rightarrow \mathfrak{g}^\ast$ is a
non-equivariant momentum map with a non-equivariance group
one-cocycle $\sigma: G\rightarrow \mathfrak{g}^\ast$,
which is defined by $\sigma(g):=\mathbf{J}(g\cdot
z)-\operatorname{Ad}^\ast_{g^{-1}}\mathbf{J}(z)$, where $g\in G$ and
$z\in T^\ast Q$. Then we know that $\sigma$ produces a new affine
action $\Theta: G\times \mathfrak{g}^\ast \rightarrow
\mathfrak{g}^\ast $ defined by
$\Theta(g,\mu):=\operatorname{Ad}^\ast_{g^{-1}}\mu + \sigma(g)$,
where $\mu \in \mathfrak{g}^\ast$, with respect to which the given
momentum map $\mathbf{J}$ is equivariant. Assume that $G$ acts
freely and properly on $T^\ast Q$, and $\mathcal{O}_\mu= G\cdot \mu
\subset \mathfrak{g}^\ast$ denotes the G-orbit of the point $\mu \in
\mathfrak{g}^\ast$ with respect to this affine action $\Theta$,
and $\mu$ is a regular value of $\mathbf{J}$.
Then the quotient space $(T^\ast
Q)_{\mathcal{O}_\mu}=\mathbf{J}^{-1}(\mathcal{O}_\mu)/G$ is also a symplectic
manifold with symplectic form $\omega_{\mathcal{O}_\mu}$ uniquely characterized by
$(6.3)$, see Ortega and Ratiu \cite{orra04}. In this case,
we can also define the $\mathbf{J}$-nonholonomic regular orbit reducible
Hamiltonian system $(T^*Q,G,\omega,\mathbf{J},\mathcal{D},H)$ with an associated
$\mathbf{J}$-nonholonomic $R_o$-reduced distributional Hamiltonian system
$(\mathcal{K}_{\mathcal{O}_\mu},\omega_{\mathcal{K}_{\mathcal{O}_\mu}},
h_{\mathcal{K}_{\mathcal{O}_\mu}})$,
and prove the Type I and Type II of
the Hamilton-Jacobi theorem for the $\mathbf{J}$-nonholonomic $R_o$-reduced
distributional Hamiltonian system
$(\mathcal{K}_{\mathcal{O}_\mu},\omega_{\mathcal{K}_{\mathcal{O}_\mu}}, h_{\mathcal{K}_{\mathcal{O}_\mu}})$
by using the above similar way,
in which the $\mathbf{J}$-nonholonomic $R_o$-reduced space $((\mathcal{K}_{\mathcal{O}_\mu},\omega_{\mathcal{K}_{\mathcal{O}_\mu}})$
is determined by the affine action and $\mathbf{J}$-nonholonomic regular orbit reduction.
\end{rema}

\section{Applications}

In this section, in order to illustrate the Hamilton-Jacobi theory for
the nonholonomic reducible Hamiltonian system with symmetry, we shall discuss the
following two examples: (1) the motion of constrained particle in
space $\mathbb{R}^3$; (2) the motion of vertical rolling disk. These
two examples are classical in the theory of nonholonomic mechanical
systems. We shall follow the notations and conventions introduced in
Bates and $\acute{S}$niatycki \cite{basn93}, Bloch et al. \cite{blkrmamu96}
and Wang \cite{wa17}.

\subsection{The constrained particle in $\mathbb{R}^3$ }

In this subsection, we consider the motion of constrained particle
in space $\mathbb{R}^3$, and give explicitly the motion equations
of this problem and derive precisely
the geometric constraint conditions of the induced distributional two-forms
for the nonholonomic dynamical vector fields, that is,
the Type I and Type II of Hamilton-Jacobi equations of this problem.
At first, the configuration
space of motion of the constrained particle in space is $Q=
\mathbb{R}^3 $, whose coordinates are denoted by $q=(x,y,z), $ its
velocity space is $T\mathbb{R}^3 $, and the phase space is
$T^*\mathbb{R}^3 $ with canonical symplectic form
$$\omega=
\mathbf{d}x \wedge \mathbf{d}p_x + \mathbf{d}y \wedge \mathbf{d}p_y
+ \mathbf{d}z \wedge \mathbf{d}p_z .$$
The constraint set on the
velocities is given by $$\mathcal{D}=\{(x,y,z,v_x,v_y,v_z)\in
T\mathbb{R}^3 |\; v_z= \sigma(y)v_x \}, $$
where $\sigma(y)$ is a
smooth function. For any $q\in Q, \; \mathcal{D}(q)=
\textrm{Span}\{\partial_x+\sigma(y)\partial_z, \; \partial_y\}.$
Note that $$[\partial_x+\sigma(y)\partial_z, \; \partial_y]=
[\partial_x,\; \partial_y]+[\sigma(y)\partial_z,\; \partial_y]=
\sigma'(y)\partial_z, $$ which is nonzero everywhere if
$\sigma'(y)\neq 0$, then $\mathcal{D}$ is nonholonomic and it is
{\bf completely nonholonomic}, that is, $\mathcal{D}$ along with all
of its iterated Lie brackets $[\mathcal{D},\mathcal{D}],
[\mathcal{D},[\mathcal{D},\mathcal{D}]], \cdots $ spans the tangent
bundle $TQ$. The Lagrangian $L: T\mathbb{R}^3 \rightarrow
\mathbb{R}$ is the kinetic energy of the Euclidean metric of
$\mathbb{R}^3$, that is, $$L= \frac{1}{2}(v_x^2+ v_y^2+ v_z^2), $$
which is simple and it is hyperregular, and hence the system is {\bf
$\mathcal{D}$-regular } automatically. The momenta are $p_x=
\frac{\partial L}{\partial \dot{x}}= v_x, \; p_y= \frac{\partial
L}{\partial \dot{y}}= v_y,$ and $ p_z=\frac{\partial L}{\partial
\dot{z}}= v_z= \sigma(y)p_x, $ and the Hamiltonian $H:
T^*\mathbb{R}^3 \rightarrow \mathbb{R}$ is given by $$H=
\frac{1}{2}(p_x^2+ p_y^2+ p_z^2). $$ The unconstrained Hamiltonian
vector field is given by $$X_H= p_x\partial_x +p_y\partial_y +p_z\partial_z .$$ By
using the Legendre transformation
$$\mathcal{F}L: T\mathbb{R}^3 \rightarrow T^*\mathbb{R}^3 , \;
\mathcal{F}L(x,y,z,v_x,v_y,v_z)= (x,y,z,p_x,p_y,p_z), $$ the
constraint submanifold $\mathcal{M}= \mathcal{F}L(\mathcal{D})$ is
given by $$\mathcal{M}= \{(x,y,z,p_x,p_y,p_z)\in T^* Q |\; p_z=
\sigma(y)p_x\}. $$ Moreover, define
$\mathcal{F}=(T\pi_Q)^{-1}(\mathcal{D})$, and the compatibility
condition $T\mathcal{M}\cap \mathcal{F}^\bot= \{0\}$ holds, where
$\mathcal{F}^\bot$ denotes the symplectic orthogonal of
$\mathcal{F}$ with respect to the canonical symplectic form
$\omega$. Then the distribution is given by
$$\mathcal{K}= \mathcal{F} \cap T\mathcal{M}=
\mathrm{span}\{\partial_x + \sigma(y)\partial_z, \partial_y,
\partial_{p_x}, \partial_{p_y} \}. $$
The induced two-form $\omega_{\mathcal{M}}= i_{\mathcal{M}}^*\cdot
\omega $ is given by
$$\omega_\mathcal{M}= \mathbf{d}x \wedge \mathbf{d}p_x + \mathbf{d}y
\wedge \mathbf{d}p_y + \mathbf{d}z \wedge (p_x \sigma'(y)\mathbf{d}y
+ \sigma(y)\mathbf{d}p_x), $$ and the non-degenerate distributional
two-form is given by $$\omega_{\mathcal{K}}= \tau_{\mathcal{K}}\cdot
\omega_{\mathcal{M}}.$$ A direct computation yields
\begin{align*}
\mathbf{i}_{\partial_x + \sigma(y)\partial_z}\omega_\mathcal{K} & =
(1+ \sigma^2(y))\mathbf{d}p_x + \sigma(y)\sigma'(y)p_x\mathbf{d}y
,\;\;\;\;\;\; \mathbf{i}_{\partial_y}\omega_\mathcal{K} =
\mathbf{d}p_y-
\sigma'(y)p_x\mathbf{d}z, \\
\mathbf{i}_{\partial_{p_x}}\omega_\mathcal{K} & =
-\sigma(y)\mathbf{d}z -\mathbf{d}x,
\;\;\;\;\;\;\;\;\;\;\;\;\;\;\;\;\;\;\;\;\;\;\;\;\;\;\;\;\;\;\;\;\;
\mathbf{i}_{\partial_{p_y}}\omega_\mathcal{K} = -\mathbf{d}y,
\end{align*} and
\begin{align*}
 \mathbf{d}H_\mathcal{K} & = p_x\mathbf{d}p_x + p_y\mathbf{d}p_y +
\sigma(y)\sigma'(y)p_x^2\mathbf{d}y + \sigma^2(y)p_x\mathbf{d}p_x
\\ & = \sigma(y)\sigma'(y)p_x^2\mathbf{d}y + (1+
\sigma^2(y))p_x\mathbf{d}p_x + p_y\mathbf{d}p_y .\end{align*} Assume
that $X_\mathcal{K}= X_1(\partial_x + \sigma(y)\partial_z)+ X_2
\partial_y + X_3 \partial_{p_x} + X_4 \partial_{p_y}, $ then we have
that
\begin{align*}
 \mathbf{i}_{X_\mathcal{K}}\omega_\mathcal{K} & = X_1((1+
\sigma^2(y))\mathbf{d}p_x + \sigma(y)\sigma'(y)p_x\mathbf{d}y ) \\&
\;\;\; + X_2(\mathbf{d}p_y- \sigma'(y) p_x\mathbf{d}z)
+ X_3(-\sigma(y)\mathbf{d}z -\mathbf{d}x) +X_4(-\mathbf{d}y)\\
& = (-X_3) \mathbf{d}x +(\sigma(y)\sigma'(y)p_x X_1-X_4) \mathbf{d}y
\\& \;\;\; + (-\sigma'(y)p_x X_2-\sigma(y)X_3)\mathbf{d}z +
(1+\sigma^2(y))X_1 \mathbf{d}p_x + X_2 \mathbf{d}p_y.
\end{align*} From the equation of distributional Hamiltonian system
$\mathbf{i}_{X_\mathcal{K}}\omega_\mathcal{K}=
\mathbf{d}H_\mathcal{K}, $ we have that $$ X_1= p_x,\;\;\;\;
X_2=p_y, \;\;\;\;\;  X_3=0, \;\;\;\;\; X_4= 0. $$ Hence, the
nonholonomic dynamical vector field is given by
$$X_\mathcal{K}= p_x(\partial_x +
\sigma(y)\partial_z) + p_y \partial_y, $$ and the motion equations of the
distributional Hamiltonian system
$(\mathcal{K},\omega_{\mathcal{K}},H )$ are expressed by
$$ \dot{x}= p_x, \;\;\; \dot{y}=p_y,\;\;\; \dot{z}
= \sigma(y)p_x, \;\;\; \dot{p}_x=0, \;\;\; \dot{p}_y=0.$$

In the following we shall derive precisely
the geometric constraint conditions of the induced distributional two-form
for the nonholonomic dynamical vector field, that is,
the Type I and Type II of Hamilton-Jacobi equations for the
distributional Hamiltonian system
$(\mathcal{K},\omega_{\mathcal{K}},H )$. Assume that $$\gamma:
\mathbb{R}^3 \rightarrow T^*\mathbb{R}^3, \;\;\; \gamma(x,y,z)= (
\gamma_1, \gamma_2, \gamma_3, \gamma_4, \gamma_5, \gamma_6), $$
and $\lambda=\gamma \cdot \pi_{Q}: T^* \mathbb{R}^3
\rightarrow T^* \mathbb{R}^3 $ given by
\begin{align*}
& \lambda(x,y,z,p_x,p_y,p_z)= (\lambda_1, \lambda_2, \lambda_3, \lambda_4, \lambda_5, \lambda_6 )\\
& = \gamma \cdot \pi_Q (x,y,z,p_x,p_y,p_z)=\gamma(x,y,z)\\
& =(\gamma_1 \cdot \pi_Q,\gamma_2 \cdot \pi_Q,\gamma_3 \cdot \pi_Q,\gamma_4 \cdot \pi_Q,\gamma_5 \cdot \pi_Q,\gamma_6 \cdot \pi_Q),
\end{align*} that is, $\lambda_i= \gamma_i \cdot \pi_Q, \; i=1, \cdots, 6,$
where $\lambda_i, \; i=1,\cdots, 6, $ are functions on
$T^*\mathbb{R}^3$, and $\gamma_i, \; i=1,\cdots, 6, $ are functions on
$\mathbb{R}^3$. We may choose $q=(x,y,z)\in \mathbb{R}^3, $ such that
$\gamma_1(q)=x, \; \gamma_2(q)=y, \; \gamma_3(q)=z, $ and
$\gamma(q)= \gamma_4(q)\mathbf{d}x+ \gamma_5(q)\mathbf{d}y+ \gamma_6(q)\mathbf{d}z. $
Note that $\mathcal{D}(q)=
\textrm{Span}\{\partial_x+\sigma(y)\partial_z, \; \partial_y\}, $
take that $\alpha= \partial_x+\sigma(y)\partial_z $ and $\beta= \partial_y, $ then we have that
\begin{align*}
\mathbf{d}\gamma(\alpha,\beta)& = \alpha(\gamma(\beta))-\beta(\gamma(\alpha))-\gamma([\alpha,\beta])\\
& = (\frac{\partial \gamma_5}{\partial x}- \frac{\partial \gamma_4}{\partial y})
- \sigma(y)(\frac{\partial \gamma_6}{\partial y}- \frac{\partial \gamma_5}{\partial z})
- 2\sigma'(y)\gamma_6.
\end{align*}
Thus, when $\mathbf{d}\gamma(\alpha,\beta)=0, $ we know that for any
$v, w \in \mathcal{F}, $ and $T\pi_{Q}(v), \; T\pi_{Q}(w) \in \mathcal{D},$  then
$\mathbf{d}\gamma(T\pi_{Q}(v),T\pi_{Q}(w))=0, $ that is,
$\gamma$ is closed on $\mathcal{D}$ with respect to $T\pi_{Q}:
TT^* \mathbb{R}^3 \rightarrow T\mathbb{R}^3. $
Note that $Im(\gamma) \subset \mathcal{M}, $
then $p_x=\gamma_4, \; p_y=\gamma_5, \; p_z=\gamma_6$ and
$\gamma_6=\sigma(y)\gamma_4, $ and hence
$$H\cdot \gamma= \frac{1}{2}((1+ \sigma^2(y))\gamma_4^2+
\gamma_5^2), $$
$$X_H \cdot \gamma= \gamma_4\partial_x + \gamma_5\partial_y + \sigma(y)\gamma_4\partial_z
= X_{\mathcal{K}} \cdot \gamma, $$
$$ X_H^\gamma= T\pi_Q \cdot
X_H \cdot \gamma= \gamma_4\partial_x + \gamma_5\partial_y + \sigma(y)\gamma_4\partial_z .$$
Thus, $T\gamma \cdot X_H^\gamma= X_{\mathcal{K}} \cdot \gamma, $
that is, the Type I of Hamilton-Jacobi equation for the
distributional Hamiltonian system
$(\mathcal{K},\omega_{\mathcal{K}},H )$ holds trivially.\\

Now, for any symplectic map $\varepsilon: T^* \mathbb{R}^3
\rightarrow T^* \mathbb{R}^3, $ from $\omega= \varepsilon^*\omega=\omega \cdot \varepsilon
= (\partial_x\varepsilon \cdot \partial_{p_x}\varepsilon) \mathbf{d}x \wedge \mathbf{d}p_x
+ (\partial_y\varepsilon \cdot \partial_{p_y}\varepsilon) \mathbf{d}y \wedge \mathbf{d}p_y
+ (\partial_z\varepsilon \cdot \partial_{p_z}\varepsilon) \mathbf{d}z \wedge \mathbf{d}p_z, $ we have that
\begin{align*}
 \partial_x \varepsilon\cdot \partial_{p_x}\varepsilon= 1, \;\;\;\;\;\;
 \partial_y \varepsilon\cdot \partial_{p_y}\varepsilon= 1, \;\;\;\;\;\;
 \partial_z \varepsilon\cdot \partial_{p_z}\varepsilon= 1.
\end{align*}
Denote by $\varepsilon(x,y,z,p_x,p_y,p_z)
=(\varepsilon_1, \varepsilon_2, \varepsilon_3, \varepsilon_4, \varepsilon_5, \varepsilon_6 ),$
then we have that
$$H\cdot \varepsilon= \frac{1}{2}(\varepsilon_4^2+
\varepsilon_5^2+\varepsilon_6^2), $$ and
$$X_H \cdot \varepsilon= \varepsilon_4\partial_x + \varepsilon_5\partial_y + \varepsilon_6\partial_z ,$$
and hence $$ X_H^\varepsilon= T\pi_Q \cdot
X_H \cdot \varepsilon= \varepsilon_4\partial_x +\varepsilon_5\partial_y +\varepsilon_6\partial_z .$$
Since $Im(\gamma) \subset \mathcal{M}, $ then
\begin{align*}
T\gamma\cdot X_H^\varepsilon
& = \varepsilon_4\partial_x +\varepsilon_5\partial_y +\sigma(y)\varepsilon_4\partial_z\\
& = \varepsilon_4(\partial_x +\sigma(y)\partial_z)+\varepsilon_5\partial_y= X_{\mathcal{K}}\cdot\varepsilon,
\end{align*}
because $\varepsilon_6=\sigma(y)\varepsilon_4. $ In the same way, note that
$\lambda=\gamma\cdot \pi_Q, $ and $Im(\lambda) \subset \mathcal{M}, $ then
$$
T\lambda\cdot X_H \cdot \varepsilon=\varepsilon_4\partial_x +\varepsilon_5\partial_y
+\sigma(y)\varepsilon_4\partial_z = X_{\mathcal{K}}\cdot\varepsilon.
$$
On the other hand, since $\varepsilon: T^* Q \rightarrow T^* Q $ is symplectic, we have that
\begin{align*}
\tau_{\mathcal{K}}\cdot T\varepsilon\cdot X_{H \cdot \varepsilon}
& =\tau_{\mathcal{K}}\cdot X_H \cdot \varepsilon\\
& =\varepsilon_4\partial_x +\varepsilon_5\partial_y +\sigma(y)\varepsilon_4\partial_z
= X_{\mathcal{K}}\cdot\varepsilon.
\end{align*}
Thus, $T\gamma\cdot X_H^\varepsilon=X_{\mathcal{K}}\cdot\varepsilon
=T\lambda\cdot X_H \cdot \varepsilon
=\tau_{\mathcal{K}}\cdot T\varepsilon\cdot X_{H \cdot \varepsilon}.$
In this case, we must have that
$\varepsilon$ is a
solution of the Type II of Hamilton-Jacobi equation $T\gamma \cdot X_H^\varepsilon= X_{\mathcal{K}}
\cdot \varepsilon, $ for the distributional Hamiltonian system
$(\mathcal{K},\omega_{\mathcal{K}},H ), $ if and only if it is a solution of the equation
$T\lambda\cdot X_H \cdot \varepsilon
=\tau_{\mathcal{K}}\cdot T\varepsilon\cdot X_{H \cdot \varepsilon}$.\\

Next, we consider the action of Lie group $G= \mathbb{R}^2$ on
$\mathbb{R}^3$, and derive precisely the motion equations and
the geometric constraint conditions of the reduced distributional two-form
for the nonholonomic reduced dynamical vector field, that is,
the Type I and Type II of Hamilton-Jacobi
equations of the nonholonomic reduced distributional Hamiltonian
system. At first, the action of Lie group $G= \mathbb{R}^2$ on
$\mathbb{R}^3$ is given by
$$\Phi: G\times \mathbb{R}^3 \rightarrow \mathbb{R}^3, \;
\Phi((r,s),(x,y,z))= (x+r,y,z+s), $$ and we have the cotangent lifted
$G$-action on $T^* \mathbb{R}^3$, such that the Hamiltonian $H: T^*
\mathbb{R}^3 \rightarrow \mathbb{R}$ is $G$-invariant. Therefore,
$$\bar{\mathcal{M}}= \{(y,p_x,p_y,p_z)\in T^*\mathbb{R}^3 /G|\; p_z=
\sigma(y)p_x\}, $$ and the reduced distribution is given by
$$\bar{\mathcal{K}}=
\mathrm{span}\{(1+\sigma^2(y))\partial_y- \sigma(y)\sigma'(y)p_x
\partial_{p_x}, \partial_{p_y} \}, $$
and the non-degenerate and the reduced
distributional two-form $\omega_{\bar{\mathcal{K}}}$
is given by $$\omega_{\bar{\mathcal{K}}}= \mathbf{d}x \wedge
\mathbf{d}p_x + \mathbf{d}y \wedge \mathbf{d}p_y + \mathbf{d}z
\wedge (p_x \sigma'(y)\mathbf{d}y + \sigma(y)\mathbf{d}p_x).$$
A direct computation yields
\begin{align*}
\mathbf{i}_{(1+\sigma^2(y))\partial_y- \sigma(y)\sigma'(y)p_x
\partial_{p_x}}\omega_{\bar{\mathcal{K}}}=& \sigma(y)\sigma'(y)p_x
\mathbf{d}x -\sigma'(y)p_x \mathbf{d}z
+(1+\sigma^2(y))\mathbf{d}p_y, \\
\mathbf{i}_{\partial_{p_y}}\omega_{\bar{\mathcal{K}}}=& -\mathbf{d}y
\end{align*}
and $$\mathbf{d}h_{\bar{\mathcal{K}}}= \mathbf{d}H_\mathcal{K}=
\sigma(y)\sigma'(y)p_x^2\mathbf{d}y + (1+
\sigma^2(y))p_x\mathbf{d}p_x + p_y\mathbf{d}p_y .$$ Assume that
$X_{\bar{\mathcal{K}}}= X_1((1+\sigma^2(y))\partial_y-
\sigma(y)\sigma'(y)p_x
\partial_{p_x}) + X_2 \partial_{p_y}, $ then we have that
$$\mathbf{i}_{X_{\bar{\mathcal{K}}}}\omega_{\bar{\mathcal{K}}}= X_1(\sigma(y)\sigma'(y)p_x
\mathbf{d}x -\sigma'(y)p_x \mathbf{d}z
+(1+\sigma^2(y))\mathbf{d}p_y) - X_2 \mathbf{d}z
$$
$$ = (X_1\sigma(y)\sigma'(y)p_x) \mathbf{d}x +(-X_2) \mathbf{d}y
+ (-X_1\sigma'(y)p_x)\mathbf{d}z + (X_1(1+\sigma^2(y))) \mathbf{d}p_y .$$
From the nonholonomic reduced distributional Hamiltonian equation
$\mathbf{i}_{X_{\bar{\mathcal{K}}}}\omega_{\bar{\mathcal{K}}}=
\mathbf{d}h_{\bar{\mathcal{K}}}, $ we have that $X_1= 0, \;\; X_2=
-\sigma(y)\sigma'(y)p_x^2. $ Hence, we get that the nonholonomic
reduced dynamical vector field is given by $$X_{\bar{\mathcal{K}}}=
-\sigma(y)\sigma'(y)p_x^2 \partial_{p_y}, $$
and the motion equations of the nonholonomic reduced distributional Hamiltonian
system $(\bar{\mathcal{K}},\omega_{\bar{\mathcal{K}}},h )$ are
expressed by $$ \dot{y}=0,\;\;\; \dot{p}_x=0, \;\;\; \dot{p}_y=
-\sigma(y)\sigma'(y)p_x^2. $$

In the following we shall derive precisely
the Type I and Type II of Hamilton-Jacobi equations of the
nonholonomic reduced distributional Hamiltonian system
$(\bar{\mathcal{K}},\omega_{\bar{\mathcal{K}}},h ). $ Assume that
$\gamma: \mathbb{R}^3 \rightarrow T^*\mathbb{R}^3, $ and $\lambda=
\gamma \cdot \pi_{Q}: T^*\mathbb{R}^3 \rightarrow T^*\mathbb{R}^3, $ and
$\textmd{Im}(\gamma)\subset \mathcal{M}, $ and it is
$G$-invariant, $ \textmd{Im}(T\gamma)\subset \mathcal{K}, $
then we have that $\bar{\gamma}=\pi_{/G}(\gamma): \mathbb{R}^3
\rightarrow T^*\mathbb{R}^3/G, \; \bar{\gamma}(x,y,z)=
(\bar{\gamma}_0,\bar{\gamma}_1,\bar{\gamma}_2,\bar{\gamma}_3)$, and
$\bar{\lambda}=\pi_{/G}(\lambda): T^* \mathbb{R}^3 \rightarrow T^*
\mathbb{R}^3/G, \; \bar{\lambda}(x,y,z,p_x,p_y,p_z)=
(\bar{\lambda}_0,\bar{\lambda}_1,\bar{\lambda}_2,\bar{\lambda}_3), $
that is, $\bar{\lambda}_i= \bar{\gamma}_i \cdot \pi_Q, \; i=0,\cdots,3, $
where $\bar{\lambda}_i, \; i=0,\cdots,3, $ are functions on $T^*\mathbb{R}^3$, and
$\bar{\gamma}_i, \; i=0,\cdots,3, $ are functions on $\mathbb{R}^3$.
Note that $h\cdot \pi_{/G}= \tau_{\mathcal{M}}\cdot H, $
since $Im(\gamma) \subset \mathcal{M}, $ and it is $G$-invariant, we
have that $Im(\bar{\gamma}) \subset \bar{\mathcal{M}}, $ and
$\bar{\gamma}_3=\sigma(y)\bar{\gamma}_1, $ and hence
$$ h\cdot \bar{\gamma}= \frac{1}{2}((1+
\sigma^2(y))\bar{\gamma}_1^2+ \bar{\gamma}_2^2), $$ and
$$
X_h\cdot \bar{\gamma}=\bar{\gamma}_2\partial_y-\sigma(y)\sigma'(y)\bar{\gamma}_1^2
\partial_{p_y}.
$$
When $\mathbf{d}\gamma(\alpha,\beta)=0, $ that is,
$\gamma$ is closed on $\mathcal{D}$ with respect to $T\pi_{Q}:
TT^* \mathbb{R}^3 \rightarrow T\mathbb{R}^3, $ we have that
$$
T\bar{\gamma}\cdot X^\gamma_H= \tau_{\bar{\mathcal{K}}}\cdot X_h\cdot \bar{\gamma}
=-\sigma(y)\sigma'(y)\bar{\gamma}_1^2
\partial_{p_y}=X_{\bar{\mathcal{K}}}\cdot \bar{\gamma},
$$
that is, the Type I of Hamilton-Jacobi equation for
the nonholonomic reduced distributional Hamiltonian system
$(\bar{\mathcal{K}},\omega_{\bar{\mathcal{K}}},h ) $ holds.\\

Now, for any $G$-invariant symplectic map $\varepsilon: T^* \mathbb{R}^3
\rightarrow T^* \mathbb{R}^3, $ $\bar{\varepsilon}=\pi_{/G}(\varepsilon): T^* \mathbb{R}^3 \rightarrow T^*
\mathbb{R}^3/G, $ is given by $\bar{\varepsilon}(x,y,z,p_x,p_y,p_z)=
(\bar{\varepsilon}_0,\bar{\varepsilon}_1,\bar{\varepsilon}_2,\bar{\varepsilon}_3), $ then we have that
$$ h\cdot \bar{\varepsilon}= \frac{1}{2}((1+
\sigma^2(y))\bar{\varepsilon}_1^2+ \bar{\varepsilon}_2^2), $$ and
$$
X_h\cdot \bar{\varepsilon}=\bar{\varepsilon}_2\partial_y-\sigma(y)\sigma'(y)\bar{\varepsilon}_1^2
\partial_{p_y}.
$$
Because $Im(\bar{\gamma}) \subset \bar{\mathcal{M}}, $ and $Im(T\bar{\gamma})\subset \bar{\mathcal{K}}, $
and hence
$$
T\bar{\gamma}\cdot X_H^\varepsilon=\tau_{\bar{\mathcal{K}}}\cdot X_h\cdot \bar{\varepsilon}
=-\sigma(y)\sigma'(y)\bar{\varepsilon}_1^2
\partial_{p_y} =X_{\bar{\mathcal{K}}}\cdot \bar{\varepsilon}.
$$
Note that $\bar{\lambda}= \bar{\gamma} \cdot \pi_Q, $ and $Im(\bar{\lambda}) \subset \bar{\mathcal{M}}, $
and $Im(T\bar{\lambda})\subset \bar{\mathcal{K}}, $
then we have that
$$
T\bar{\lambda}\cdot X_H \cdot \varepsilon=\tau_{\bar{\mathcal{K}}}\cdot X_h\cdot \bar{\varepsilon}
=X_{\bar{\mathcal{K}}}\cdot \bar{\varepsilon}.
$$
On the other hand, since $\varepsilon: T^* \mathbb{R}^3 \rightarrow T^* \mathbb{R}^3 $ is symplectic, and
$\bar{\varepsilon}^*= \varepsilon^*\cdot \pi_{/G}^*: T^*(T^* \mathbb{R}^3)/G
\rightarrow T^*T^* \mathbb{R}^3$ is also symplectic along $\bar{\varepsilon}$, then we have that
\begin{align*}
\tau_{\bar{\mathcal{K}}}\cdot T\bar{\varepsilon}\cdot X_{h \cdot \bar{\varepsilon}} &
= \tau_{\bar{\mathcal{K}}}\cdot X_h\cdot \bar{\varepsilon}\\
& = -\sigma(y)\sigma'(y)\bar{\varepsilon}_1^2
\partial_{p_y}= X_{\bar{\mathcal{K}}}\cdot \bar{\varepsilon}.
\end{align*}
Thus, $T\bar{\gamma}\cdot X_H^\varepsilon=X_{\bar{\mathcal{K}}}\cdot\bar{\varepsilon}
=T\bar{\lambda}\cdot X_H \cdot \varepsilon=\tau_{\bar{\mathcal{K}}}\cdot T\bar{\varepsilon}\cdot X_{h \cdot \bar{\varepsilon}}.$
In this case, we must have that
$\varepsilon$ and $\bar{\varepsilon}$ are the
solution of the Type II of Hamilton-Jacobi equation
$T\bar{\gamma} \cdot X_H^\varepsilon= X_{\bar{\mathcal{K}}}
\cdot \bar{\varepsilon}, $ for the nonholonomic reduced distributional Hamiltonian system
$(\bar{\mathcal{K}},\omega_{\bar{\mathcal{K}}},h ), $ if and only if they satisfy the equation
$T\bar{\lambda}\cdot X_H \cdot \varepsilon=\tau_{\bar{\mathcal{K}}}\cdot T\bar{\varepsilon}\cdot X_{h \cdot \bar{\varepsilon}}$.

\subsection{The vertical rolling disk }

In this subsection, we consider the motion of a vertical rolling disk,
and give explicitly the motion equations of this problem and
derive precisely
the geometric constraint conditions of the induced distributional two-forms
for the nonholonomic dynamical vector fields,
that is, the Type I and Type II of Hamilton-Jacobi
equations of this problem. Assume that a vertical disk of zero width
rolls without slipping on a horizontal plane and it rotates freely
about its vertical axis. Let $x$ and $y$ denote the position of
contact point of the disk in the plane, and the variables $\theta$
and $\varphi$ denote the orientations of a chosen material point
with respect to the vertical plane and the "heading angle" of the
disk, see \cite{blkrmamu96}. Thus, the configuration space of motion
for the vertical rolling disk is $Q= \mathbb{R}^2 \times
\mathbb{S}^1 \times \mathbb{S}^1 $ whose coordinates are denoted by
$q=(x,y,\theta,\varphi), $ and its velocity space is $TQ$, and the
phase space is $T^* Q$ with canonical symplectic form
$$\omega= \mathbf{d}x \wedge \mathbf{d}p_x + \mathbf{d}y \wedge \mathbf{d}p_y
+ \mathbf{d}\theta \wedge \mathbf{d}p_\theta+ \mathbf{d}\varphi \wedge \mathbf{d}p_\varphi .$$
The rolling constraint set on the velocities is given by
$$\mathcal{D}=\{(x,y,\theta,\varphi,\dot{x},\dot{y},\dot{\theta},\dot{\varphi})\in
TQ |\; \dot{x}= R\dot{\theta}\cos\varphi, \; \dot{y}=
R\dot{\theta}\sin\varphi \}, $$ where $R$ denotes the radius of the
disk. For any $q\in Q, $ we have that $$\mathcal{D}(q)=
\textrm{Span}\{R\cos\varphi\partial_x+ R\sin\varphi\partial_y+
\partial_{\theta}, \; \partial_{\varphi}\}. $$
Note that
\begin{align*} [R\cos\varphi\partial_x+ R\sin\varphi\partial_y+
\partial_{\theta}, \;
\partial_{\varphi}] &= [R\cos\varphi\partial_x,\; \partial_{\varphi}]+[R\sin\varphi\partial_y,\;
\partial_{\varphi}]+ [\partial_{\theta}, \; \partial_{\varphi}] \\ &= -R\sin\varphi\partial_x + R\cos\varphi\partial_y,
\end{align*} which is nonzero everywhere and it is not in
$\mathcal{D}$, then $\mathcal{D}$ is nonholonomic and it is {\bf
completely nonholonomic}, that is, $\mathcal{D}$ along with all of
its iterated Lie brackets $[\mathcal{D},\mathcal{D}],
[\mathcal{D},[\mathcal{D},\mathcal{D}]], \cdots $ spans the tangent
bundle $TQ$. The Lagrangian $L: TQ \rightarrow \mathbb{R}$ is the
kinetic energy, that is, $$L= \frac{1}{2}m(\dot{x}^2+ \dot{y}^2) +
\frac{1}{2}I\dot{\theta}^2 + \frac{1}{2}J\dot{\varphi}^2, $$ where
$m$ is the mass of the disk, and $I$ and $J$ are its moments of
inertia. Note that $L$ is simple and it is hyperregular, and hence
the system is {\bf $\mathcal{D}$-regular } automatically. The
momenta are $p_x= \frac{\partial L}{\partial \dot{x}}= m\dot{x}, \;
p_y= \frac{\partial L}{\partial \dot{y}}= m\dot{y}, \; p_\theta=
\frac{\partial L}{\partial \dot{\theta}}= I\dot{\theta}, \;
p_\varphi= \frac{\partial L}{\partial \dot{\varphi}}=
J\dot{\varphi}, $ and the Hamiltonian $H: T^*Q \rightarrow
\mathbb{R}$ is given by $$H= \frac{1}{2m}(p_x^2 + p_y^2)+
\frac{1}{2I}p_\theta^2 + \frac{1}{2J}p_\varphi^2. $$ The
unconstrained Hamiltonian vector field is given by $$X_H=\frac{1}{m}
p_x\partial_x + \frac{1}{m}p_y\partial_y +
\frac{1}{I}p_\theta\partial_\theta +
\frac{1}{J}p_\varphi\partial_\varphi .$$ By using the Legendre
transformation
$$\mathcal{F}L: TQ \rightarrow T^* Q, \;
\mathcal{F}L(x,y,\theta,\varphi,\dot{x},\dot{y},\dot{\theta},\dot{\varphi})=
(x,y,\theta,\varphi,p_x,p_y,p_\theta,p_\varphi), $$ we obtain the constraint
submanifold $\mathcal{M}= \mathcal{F}L(\mathcal{D})$ given by
$$\mathcal{M}= \{(x,y,\theta,\varphi,p_x,p_y,p_\theta,p_\varphi)\in
T^* Q |\; p_x= \frac{mR}{I}p_\theta \cos\varphi, \; p_y=
\frac{mR}{I}p_\theta \sin\varphi \}. $$ Moreover, if we define
$\mathcal{F}=(T\pi_Q)^{-1}(\mathcal{D})$, then the compatibility
condition $T\mathcal{M}\cap \mathcal{F}^\bot= \{0\}$ holds, where
$\mathcal{F}^\bot$ denotes the symplectic orthogonal of
$\mathcal{F}$ with respect to the canonical symplectic form
$\omega$. Thus, the distribution is given by
$$\mathcal{K}= \mathcal{F} \cap T\mathcal{M}=
\mathrm{span}\{\partial_\theta+ R\cos\varphi
\partial_x+ R\sin\varphi \partial_y,
\partial_\varphi, \partial_{p_\theta}, \partial_{p_\varphi} \}. $$
The induced two-form $\omega_{\mathcal{M}}=
i_{\mathcal{M}}^*\cdot\omega $ is given by
$$\omega_\mathcal{M}= \mathbf{d}x \wedge
(\frac{mR\cos\varphi}{I}\mathbf{d}p_\theta -
\frac{mR\sin\varphi}{I}p_\theta\mathbf{d}\varphi) + \mathbf{d}y
\wedge (\frac{mR\sin\varphi}{I}\mathbf{d}p_\theta +
\frac{mR\cos\varphi}{I}p_\theta \mathbf{d}\varphi) +
\mathbf{d}\theta \wedge \mathbf{d}p_\theta + \mathbf{d}\varphi
\wedge \mathbf{d}p_\varphi, $$ and hence we have the non-degenerate distributional
two-form $\omega_{\mathcal{K}}= \tau_{\mathcal{K}}\cdot
\omega_{\mathcal{M}}$. A direct computation yields
\begin{align*}
& \mathbf{i}_{\partial_\theta+ R\cos\varphi
\partial_x+ R\sin\varphi \partial_y}\omega_\mathcal{K}= (1+ \frac{mR^2}{I})\mathbf{d}p_\theta,
\;\;\;\;\;\;\;\;\;\;\;
\mathbf{i}_{\partial_\varphi}\omega_\mathcal{K}=
\mathbf{d}p_\varphi+ \frac{mR\sin\varphi}{I}p_\theta \mathbf{d}x -
\frac{mR\cos\varphi}{I}p_\theta\mathbf{d}y, \\ &
\mathbf{i}_{\partial_{p_\theta}}\omega_\mathcal{K}=
-\frac{mR\cos\varphi}{I}\mathbf{d}x -
\frac{mR\sin\varphi}{I}\mathbf{d}y - \mathbf{d}\theta,
\;\;\;\;\;\;\;\;\;\;\;
\mathbf{i}_{\partial_{p_\varphi}}\omega_\mathcal{K}=
-\mathbf{d}\varphi,
\end{align*} and
\begin{align*}
 \mathbf{d}H_\mathcal{K} & = \frac{1}{m}\frac{mR}{I}p_\theta
\cos\varphi (\frac{mR\cos\varphi}{I}\mathbf{d}p_\theta -
\frac{mR\sin\varphi}{I}p_\theta\mathbf{d}\varphi)\\ & +
\frac{1}{m}\frac{mR}{I}p_\theta \sin\varphi
(\frac{mR\sin\varphi}{I}\mathbf{d}p_\theta +
\frac{mR\cos\varphi}{I}p_\theta \mathbf{d}\varphi) +
\frac{1}{I}p_\theta \mathbf{d}p_\theta + \frac{1}{J}p_\varphi
\mathbf{d}p_\varphi \\ & = \frac{1}{I}(1+
\frac{mR^2}{I})p_\theta\mathbf{d}p_\theta + \frac{1}{J}p_\varphi
\mathbf{d}p_\varphi .\end{align*} Assume that $X_\mathcal{K}=
X_1(\partial_\theta+ R\cos\varphi
\partial_x+ R\sin\varphi \partial_y)+ X_2
\partial_\varphi + X_3 \partial_{p_\theta} + X_4 \partial_{p_\varphi}, $ then
\begin{align*}
\mathbf{i}_{X_\mathcal{K}}\omega_\mathcal{K} & = X_1((1+
\frac{mR^2}{I})\mathbf{d}p_\theta ) + X_2(\mathbf{d}p_\varphi+
\frac{mR\sin\varphi}{I}p_\theta \mathbf{d}x -
\frac{mR\cos\varphi}{I}p_\theta\mathbf{d}y)\\
& \;\;\;\;\; + X_3(-\frac{mR\cos\varphi}{I}\mathbf{d}x -
\frac{mR\sin\varphi}{I}\mathbf{d}y - \mathbf{d}\theta ) +
X_4(-\mathbf{d}\varphi )\\
& = (X_2\frac{mR\sin\varphi}{I}p_\theta- X_3
\frac{mR\cos\varphi}{I}) \mathbf{d}x
+(-X_2\frac{mR\cos\varphi}{I}p_\theta- X_3 \frac{mR\sin\varphi}{I})
\mathbf{d}y\\ & \;\;\;\;\; + (-X_3)\mathbf{d}\theta +
(-X_4)\mathbf{d}\varphi + (X_1 (1+ \frac{mR^2}{I}))
\mathbf{d}p_\theta + (X_2) \mathbf{d}p_\varphi . \end{align*} From
the distributional Hamiltonian equation
$\mathbf{i}_{X_\mathcal{K}}\omega_\mathcal{K}=
\mathbf{d}H_\mathcal{K}, $ we have that
$$X_1=\frac{1}{I}p_\theta, \;\;\;  X_2=\frac{1}{J}p_\varphi,
\;\;\; X_3= 0, \;\;\; X_4= 0.$$ Hence, we get that the nonholonomic
dynamical vector field is given by
$$X_\mathcal{K}=\frac{1}{I}p_\theta(\partial_\theta+ R\cos\varphi
\partial_x+ R\sin\varphi \partial_y) +
\frac{1}{J}p_\varphi\partial_{\varphi}, $$ and the motion equations
of the distributional Hamiltonian system
$(\mathcal{K},\omega_{\mathcal{K}},H )$ are given by
$$\dot{x}= \frac{R\cos\varphi}{I}p_\theta, \;\;\; \dot{y}= \frac{R\sin\varphi}{I}p_\theta,
\;\;\; \dot{\theta}=\frac{1}{I}p_\theta,\;\;\; \dot{\varphi}=
\frac{1}{J}p_\varphi, \;\;\; \dot{p}_{\theta}=0, \;\;\;
\dot{p}_{\varphi}=0. $$

In the following we shall derive precisely
the geometric constraint conditions of the induced distributional two-forms
for the nonholonomic dynamical vector fields, that is,
the Type I and Type II of Hamilton-Jacobi equations for the
distributional Hamiltonian system
$(\mathcal{K},\omega_{\mathcal{K}},H )$. Assume that $$\gamma: Q
\rightarrow T^* Q, \;\;\; \gamma(x,y,\theta,\varphi)= (
\gamma_1, \gamma_2, \gamma_3, \gamma_4, \gamma_5, \gamma_6, \gamma_7, \gamma_8), $$
then $\lambda=\gamma
\cdot \pi_{Q}: T^* Q \rightarrow T^* Q $ given by
\begin{align*}
& \lambda(x,y,\theta,\varphi,p_x,p_y,p_\theta,p_\varphi)
=(\lambda_1, \lambda_2, \lambda_3, \lambda_4, \lambda_5, \lambda_6, \lambda_7, \lambda_8)\\
& = \gamma \cdot \pi_{Q}(x,y,\theta,\varphi,p_x,p_y,p_\theta,p_\varphi)= \gamma(x,y,\theta,\varphi)\\
& =(\gamma_1 \cdot \pi_{Q}, \gamma_2 \cdot \pi_{Q}, \gamma_3 \cdot \pi_{Q},
 \gamma_4 \cdot \pi_{Q}, \gamma_5 \cdot \pi_{Q}, \gamma_6 \cdot \pi_{Q},
\gamma_7 \cdot \pi_{Q}, \gamma_8 \cdot \pi_{Q}),
\end{align*}
that is, $\lambda_i= \gamma_i \cdot \pi_Q, \; i=1, \cdots, 8,$
where $\lambda_i, \; i=1,\cdots,8,$ are functions on
$T^* Q$, and $\gamma_i, \; i=1,\cdots,8,$ are functions on $Q$.
We may choose $q=(x,y,\theta,\varphi)\in Q, $ such that
$\gamma_1(q)=x, \; \gamma_2(q)=y, \; \gamma_3(q)=\theta, \; \gamma_4(q)=\varphi, $ and
$\gamma(q)= \gamma_5(q)\mathbf{d}x+ \gamma_6(q)\mathbf{d}y
+ \gamma_7(q)\mathbf{d}\theta+ \gamma_8(q)\mathbf{d}\varphi. $
Note that $\mathcal{D}(q)=
\textrm{Span}\{R\cos\varphi\partial_x+ R\sin\varphi\partial_y+
\partial_{\theta}, \; \partial_{\varphi}\}, $
take that $\alpha= R\cos\varphi\partial_x+ R\sin\varphi\partial_y+
\partial_{\theta} $ and $\beta= \partial_{\varphi}, $ then we have that
\begin{align*}
\mathbf{d}\gamma(\alpha,\beta)& = \alpha(\gamma(\beta))-\beta(\gamma(\alpha))-\gamma([\alpha,\beta])\\
& = R\cos\varphi(\frac{\partial \gamma_8}{\partial x}- \frac{\partial \gamma_5}{\partial \varphi})
+ R\sin\varphi(\frac{\partial \gamma_8}{\partial y}- \frac{\partial \gamma_6}{\partial \varphi})\\
& \;\;\;\;\; + (\frac{\partial \gamma_8}{\partial \theta}- \frac{\partial \gamma_7}{\partial \varphi})
+ 2(R\sin\varphi \gamma_5- R\cos\varphi \gamma_6).
\end{align*}
Thus, when $\mathbf{d}\gamma(\alpha,\beta)=0, $ we know that for any
$v, w \in \mathcal{F}, $ and $T\pi_{Q}(v), \; T\pi_{Q}(w) \in \mathcal{D},$  then
$\mathbf{d}\gamma(T\pi_{Q}(v),T\pi_{Q}(w))=0, $ that is,
$\gamma$ is closed on $\mathcal{D}$ with respect to $T\pi_{Q}:
TT^* Q \rightarrow TQ. $
Note that $Im(\gamma) \subset \mathcal{M}, $ then we have that
$p_x=\gamma_5,\; p_y=\gamma_6, \; p_\theta=\gamma_7.\; p_\varphi=\gamma_8, $ and
$\gamma_5= \frac{mR}{I}\gamma_7 \cos\varphi, \;\; \gamma_6=
\frac{mR}{I}\gamma_7 \sin\varphi. $
Hence,
\begin{align*} H\cdot \gamma =
\frac{1}{2m}(\frac{m^2R^2}{I^2}\gamma_7^2)+ \frac{1}{2I} \gamma_7^2
+ \frac{1}{2J} \gamma_8^2 = \frac{1}{2I}(1+
\frac{mR^2}{I})\gamma_7^2 +\frac{1}{2J} \gamma_8^2 ,
\end{align*}
and \begin{align*} X_H \cdot \gamma &= \frac{R}{I}\gamma_7 \cos\varphi \partial_x
+ \frac{R}{I}\gamma_7 \sin\varphi \partial_y +
\frac{1}{I}\gamma_7 \partial_\theta +
\frac{1}{J}\gamma_8 \partial_\varphi\\
&= \frac{\gamma_7}{I}(R\cos\varphi \partial_x
+ R \sin\varphi \partial_y + \partial_\theta )+ \frac{\gamma_8}{J}\partial_\varphi
= X_{\mathcal{K}} \cdot \gamma,
\end{align*}
$$ X_H^\gamma= T\pi_Q \cdot
X_H \cdot \gamma= \frac{\gamma_7}{I}(R\cos\varphi \partial_x
+ R \sin\varphi \partial_y + \partial_\theta )+ \frac{\gamma_8}{J}\partial_\varphi .$$
Thus, $T\gamma \cdot X_H^\gamma= X_{\mathcal{K}} \cdot \gamma, $
that is, the Type I of Hamilton-Jacobi equation for the
distributional Hamiltonian system
$(\mathcal{K},\omega_{\mathcal{K}},H )$ holds trivially.\\

Now, for any symplectic map $\varepsilon: T^* Q \rightarrow T^* Q, $ from
$\omega= \varepsilon^*\omega=\omega \cdot \varepsilon
= (\partial_x \varepsilon\cdot \partial_{p_x}\varepsilon) \mathbf{d}x \wedge \mathbf{d}p_x
+ (\partial_y \varepsilon\cdot \partial_{p_y}\varepsilon) \mathbf{d}y \wedge \mathbf{d}p_y
+ (\partial_\theta \varepsilon\cdot \partial_{p_\theta}\varepsilon) \mathbf{d}\theta \wedge \mathbf{d}p_\theta
+ (\partial_\varphi \varepsilon\cdot \partial_{p_\varphi}\varepsilon) \mathbf{d}\varphi \wedge \mathbf{d}p_\varphi, $ we have that
\begin{align*}
 \partial_x \varepsilon\cdot \partial_{p_x}\varepsilon= 1, \;\;\;\;\;
 \partial_y \varepsilon\cdot \partial_{p_y}\varepsilon= 1, \;\;\;\;\;
 \partial_\theta \varepsilon\cdot \partial_{p_\theta}\varepsilon= 1, \;\;\;\;\;
 \partial_\varphi \varepsilon\cdot \partial_{p_\varphi}\varepsilon= 1.
\end{align*}
Denote by $\varepsilon(x,y,\theta,\varphi,p_x,p_y,p_\theta,p_\varphi)
=(\varepsilon_1, \varepsilon_2, \varepsilon_3, \varepsilon_4, \varepsilon_5, \varepsilon_6, \varepsilon_7, \varepsilon_8 ),$
then we have that
$$H\cdot \varepsilon= \frac{1}{2m}(\varepsilon_5^2 + \varepsilon_6^2)+
\frac{1}{2I}\varepsilon_7^2 + \frac{1}{2J}\varepsilon_8^2, $$ and
 $$X_H \cdot \varepsilon=\frac{1}{m}
\varepsilon_5\partial_x + \frac{1}{m}\varepsilon_6\partial_y +
\frac{1}{I}\varepsilon_7\partial_\theta +
\frac{1}{J}\varepsilon_8\partial_\varphi ,$$
and hence $$ X_H^\varepsilon= T\pi_Q \cdot
X_H \cdot \varepsilon= \frac{1}{m}
\varepsilon_5\partial_x + \frac{1}{m}\varepsilon_6\partial_y +
\frac{1}{I}\varepsilon_7\partial_\theta +
\frac{1}{J}\varepsilon_8\partial_\varphi .$$
Since $Im(\gamma) \subset \mathcal{M}, $ then
\begin{align*}
T\gamma\cdot X_H^\varepsilon
& = T\gamma\cdot(\frac{1}{m}
\varepsilon_5\partial_x + \frac{1}{m}\varepsilon_6\partial_y +
\frac{1}{I}\varepsilon_7\partial_\theta +
\frac{1}{J}\varepsilon_8\partial_\varphi)\\
&= \frac{R}{I}\varepsilon_7 \cos\varphi \partial_x
+ \frac{R}{I}\varepsilon_7 \sin\varphi \partial_y +
\frac{1}{I}\varepsilon_7 \partial_\theta +
\frac{1}{J}\varepsilon_8 \partial_\varphi\\
& = \frac{\varepsilon_7}{I}(R\cos\varphi \partial_x
+ R \sin\varphi \partial_y + \partial_\theta )+ \frac{\varepsilon_8}{J}\partial_\varphi\\
& = X_{\mathcal{K}}\cdot\varepsilon,
\end{align*}
because $\varepsilon_5= \frac{mR}{I}\varepsilon_7 \cos\varphi, \; \varepsilon_6=
\frac{mR}{I}\varepsilon_7 \sin\varphi. $ In the same way, note that
$\lambda=\gamma\cdot \pi_Q, $ and $Im(\lambda) \subset \mathcal{M}, $ then
\begin{align*}
T\lambda\cdot X_H \cdot \varepsilon &= \frac{R}{I}\varepsilon_7 \cos\varphi \partial_x
+ \frac{R}{I}\varepsilon_7 \sin\varphi \partial_y +
\frac{1}{I}\varepsilon_7 \partial_\theta +
\frac{1}{J}\varepsilon_8 \partial_\varphi\\
& = \frac{\varepsilon_7}{I}(R\cos\varphi \partial_x
+ R \sin\varphi \partial_y + \partial_\theta )+ \frac{\varepsilon_8}{J}\partial_\varphi
= X_{\mathcal{K}}\cdot\varepsilon.
\end{align*}
On the other hand, since $\varepsilon: T^* Q \rightarrow T^* Q $ is symplectic, we have that
\begin{align*}
\tau_{\mathcal{K}}\cdot T\varepsilon\cdot X_{H \cdot \varepsilon} & =\tau_{\mathcal{K}}\cdot X_H \cdot \varepsilon\\
& = \frac{R}{I}\varepsilon_7 \cos\varphi \partial_x
+ \frac{R}{I}\varepsilon_7 \sin\varphi \partial_y +
\frac{1}{I}\varepsilon_7 \partial_\theta +
\frac{1}{J}\varepsilon_8 \partial_\varphi\\
& = \frac{\varepsilon_7}{I}(R\cos\varphi \partial_x
+ R \sin\varphi \partial_y + \partial_\theta )+ \frac{\varepsilon_8}{J}\partial_\varphi
= X_{\mathcal{K}}\cdot\varepsilon.
\end{align*}
Thus, $T\gamma\cdot X_H^\varepsilon=X_{\mathcal{K}}\cdot\varepsilon
=T\lambda\cdot X_H \cdot \varepsilon=\tau_{\mathcal{K}}\cdot T\varepsilon\cdot X_{H \cdot \varepsilon}.$
In this case, we must have that $\varepsilon$ is a
solution of the Type II of Hamilton-Jacobi equation $T\gamma \cdot X_H^\varepsilon= X_{\mathcal{K}}
\cdot \varepsilon, $ for the distributional Hamiltonian system
$(\mathcal{K},\omega_{\mathcal{K}},H ), $ if and only if it is a solution of the equation
$T\lambda\cdot X_H \cdot \varepsilon
=\tau_{\mathcal{K}}\cdot T\varepsilon\cdot X_{H \cdot \varepsilon}$.\\

In the following we consider respectively the actions of two Lie
groups $G= \mathbb{R}^2$ and $G= SE(2)$ on $Q$, and derive precisely the motion equations
and the geometric constraint conditions of the reduced distributional two-forms
for the nonholonomic reduced dynamical vector fields, that is,
the Type I and Type II of Hamilton-Jacobi
equations of the nonholonomic reduced distributional Hamiltonian
systems. Firstly, we consider the action of Lie group $G=
\mathbb{R}^2$ on $Q$, which is given by $$\Phi: G\times Q
\rightarrow Q, \; \Phi((r,s),(x,y,\theta,\varphi))=
(x+r,y+s,\theta,\varphi), $$ and we have the cotangent lifted
$G$-action on $T^* Q$, and the Hamiltonian $H: T^* Q \rightarrow
\mathbb{R}$ is $G$-invariant. In this case we have that
$$\bar{\mathcal{M}}= \{(\theta,\varphi,p_x,p_y,p_\theta,p_\varphi)\in
T^* Q /G |\; p_x= \frac{mR}{I}p_\theta \cos\varphi, \; p_y=
\frac{mR}{I}p_\theta \sin\varphi \}, $$ and the reduced distribution is given by
$$\bar{\mathcal{K}}= \mathrm{span}\{\partial_\theta,
\partial_\varphi, \partial_{p_\theta},
\partial_{p_\varphi} \}, $$ and the non-degenerate and the reduced distributional two-form
$\omega_{\bar{\mathcal{K}}}$ is given by
$$\omega_{\bar{\mathcal{K}}}= (1+ \frac{mR^2}{I})\mathbf{d}\theta
\wedge \mathbf{d}p_\theta + \mathbf{d}\varphi \wedge
\mathbf{d}p_\varphi .$$ A direct computation yields
$$ \mathbf{i}_{\partial_\theta}\omega_{\bar{\mathcal{K}}}=(1+ \frac{mR^2}{I})
\mathbf{d}p_\theta, \;\;\;\;\;\;
\mathbf{i}_{\partial_\varphi}\omega_{\bar{\mathcal{K}}}=
\mathbf{d}p_\varphi, \;\;\;\;\;\;
\mathbf{i}_{\partial_{p_\theta}}\omega_{\bar{\mathcal{K}}}=-(1+
\frac{mR^2}{I}) \mathbf{d}\theta, \;\;\;\;\;\;
\mathbf{i}_{\partial_{p_\varphi}}\omega_{\bar{\mathcal{K}}}=-
\mathbf{d}\varphi,
$$
and $$\mathbf{d}h_{\bar{\mathcal{K}}}=\mathbf{d}H_\mathcal{K} =
\frac{1}{I}(1+ \frac{mR^2}{I})p_\theta\mathbf{d}p_\theta +
\frac{1}{J}p_\varphi \mathbf{d}p_\varphi .$$ Assume that
$X_{\bar{\mathcal{K}}}= X_1\partial_\theta + X_2\partial_\varphi +
X_3
\partial_{p_\theta} + X_4 \partial_{p_\varphi}, $ then
$$\mathbf{i}_{X_{\bar{\mathcal{K}}}}\omega_{\bar{\mathcal{K}}}= X_1((1+ \frac{mR^2}{I})
\mathbf{d}p_\theta ) + X_2\mathbf{d}p_\varphi - X_3 (1+
\frac{mR^2}{I})\mathbf{d}\theta - X_4\mathbf{d}\varphi$$
$$ = - X_3 (1+
\frac{mR^2}{I})\mathbf{d}\theta - X_4\mathbf{d}\varphi + X_1(1+
\frac{mR^2}{I})\mathbf{d}p_\theta + X_2\mathbf{d}p_\varphi .$$ From
the nonholonomic reduced distributional Hamiltonian equation
$\mathbf{i}_{X_{\bar{\mathcal{K}}}}\omega_{\bar{\mathcal{K}}}=
\mathbf{d}h_{\bar{\mathcal{K}}}, $ we have that
$$X_1=\frac{1}{I}p_\theta, \;\;\;  X_2=\frac{1}{J}p_\varphi, \;\;\;
X_3= 0, \;\;\; X_4= 0.$$ Hence, the nonholonomic reduced dynamical vector
field is given by $$X_{\bar{\mathcal{K}}}= \frac{1}{I}p_\theta\partial_\theta +
\frac{1}{J}p_\varphi\partial_{\varphi},$$ and the motion equations of
the nonholonomic reduced distributional Hamiltonian system
$(\bar{\mathcal{K}},\omega_{\bar{\mathcal{K}}},h )$ are expressed by
$$ \dot{\theta}=\frac{1}{I}p_\theta,\;\;\; \dot{\varphi}=
\frac{1}{J}p_\varphi, \;\;\; \dot{p}_{\theta}=0, \;\;\;
\dot{p}_{\varphi}=0.$$

In the following we shall derive precisely
the Type I and Type II of Hamilton-Jacobi equations for
the nonholonomic reduced distributional Hamiltonian system
$(\bar{\mathcal{K}},\omega_{\bar{\mathcal{K}}},h ). $ As above
$\gamma: Q \rightarrow T^* Q, $ and $\lambda= \gamma \cdot \pi_{Q}:
T^* Q \rightarrow T^* Q, $ and assume that $\textmd{Im}(\gamma)\subset \mathcal{M}, $ and it is
$G$-invariant, $ \textmd{Im}(T\gamma)\subset \mathcal{K}, $ then we have that
$\bar{\gamma}=\pi_{/G}(\gamma): Q \rightarrow T^* Q/G, \;
\bar{\gamma}(x,y,\theta,\varphi)=
(\bar{\gamma}_1,\bar{\gamma}_2,\bar{\gamma}_3,\bar{\gamma}_4,\bar{\gamma}_5,\bar{\gamma}_6), $
and $\bar{\lambda}=\pi_{/G}(\lambda): T^* Q \rightarrow T^* Q/G, \;
\bar{\lambda}(x,y,\theta,\varphi,p_x,p_y,p_\theta,p_\varphi)=
(\bar{\lambda}_1,\bar{\lambda}_2,\bar{\lambda}_3,\bar{\lambda}_4,\bar{\lambda}_5,\bar{\lambda}_6),
$ that is, $\bar{\lambda}_i= \bar{\gamma}_i \cdot \pi_Q, \; i=1,\cdots,6, $
where $\bar{\lambda}_i, \; i=1,\cdots,6, $ are functions on $T^*Q$,
and $\bar{\gamma}_j, \; i=1,\cdots,6, $ are functions on $Q$.
Note that $h\cdot \pi_{/G}= \tau_{\mathcal{M}}\cdot H, $ and the reduced distribution
$\bar{\mathcal{K}}= \mathrm{span}\{\partial_\theta, \partial_\varphi, \partial_{p_\theta},
\partial_{p_\varphi} \}, $ then we have that
$$ h\cdot \bar{\gamma}= \frac{1}{2I}\bar{\gamma}_5^2 +\frac{1}{2J} \bar{\gamma}_6^2, \;\;\;\;\;\;\;\;\;\;
X_h\cdot \bar{\gamma}=\frac{1}{I}\bar{\gamma}_5 \partial_\theta +
\frac{1}{J}\bar{\gamma}_6 \partial_\varphi.
$$
When $\mathbf{d}\gamma(\alpha,\beta)=0, $ that is,
$\gamma$ is closed on $\mathcal{D}$ with respect to $T\pi_{Q}:
TT^* Q \rightarrow TQ, $ we have that
$$
T\bar{\gamma}\cdot X^\gamma_H= \tau_{\bar{\mathcal{K}}}\cdot X_h\cdot \bar{\gamma}
=X_{\bar{\mathcal{K}}}\cdot \bar{\gamma},
$$
that is, the Type I of Hamilton-Jacobi equation for
the nonholonomic reduced distributional Hamiltonian system
$(\bar{\mathcal{K}},\omega_{\bar{\mathcal{K}}},h ) $ holds.\\

Now, for any $G$-invariant symplectic map $\varepsilon: T^* Q\rightarrow T^* Q, $
$\bar{\varepsilon}=\pi_{/G}(\varepsilon): T^* Q \rightarrow T^*
Q/G,$ is given by $\bar{\varepsilon}(x,y,\theta,\varphi,p_x,p_y,p_\theta,p_\varphi)=
(\bar{\varepsilon}_1,\bar{\varepsilon}_2,\bar{\varepsilon}_3,\bar{\varepsilon}_4,\bar{\varepsilon}_5,\bar{\varepsilon}_6), $ then we have that
$$ h\cdot \bar{\varepsilon}= \frac{1}{2I}\bar{\varepsilon}_5^2 +\frac{1}{2J} \bar{\varepsilon}_6^2, \;\;\;\;\;\;\;\;\;\;
X_h\cdot \bar{\varepsilon}=\frac{1}{I}\bar{\varepsilon}_5 \partial_\theta +
\frac{1}{J}\bar{\varepsilon}_6 \partial_\varphi.
$$
Since $\textmd{Im}(\gamma)\subset \mathcal{M}, $ and $ \textmd{Im}(T\gamma)\subset \mathcal{K}, $
and hence $Im(\bar{\gamma}) \subset \bar{\mathcal{M}} $, $Im(T\bar{\gamma})\subset \bar{\mathcal{K}}. $
Thus,
$$
T\bar{\gamma}\cdot X_H^\varepsilon=\tau_{\bar{\mathcal{K}}}\cdot X_h\cdot \bar{\varepsilon}
=\frac{1}{I}\bar{\varepsilon}_5 \partial_\theta +
\frac{1}{J}\bar{\varepsilon}_6 \partial_\varphi
=X_{\bar{\mathcal{K}}}\cdot \bar{\varepsilon}.
$$
Note that $\bar{\lambda}= \bar{\gamma} \cdot \pi_Q, $ and $Im(\bar{\lambda}) \subset \bar{\mathcal{M}}, $
and $Im(T\bar{\lambda})\subset \bar{\mathcal{K}}, $
then we have that
$$
T\bar{\lambda}\cdot X_H \cdot \varepsilon=\tau_{\bar{\mathcal{K}}}\cdot X_h\cdot \bar{\varepsilon}
=X_{\bar{\mathcal{K}}}\cdot \bar{\varepsilon}.
$$
On the other hand, since $\varepsilon: T^* Q \rightarrow T^* Q $ is symplectic, and
$\bar{\varepsilon}^*= \varepsilon^*\cdot \pi_{/G}^*: T^*(T^* Q)/G
\rightarrow T^*T^* Q$ is also symplectic along $\bar{\varepsilon}$, then we have that
\begin{align*}
\tau_{\bar{\mathcal{K}}}\cdot T\bar{\varepsilon}\cdot X_{h \cdot \bar{\varepsilon}}
& = \tau_{\bar{\mathcal{K}}}\cdot X_h\cdot \bar{\varepsilon}\\
& = \frac{1}{I}\bar{\varepsilon}_5 \partial_\theta +
\frac{1}{J}\bar{\varepsilon}_6 \partial_\varphi= X_{\bar{\mathcal{K}}}\cdot \bar{\varepsilon}.
\end{align*}
Thus, $T\bar{\gamma}\cdot X_H^\varepsilon=X_{\bar{\mathcal{K}}}\cdot\bar{\varepsilon}
=T\bar{\lambda}\cdot X_H \cdot \varepsilon=\tau_{\bar{\mathcal{K}}}\cdot T\bar{\varepsilon}\cdot X_{h \cdot \bar{\varepsilon}}.$
In this case, we must have that
$\varepsilon$ and $\bar{\varepsilon}$ are the
solution of the Type II of Hamilton-Jacobi equation $T\bar{\gamma} \cdot X_H^\varepsilon= X_{\bar{\mathcal{K}}}
\cdot \bar{\varepsilon}, $ for the nonholonomic reduced distributional Hamiltonian system
$(\bar{\mathcal{K}},\omega_{\bar{\mathcal{K}}},h ), $ if and only if they satisfy the equation
$T\bar{\lambda}\cdot X_H \cdot \varepsilon=\tau_{\bar{\mathcal{K}}}\cdot T\bar{\varepsilon}\cdot X_{h \cdot \bar{\varepsilon}}$.\\

Next, we consider the action of the Lie group $G= SE(2)\cong
SO(2)\circledS \mathbb{R}^2$ on $Q$, which is given by $$\Phi:
G\times Q \rightarrow Q, \; \Phi((\alpha,r,s),(x,y,\theta,\varphi))=
(x\cos\alpha- y\sin\alpha+ r,x\sin\alpha +y\cos\alpha
+s,\theta,\varphi+\alpha), $$ and we have the cotangent lifted
$G$-action on $T^* Q$; and the Hamiltonian $H: T^* Q \rightarrow
\mathbb{R}$ is $G$-invariant. In this case we have that
$$\bar{\mathcal{M}}= \{(\theta,p_x,p_y,p_\theta)\in T^* Q /G |\;
p_x= \frac{mR}{I}p_\theta \cos\varphi, \; p_y= \frac{mR}{I}p_\theta
\sin\varphi \}, $$ and the reduced distribution is given by $\bar{\mathcal{K}}=
\mathrm{span}\{\partial_\theta,
\partial_{p_\theta} \}, $ and the non-degenerate and the reduced distributional two-form
$\omega_{\bar{\mathcal{K}}}$ is given by
$$\omega_{\bar{\mathcal{K}}}= (1+ \frac{mR^2}{I})\mathbf{d}\theta
\wedge \mathbf{d}p_\theta .$$ A direct computation yields
$$ \mathbf{i}_{\partial_\theta}\omega_{\bar{\mathcal{K}}}=(1+ \frac{mR^2}{I})
\mathbf{d}p_\theta, \;\;\;\;\;\;\;\;\;
\mathbf{i}_{\partial_{p_\theta}}\omega_{\bar{\mathcal{K}}}= -(1+
\frac{mR^2}{I}) \mathbf{d}\theta,
$$
and $$\mathbf{d}h_{\bar{\mathcal{K}}}=\mathbf{d}H_\mathcal{K} =
\frac{1}{I}(1+ \frac{mR^2}{I})p_\theta\mathbf{d}p_\theta +
\frac{1}{J}p_\varphi \mathbf{d}p_\varphi .$$ Assume that
$X_{\bar{\mathcal{K}}}= X_1\partial_\theta + X_2
\partial_{p_\theta}, $ then we have that
$$\mathbf{i}_{X_{\bar{\mathcal{K}}}}\omega_{\bar{\mathcal{K}}}= X_1((1+ \frac{mR^2}{I})
\mathbf{d}p_\theta ) + X_2(-(1+ \frac{mR^2}{I}) \mathbf{d}\theta)
.$$ From the nonholonomic reduced distributional Hamiltonian equation
$\mathbf{i}_{X_{\bar{\mathcal{K}}}}\omega_{\bar{\mathcal{K}}}=
\mathbf{d}h_{\bar{\mathcal{K}}}, $ we have that
$X_1=\frac{1}{I}p_\theta, \;\; X_2=0.$ Hence, the nonholonomic
reduced dynamical vector field is given by $$X_{\bar{\mathcal{K}}}=
\frac{1}{I}p_\theta\partial_\theta ,$$ and the motion equations of the
nonholonomic reduced distributional Hamiltonian system
$(\bar{\mathcal{K}},\omega_{\bar{\mathcal{K}}},h )$ are expressed by
$$ \dot{\theta}=\frac{1}{I}p_\theta,\;\;\; \dot{p}_{\theta}=0.$$

In the following we shall derive precisely the Type I and Type II of Hamilton-Jacobi equations for the
nonholonomic reduced distributional Hamiltonian system
$(\bar{\mathcal{K}},\omega_{\bar{\mathcal{K}}},h ). $ As above
$\gamma: Q \rightarrow T^* Q, $ and $\lambda= \gamma \cdot \pi_{Q}:
T^* Q \rightarrow T^* Q, $ and assume that $\textmd{Im}(\gamma)\subset \mathcal{M}, $ and it is
$G$-invariant, $ \textmd{Im}(T\gamma)\subset \mathcal{K}, $ then we have that
$\bar{\gamma}=\pi_{/G}(\gamma): Q \rightarrow T^* Q/G, \;
\bar{\gamma}(x,y,\theta,\varphi)=
(\bar{\gamma}_1,\bar{\gamma}_2,\bar{\gamma}_3,\bar{\gamma}_4), $ and
$\bar{\lambda}=\pi_{/G}(\lambda): T^* Q \rightarrow T^* Q/G, \;
\bar{\lambda}(x,y,\theta,\varphi,p_x,p_y,p_\theta,p_\varphi)=
(\bar{\lambda}_1,\bar{\lambda}_2,\bar{\lambda}_3,\bar{\lambda}_4), $
that is, $\bar{\lambda}_i= \bar{\gamma}_i \cdot \pi_Q, \; i=1,\cdots,4, $
where $\bar{\lambda}_i, \; i=1,\cdots,4, $ are functions on $T^*Q$, and $\bar{\gamma}_i,
\; i=1,\cdots,4,$ are functions on $Q$.
Note that $h\cdot \pi_{/G}= \tau_{\mathcal{M}}\cdot H, $ and the reduced distribution is given by
$\bar{\mathcal{K}}= \mathrm{span}\{\partial_\theta, \partial_{p_\theta}\}, $ then we have that
$$ h\cdot \bar{\gamma}= \frac{1}{2I}\bar{\gamma}_3^2,  \;\;\;\;\;\;\;\;\;\;\;\;\;\;\;
X_h\cdot \bar{\gamma}=\frac{1}{I}\bar{\gamma}_3 \partial_\theta .
$$
When $\mathbf{d}\gamma(\alpha,\beta)=0, $ that is,
$\gamma$ is closed on $\mathcal{D}$ with respect to $T\pi_{Q}:
TT^* Q \rightarrow TQ, $ we have that
$$
T\bar{\gamma}\cdot X^\gamma_H= \tau_{\bar{\mathcal{K}}}\cdot X_h\cdot \bar{\gamma}
=X_{\bar{\mathcal{K}}}\cdot \bar{\gamma},
$$
that is, the Type I of Hamilton-Jacobi equation for
the nonholonomic reduced distributional Hamiltonian system
$(\bar{\mathcal{K}},\omega_{\bar{\mathcal{K}}},h ) $ holds.\\

Now, for any $G$-invariant symplectic map $\varepsilon: T^* Q\rightarrow T^* Q, $
$\bar{\varepsilon}=\pi_{/G}(\varepsilon): T^* Q \rightarrow T^*
Q/G,$ is given by $\bar{\varepsilon}(x,y,\theta,\varphi,p_x,p_y,p_\theta,p_\varphi)=
(\bar{\varepsilon}_1,\bar{\varepsilon}_2,\bar{\varepsilon}_3,\bar{\varepsilon}_4), $ then we have that
$$ h\cdot \bar{\varepsilon}= \frac{1}{2I}\bar{\varepsilon}_3^2, \;\;\;\;\;\;\;\;\;\;\;\;\;\;\;
X_h\cdot \bar{\varepsilon}=\frac{1}{I}\bar{\varepsilon}_3 \partial_\theta .
$$
Since $\textmd{Im}(\gamma)\subset \mathcal{M}, $ and $ \textmd{Im}(T\gamma)\subset \mathcal{K}, $
and hence $Im(\bar{\gamma}) \subset \bar{\mathcal{M}} $, $Im(T\bar{\gamma})\subset \bar{\mathcal{K}}. $
Thus,
$$
T\bar{\gamma}\cdot X_H^\varepsilon=\tau_{\bar{\mathcal{K}}}\cdot X_h\cdot \bar{\varepsilon}
=\frac{1}{I}\bar{\varepsilon}_3 \partial_\theta
=X_{\bar{\mathcal{K}}}\cdot \bar{\varepsilon}.
$$
Note that $\bar{\lambda}= \bar{\gamma} \cdot \pi_Q, $ and $Im(\bar{\lambda}) \subset \bar{\mathcal{M}}, $
and $Im(T\bar{\lambda})\subset \bar{\mathcal{K}}, $
then we have that
$$
T\bar{\lambda}\cdot X_H \cdot \varepsilon=\tau_{\bar{\mathcal{K}}}\cdot X_h\cdot \bar{\varepsilon}
=X_{\bar{\mathcal{K}}}\cdot \bar{\varepsilon}.
$$
On the other hand, since $\varepsilon: T^* Q \rightarrow T^* Q $ is symplectic, and
$\bar{\varepsilon}^*= \varepsilon^*\cdot \pi_{/G}^*: T^*(T^* Q)/G
\rightarrow T^*T^* Q$ is also symplectic along $\bar{\varepsilon}$, then we have that
\begin{align*}
\tau_{\bar{\mathcal{K}}}\cdot T\bar{\varepsilon}\cdot X_{h \cdot \bar{\varepsilon}}
= \tau_{\bar{\mathcal{K}}}\cdot X_h\cdot \bar{\varepsilon}
= \frac{1}{I}\bar{\varepsilon}_3 \partial_\theta = X_{\bar{\mathcal{K}}}\cdot \bar{\varepsilon}.
\end{align*}
Thus, $T\bar{\gamma}\cdot X_H^\varepsilon=X_{\bar{\mathcal{K}}}\cdot\bar{\varepsilon}
=T\bar{\lambda}\cdot X_H \cdot \varepsilon=\tau_{\bar{\mathcal{K}}}\cdot T\bar{\varepsilon}\cdot X_{h \cdot \bar{\varepsilon}}.$
In this case, we must have that
$\varepsilon$ and $\bar{\varepsilon}$ are the
solution of the Type II of Hamilton-Jacobi equation $T\bar{\gamma} \cdot X_H^\varepsilon= X_{\bar{\mathcal{K}}}
\cdot \bar{\varepsilon}, $ for the nonholonomic reduced distributional Hamiltonian system
$(\bar{\mathcal{K}},\omega_{\bar{\mathcal{K}}},h ), $ if and only if they satisfy the equation
$T\bar{\lambda}\cdot X_H \cdot \varepsilon=\tau_{\bar{\mathcal{K}}}\cdot T\bar{\varepsilon}\cdot X_{h \cdot \bar{\varepsilon}}$.\\

In this paper, we study the
Hamilton-Jacobi theory for
the nonholonomic Hamiltonian system and
the nonholonomic reducible Hamiltonian system
on a cotangent bundle, by using the distributional Hamiltonian system
and the reduced distributional Hamiltonian system. These researches
reveal from the geometrical point of view the internal relationships of
nonholonomic constraints, distributional two forms and
nonholonomic dynamical vector fields of a
mechanical system and its nonholonomic reduced systems.
It is well known that the theory of controlled mechanical systems
became an important subject in recent years. Its research
gathers together some separate areas of research such as mechanics,
differential geometry and nonlinear control theory, etc., and the
emphasis of this research on geometry is motivated by the aim of
understanding the structure of equations of motion of the system, in
a way that helps both for analysis and design. Thus, it is natural to
study controlled mechanical systems by combining with the analysis
of dynamic systems and the geometric reduction theory of Hamiltonian
and Lagrangian systems. In particular, Marsden et al. in \cite{mawazh10} set up the
regular reduction theory of a regular controlled Hamiltonian system on a symplectic fiber
bundle, by using momentum map and the associated reduced symplectic
form, and from the viewpoint of completeness of Marsden-Weinstein symplectic
reduction, and some developments around the above work are given in Wang and
Zhang \cite{wazh12}, Ratiu and Wang \cite{rawa12}, and Wang \cite{wa15a}.
Since the Hamilton-Jacobi theory
is developed based on the Hamiltonian picture of dynamics, it is
natural idea to extend the Hamilton-Jacobi theory to the (regular)
controlled Hamiltonian system and its a variety of reduced systems,
and it is also possible to describe the relationship between the
CH-equivalence for the controlled Hamiltonian systems and the solutions
of corresponding Hamilton-Jacobi equations, see Wang \cite{wa13d,
wa20a, wa13e} for more details.
In particular, it is the key thought of the researches of geometrical mechanics
of the Professor Jerrold E. Marsden to explore and reveal the deeply internal
relationship between the geometrical structure of phase space and the dynamical
vector field of a mechanical system. It is also our goal of pursuing and inheriting.\\

\end{document}